\documentclass{amsart}
\usepackage{amssymb}
\newtheorem{thm}{Theorem}[section]
\newtheorem{lemma}[thm]{Lemma}
\newtheorem{prop}[thm]{Proposition}
\newtheorem{cor}[thm]{Corollary}

\theoremstyle{definition}

\newtheorem{remark}[thm]{Remark}
\def \intx {\stackrel{\circ}{X}}
\def \defi {\stackrel{\operatorname{def}}{=}}

\def \restr {{\left.\right|}}
\def \mrn {{\mathbb R}^n}

\def \mr {{\mathbb R}}

\def \wtt {{\widetilde{T}}}

\def \mcs {{\mathcal S}}
\def \mcp {{\mathcal P}}
\def \mcm {{\mathcal M}}
\def \mca {{\mathcal A}}

\def \mcf {{\mathcal F}}
\def \mch {{\mathcal H}}

\def \mcr {{\mathcal R}}

\def \mcv {{\mathcal V}}
\def \mcq {{\mathcal Q}}

\def \mc {{\mathbb C}}
\def \mh {{\mathbb H}}
\def \mn {{\mathbb N}}

\def \eps {\epsilon}   
\def \la {\lambda}   
\def \La {\Lambda}   
\def \lan {\langle}   
\def \ran {\rangle}   
\def \del {\delta}   

\def \det {\operatorname{det}}
\def \tW {\widetilde{W}}

\def \tU {\widetilde{U}}

\newcommand{\ff}{f\!f}
\newcommand{\Id}{\operatorname{Id}}
\newcommand{\comp}{\operatorname{comp}}
\newcommand{\Res}{\operatorname{Res}}
\newcommand{\vol}{\operatorname{vol}}
\newcommand{\dvol}{d\operatorname{vol}}
\newcommand{\ac}{\operatorname{ac}}
\newcommand{\pp}{\operatorname{pp}}
\newcommand{\pan}{\operatorname{span}}
\newcommand{\loc}{\operatorname{loc}}

\def \ha{ {\frac{1}{2}}}
\def \tha{ {\frac{3}{2}}}
\def \oq {\frac{1}{4}}

\def \p {\partial}

\def \rao#1 {\frac{\p}{\p #1} #1}

\numberwithin{equation}{section}
\title[Radiation Fields]{Radiation Fields, Scattering 
 and Inverse Scattering on 
Asymptotically Hyperbolic Manifolds }
\author[S\'{a} Barreto]{Ant\^{o}nio S\'{a} Barreto}
\address{Department of Mathematics, Purdue University,
150 North University Street, West Lafayette IN  47907, USA}
\email{sabarre@math.purdue.edu}
\thanks{Supported by the NSF under grant DMS-0140657.}
\input epsf
\begin{document}
\maketitle

\section{Introduction }\label{intro}
 The purpose of this article is to define the radiation fields 
on asymptotically hyperbolic manifolds and to use them to
study scattering theory. The radiation fields on $\mrn$ and on
asymptotically Euclidean manifolds were introduced by F.G. Friedlander in 
a series of papers starting in the early 1960's 
\cite{fried0,fried01,fried02,fried1,fried2}. His program of using the 
radiation fields to obtain the scattering matrix in that general setting
was completed in 
\cite{sabradf}. Here we carry out the analogous construction on 
asymptotically hyperbolic manifolds.
After defining the radiation fields, we use them to give a {\it unitary}
translation 
representation of the wave group and to obtain the scattering matrix for 
such manifolds.  As an application, we use them to study the inverse 
problem of determining
the  manifold and the metric from the scattering matrix {\it at all energies.}

Asymptotically hyperbolic manifolds are smooth compact manifolds with 
boundary equipped with a 
complete 
metric that resembles the hyperbolic space near the boundary. 
The basic examples of such manifolds are  
the hyperbolic space and its quotients by certain discrete group actions, 
see \cite{mame}, but any $C^\infty$ compact manifold with boundary can be 
equipped with such a metric.

There is a history of interest in
scattering theory for this class of manifolds,
motivated by several problems of mathematics and physics, which goes back to
the work of Fadeev and Pavlov \cite{fad}, followed by Lax and Phillips
\cite{lptr,lpaut}, and later by several people, see for example
\cite{ag,bope,gz0,jsb,ma,mame,per} and references cited there.
More recently there has been interest in this class of manifolds in connection
to conformal field theory, see \cite{fefgr,graham,grazw} and references cited 
there.  

 Mazzeo,   Mazzeo and Melrose \cite{ma,ma1,mame} first studied the 
spectral and scattering theory
of the Laplacian in this general setting and gave a thorough description of
the resolvent and its meromorphic continuation.
 Their methods have been applied in \cite{bope,grazw,guil,gz0,jsb}
to study the scattering matrix   starting from a careful understanding of the 
structure of the solutions to the Schr\"odinger equation on a neighborhood of 
infinity.

 We will develop the scattering theory for this class of manifolds
using a dynamical approach in the style of Lax and Phillips
\cite{lp,lptr,lpaut},  but we do this
by following Friedlander \cite{fried0,fried01,fried02,fried1,fried2}.

 We define the forward radiation field for 
asymptotically hyperbolic manifolds as the 
limit, as times goes to infinity, of the forward  fundamental solution of the 
wave operator along certain light rays. The backward radiation field is
defined by reversing the time direction.
These are 
generalizations of the Lax-Phillips transform \cite{lax,lptr}
to this class of manifolds. We will 
show that this leads to a {\it unitary}
translation representation of the wave group and 
a dynamical definition of the scattering matrix as in \cite{lp}.

 In section \ref{mmc} we 
establish  the connection of the radiation fields and the
Poisson operator, which in this context is also 
called the Eisenstein function.  This is then used to
show that the stationary (via Schr\"odinger's equation) and the dynamical 
(via radiation fields) definitions of the scattering matrix are equivalent.

To show the existence of the radiation fields, we adapt
the techniques  of \cite{fried1,fried2} to this setting. 
To connect the radiation fields,
the Poisson operator and the scattering matrix, we use the construction 
of the resolvent of the  Laplacian due to Mazzeo and Melrose \cite{mame}
and the construction of the Eisenstein function from the resolvent from 
\cite{gz0,jsb}.  

In section \ref{ccurv}, as an example, we compute the
forward radiation field for the three dimensional hyperbolic space $\mh^3$ 
and show that it is given by the Lax-Phillips transform, which is based on the 
horocyclic Radon transform. The case  on  $\mh^n$ is treated in 
\cite{laxptr0}.

In section \ref{supth} we 
prove a precise support theorem -- in the terminology of Helgason 
\cite{helgason1,helgrt} -- for the radiation fields.  Theorem \ref{L5} below
 generalizes to this setting a theorem of
 Lax-Phillips, Theorem 3.13 of \cite{lptr}, see also \cite{lax}, 
obtained for the horocyclic Radon transform.
 Helgason \cite{helg5} proved this result for
compactly supported functions, but in more general symmetric spaces. 
Theorem \ref{L5} below extends this to asymptotically hyperbolic 
manifolds, where the horocyclic Radon transform is replaced by the radiation 
field. 
This can be thought  of
as a result in control theory which, roughly speaking,
says that the support of a function is controlled by the support of its
 radiation field.

 Radon type transforms are often used to study properties of
solutions of hyperbolic equation, but here we use the equation to study 
support properties of the radiation field.  This allows the use of uniqueness
theorems for partial differential equations to establish support properties
of these transforms. The main ingredients of 
the proof of the support theorem are
 H\"ormander's uniqueness theorem for
the Cauchy  problem, see Theorem 28.3.4 of\cite{Ho}, and two of its 
refinements, one due to Alinhac
\cite{Ali} and another one which is due to Tataru 
\cite{tat}.  
 The study of support properties of Radon transforms is a topic of 
interest in its own, see for example \cite{helgason1,lax} and references
cited there.

In section \ref{prometdet} we use 
the characterization of the scattering matrix through
the radiation fields and the
boundary control method of Belishev \cite{be}, 
see also \cite{bk1,kaku,kaku1}, and the book by  Katchalov, Kurylev and Lassas
\cite{lassas}, to study the inverse problem of
 determining the manifold and the metric from the scattering matrix 
{\it at all energies}.
We prove that the scattering matrix of an asymptotically hyperbolic manifold 
determines the manifold and the metric up to invariants.

\section{Asymptotically Hyperbolic Manifolds and Radiation
 Fields}\label{state}
A smooth compact manifold $X$ with 
boundary, $\p X,$ is called  asymptotically hyperbolic, see \cite{mame},
when it is equipped 
with a Riemannian metric $g,$ which is smooth in
the interior of $X,$ denoted by $\intx,$ and is such that
for  a smooth defining function $x$ of $\p X,$ that is $x>0$ in the interior 
of $X,$ $x=0$ on $\p X$ and $dx\not=0$ at $\p X,$
\begin{gather}
x^2 g=H, \label{hmet}
\end{gather}
is a smooth Riemannian metric on $X$ non-degenerate up to $\p X.$
Furthermore we assume that 
$$|dx|_{H}=1 \;\ \text{ at } \;\ \p X.$$
It can be shown, see \cite{ma,mame}, that under these assumptions, the 
sectional
curvature approaches $-1$ at $\p X.$ 

Observe that $g$ determines $x$ and $H$ only up to  a positive
factor.   Hence $g$ induces a conformal structure at $\p X.$

Throughout this paper, $X$ denotes a $n+1$ dimensional smooth compact
manifold with 
boundary,  and $n\geq 1.$ $g$ will be a Riemannian metric on $X$ satisfying 
\eqref{hmet} and $\Delta$ will denote the (positive) Laplace operator with 
respect to the metric $g.$ 

As  stated in \cite{megs}, 
see \cite{jsb} for a  proof in this general setting, fixed a defining function 
$x$ of $\p X,$
then for all $f \in C^\infty(\p X)$ and 
$\la \in \mr,$ $\la \not=0,$ there exists
{\it a unique} $u \in C^\infty(\stackrel{\circ}{X})$ satisfying
\begin{gather}
\begin{gathered}
(\Delta-\la^2-\frac{n^2}{4})u=0 \text{ in } \stackrel{\circ}{X}, \\
u=x^{i\la+\frac{n}{2}} f_+ + x^{-i\la+\frac{n}{2}} f_-, \;\ f_{\pm} \in 
C^\infty(\overline{X}), \;\
f_+|_{\p X}=f.
\end{gathered}\label{scat}
\end{gather}

This leads to the {\it stationary}\; definition of the scattering matrix at 
energy $\lambda\not=0,$  see for example \cite{gz0,jsb,megs}, as the operator
\begin{gather}
\begin{gathered}
A(\la): C^\infty(\p X) \longrightarrow C^\infty(\p X) \\
f \longmapsto f_{-}|_{\p X}.
\end{gathered}\label{scatmat}
\end{gather}

As pointed out in \cite{megs}, the expansion \eqref{scat}
gives two parametrizations, corresponding to  $\pm \la,$ of the 
generalized eigenspace of $\Delta-\frac{n^2}{4}$ with eigenvalue $\la^2$ by 
distributions on $\p X.$ The scattering matrix is the operator that
intertwines them.

Notice that if $\Psi$ is a diffeomorphism of $\overline{X},$ fixing $x$ and
$dx$ on $\p X,$ then the scattering matrix will be invariant under
the pull back of the metric by $\Psi.$  Moreover,
this definition of the scattering matrix depends on the choice of
the function $x.$  It can be invariantly defined as acting on appropriate
bundles, see for example \cite{jsb}.

It is shown in \cite{jsb} and \cite{graham}, that if $g$ satisfies 
\eqref{hmet}, and fixed a representative $h_0$ of the conformal class of
$g$ at $\p X,$  there exists
$\eps>0$ and a {\it unique } product structure 
$X \sim [0,\eps)\times \p X$ in which
\begin{gather}
g=\frac{dx^2}{x^2} +\frac{h(x,y,dy)}{x^2}, \;\ h_0=h(0,y,dy). \label{hmet1}
\end{gather}

We will fix such a decomposition, and from now on $x \in C^\infty(X)$ will
be as in \eqref{hmet1}.  This is equivalent to  fixing a conformal 
representative of $x^2g|_{\p X}.$ We will also work with the $A(\la)$ defined
by \eqref{scatmat} where $x$ is given by \eqref{hmet1}.

Here, as in  \cite{fried1,fried2},  we will use the wave equation to define 
the radiation fields and arrive 
at an equivalent definition of the scattering matrix. 

 We will prove
\begin{thm}\label{rad} For $f_1, f_2 \in
C_0^\infty(\intx),$ compactly
supported in the interior of $X,$ let 
$u(t,z) \in C^\infty(\mr_+ \times \intx)$  satisfy
\begin{gather}
\begin{gathered}
\left(D_{t}^2-\Delta+\frac{n^2}{4}\right)u(t,z)=0,  \text{ on } \mr_+
\times \intx , \\ 
u(0,z)=f_1(z), \;\ D_tu(0,z)=f_2(z).
\end{gathered}\label{we2}
\end{gather} 
Let
$z=(x,y)\in (0,\eps)\times \p X$ be local coordinates near $\p X$
 in which \eqref{hmet1} 
hold. Then there exist
 $w_k \in C^\infty(\mr \times \p X),$ such that
\begin{gather*}
x^{-\frac{n}{2}}u(s-\log x, x, y) \sim 
\sum_{k=0}^\infty  w_k(s,y)x^k.
\end{gather*}
\end{thm}

Clearly $w_k(s,y)$ depends on the choice of 
$x,$ and we make no attempt
to define a bundle where it would be invariant. We refer the reader to
Lemma 2.2 of \cite{graham} for the relationship between two functions that
satisfy \eqref{hmet1} corresponding to two different conformal representatives.

 Theorem \ref{rad} defines a map
\begin{gather}
\begin{gathered}
\mcr_{+} : C_0^\infty(\intx)\times C_0^\infty(\intx)
\longrightarrow
\ C^\infty(\mr \times \p X) \\
\mcr_+(f)(s,y)=x^{-\frac{n}{2}}D_t u(s-\log x,x,y)|_{x=0}=D_sw_0(s,y),
\end{gathered}\label{mcr}
\end{gather}
which will be called the {\it forward radiation field.}

Similarly one can prove that if $u_-$ satisfies
 \eqref{we2} in $\mr_-\times \intx$ then
\begin{gather*}
\lim_{x\rightarrow 0}x^{-\frac{n}{2}}u_-(s+\log x,x,y)= w_{-}(s,y) 
\end{gather*} 
exists, and thus define the {\it backward radiation field}
\begin{gather}
\begin{gathered}
\mcr_{-} : C_0^\infty(\intx)\times C_0^\infty(\intx)
\longrightarrow
\ C^\infty(\mr \times \p X)  \\
\mcr_-(f)(s,y)=x^{-\frac{n}{2}}u_-(s+\log x,x,y)|_{x=0}=D_sw_-(s,y).
\end{gathered}\label{mcr-}
\end{gather}

Finally we remark that, since the Lorentzian metric associated to $g$ is
\begin{gather*}
\sigma=dt^2-\frac{dx^2}{x^2}-\frac{h(x,y,dy)}{x^2}=
d(t-\log x)d(t+\log x)-\frac{h(x,y,dy)}{x^2},
\end{gather*}
the surfaces
\begin{gather*}
\{t-\log x=C\}, \;\ \{t+\log x=C\}
\end{gather*}
are characteristic for the wave operator, and thus a point
$(t',z'),$  $z'=(x',y'),$ has a past domain of dependence,
$\Delta^-(t',z')$ satisfying
\begin{gather}
\Delta^-(t',z')\subset\{ (t,x,y): \; t-\log x \geq t'-\log x', \;\
t+\log x \leq t'+\log x'\}.\label{domdep}
\end{gather}

\section{The Three Dimensional Space with Constant Negative 
Curvature}\label{ccurv}
In this section, as an example, we study the particular case of the
 hyperbolic space
\begin{gather}
{\mh}^{3}=\{(x,y): y\in \mr^2, \; x \in \mr, \;\ x>0\}, 
\text{ with the metric } g=\frac{dx^2}{x^2}+\frac{|dy|^2}{x^2}. \label{hmetc}
\end{gather}
In this case  the radiation fields can be explicitly computed.
The  formul\ae\; obtained in  \cite{russians}, see also \cite{helgason},
can be used in the same  way to compute the radiation fields in ${\mh}^n.$ 
This is  done in \cite{laxptr0}.

For convenience, we will work in the non-compact model, which does not 
quite fit the framework of section \ref{intro}, but it is isometric to the 
compact model given by the interior of the ball with  Poincar\'e's metric
$X={\mathbb B}^{3}=\{ z \in \mr^{3}: |z|<1\}$ and 
$g=\frac{4|d z|^2}{(1-|z|^2)^2}.$

For $z \in {\mathbb H}^3,$ let $S(z,t)$ denote the set of points in
${\mh}^{3}$ whose geodesic distance to $z$ is $t.$ Let $A(t)$ be the
area of $S(z,t),$ which is independent of $z.$  Given 
$f\in C_0^\infty({\mh}^3),$
supported in the interior of ${\mh}^3,$ let
\begin{gather*}
M(f,t,z)=\frac{1}{A(t)}\int_{S(z,t)} f(y) d\sigma_y
\end{gather*}
be the mean of $f$ over the sphere $S(z,t).$
 
According to \cite{helgason,russians},
the solution to \eqref{we2} with $f_1=0$ and 
$f_2=f$ is given by
\begin{gather}
\begin{gathered}
u(t,z)= (\sinh t) M(f,t,z).
\end{gathered}\label{fund}
\end{gather}
Therefore, the forward fundamental solution is, for $t>0,$
\begin{gather*}
U(t,z,z')=\frac{\sinh t}{A(t)}\delta\left(t-d(z',z)\right),
\end{gather*}
where the distance is given by
\begin{gather}
\cosh d(z,z')= \frac{x^2+{x'}^2+|y-y'|^2}{2x x'}. \label{cosh}
\end{gather}

Then the sphere in the hyperbolic metric
centered at $(x,y)$ with radius $t$ corresponds to the Euclidean sphere with
center $(x\cosh t, y)$ and radius $x\sinh t.$  We find that
$A(t)=4\pi \sinh^2 t,$ and hence
$$U(t,z,z')= \frac{1}{4\pi \sinh t}
\delta(t-d(z,z')).$$
But $\lim_{x\rightarrow 0} x\sinh(s-\log {x})=
\lim_{x\rightarrow 0} x\cosh(s-\log {x})=\frac{e^s}{2}.$
So the  Schwartz kernel of $\mcr_+$ is
\begin{gather}
\mcr_{+}(s,y,z')=\lim_{x\rightarrow 0} 
x^{-1}\frac{\p }{\p s}U(s-\log {x}, z,z')= 
\frac{\p }{\p s}\left(\frac{1}{4\pi e^s}
\delta\left( \frac{e^s}{2}-d_e\left[ \left(\frac{e^s}{2},y\right), z'\right]
\right)\right) \label{lptransf}
\end{gather}
where $d_e$ denotes the distance in the Euclidean metric.

So, for $f\in C_0^\infty(\mh^3),$ one obtains $\mcr_{+}f(s,y)$ 
by integrating $f$ over the surface of the Euclidean
sphere with center $(\frac{e^s}{2},y)$ and radius
$\frac{e^s}{2},$ with respect to the measure induced by
the metric $g.$ This sphere is tangent to $\{x=0\}$ at $(0,y)$ and
the integration of $f$ on those spheres is known as the horocyclic Radon 
transform, see for example \cite{helgrt}. The transformation given by
\eqref{lptransf} is called the Lax-Phillips transform, see 
\cite{lax}.

\section{The Proof of Theorem \ref{rad}}\label{prad}

In this section we work with the forward radiation field and
we will drop the index $\pm$ from the notation.
 We will also denote
\begin{gather}
h(x,y,dy)=\sum_{i,j=1}^n h_{ij}(x,y) dy_i dy_j, \;\ 
|h|=\det\left(h_{ij}\right),\text{ and }
h^{-1}=\left(h^{ij}\right), \text{ the inverse of the matrix } h_{ij}.
\label{invmat}
\end{gather}
\begin{proof}

Without loss of generality, we may assume  that
$f_1=0$ and that $f_2$ is compactly supported in $\{x>x_0\},$
 with $x_0$ small enough.
Hence $u(t,x,y) \in C^\infty(\mr \times (0,\eps) \times \p X)$
 satisfies 
\begin{gather}
\begin{gathered}
\left(D_t^2-\Delta+\frac{n^2}{4}\right)u(t,x,y)=0, \\
u(0,x,y)=0, \;\ D_tu(0,x,y)=f(x,y). \label{werad}
\end{gathered}
\end{gather}
Here we are considering the solution in $t\in \mr,$ and in this case the 
solution $u$ is odd in $t.$ Moreover \eqref{domdep}  
and the finite speed of propagation guarantee that for $x$ small, 
\begin{gather}
u=0 \; \text{ if } \;\ \log x-\log x_0<t<\log x_0-\log x, 
\;\  x<x_0<\eps .\label{werad1}
\end{gather}

We will show that
\begin{gather}
v(s,x,y)=x^{-\frac{n}{2}} u(s-\log x,x,y) \label{defv}
\end{gather}
is smooth up to $x=0.$ We denote
\begin{gather*}
P=x^{-\frac{n}{2}-1}\left(D_t^2-\Delta+\frac{n^2}{4}\right)x^{\frac{n}{2}}
\end{gather*}
and substitute $t=s-\log x.$ Then using \eqref{werad}, we find
\begin{gather}
\begin{gathered}
P v(s,x,y)=0, \text{ in } \mr \times (0,\eps)\times \p X, \;\
P=\frac{\p }{\p x }\left(2\frac{\p }{\p s }+x\frac{\p }{\p x}\right)+
x\Delta_h+ A\frac{\p }{\p s } + Ax \frac{\p }{\p x } +\frac{n}{2}A \\
v(\log x, x, y)=0, \;\ D_s v(\log x, x,y)=x^{-\frac{n}{2}}f(x,y)
\end{gathered}\label{neqq}
\end{gather}

Here $\Delta_h$ is the Laplacian on $\p X$ associated with the metric 
$h(x,y,dy),$ 
which in local coordinates is $\Delta_h=\frac{1}{|h|^\ha}\sum_{i,j=1}^{n} 
\frac{\p }{\p y_{i}}\left(|h|^{\ha} h^{ij} \frac{\p }{\p y_{j}}\right),$
and $A=\ha |h|^{-1}\frac{\p |h|}{\p x}.$ 

Equation \eqref{werad1} implies
\begin{gather}
v=0 \text{  for } \;\ 2\log x-\log x_0<s<\log x_0, \;\ x<x_0<\eps.
\label{supportv}
\end{gather}

The operator $P$ is not strictly hyperbolic at $x=0,$ so the argument
of Friedlander \cite{fried2}, which is based on a theorem of Leray,
 cannot be used directly.  We will refine the method of \cite{fried1}
and obtain energy estimates which hold uniformly up to
$\{x=0\}.$  

It is worth observing that if the tensor $h(x,y,dy)$ is an even function of
$x,$ then after setting $r=x^2,$ the
operator $r^{-\ha}P$ is smooth and strictly hyperbolic. So Friedlander's
method can be applied directly. Such metrics have been recently studied 
by Guillarmou in \cite{guil} and include the case where the metric has 
constant
curvature near $\p X,$ see \cite{gz1}. In the general case
this trick does not work because the resulting operator would not be smooth
at $r=0.$

In section \ref{supth}
we will need to understand the behavior of the forward radiation field
as $s\rightarrow -\infty,$ so we will compactify the problem and obtain
uniform estimates as $s\rightarrow -\infty.$  So we
 make the change of variables
\begin{gather}
s=2\log t', \;\ x=x't'. \label{x't'}
\end{gather}
This choice of coordinates is designed to do two things: first it transforms 
the
operator $\frac{\p}{\p x}\left(2\frac{\p}{\p s}+x\frac{\p}{\p x}\right)$
into $\frac{\p}{\p x'}\frac{\p}{\p t'},$ and secondly it compactifies
the half-line $(-\infty,0].$
 Thus let
$$V(x',t',y)=v(2\log t',x't',y)=
(x't')^{-\frac{n}{2}}u\left(\log\left(\frac{t'}{x'}\right),
x't',y\right).$$
Then $V$ is smooth in 
$t'>0,$ $x'>0$ and, as $u$ is odd in $t,$
$V(x',t',y)=-V(t',x',y).$
Moreover
for $$P'=\frac{\p}{\p x'}\frac{\p}{\p t'}+x't'\Delta_h+ 
\ha A(x't',y)\left(t'\frac{\p}{\p t'}+x'\frac{\p}{\p x'}\right)
+\frac{n}{2}A(x't',y),$$
\begin{gather}
\begin{gathered}
P'V=0, \;\ x'>0, \;\ t'>0, \\
V(x',x',y)=0, \;\ \frac{\p}{\p x'}V(x',x',y)=-{x'}^{-n-1}f({x'}^2,y).
\end{gathered}\label{eqq'}
\end{gather}
Here $\Delta_h$ is the Laplacian with respect to the metric $h(x't',y,dy)$
which in local coordinates is  $\Delta_h=
\frac{1}{|h(x't',y)|^{\ha}}\sum_{i,j}\frac{\p}{\p y_i}\left(|h(x't',y)|^{\ha}
h^{ij}(x't',y)\frac{\p}{\p y_j}\right).$

The support properties of $v$ given in \eqref{supportv} translate into
\begin{gather}
V(x',t',y)=0 \;\ \text{ if } \;\ x'< \sqrt{x_0}, \;\ t'< \sqrt{x_0}.
\label{supportV}
\end{gather}

The coefficients of $P'$ are smooth up to $\{x'=0\}\cup \{t'=0\}$ and therefore
can be extended smoothly, although not uniquely,
to a neighborhood $\{|x'|< \sqrt{x_0}\}\cup
\{|t'|< \sqrt{x_0}\}.$
To obtain the desired regularity of $V$ we will differentiate the equation 
\eqref{eqq'}  and obtain energy estimates for the resulting
system of differential equations. We begin by proving
\begin{lemma}\label{enerest} For $T>0,$ let 
$\Omega=[0,T]\times [0,T]\times \p X.$ 
Let $V(x',t',y)$ be a $N\times 1$ vector, 
smooth in $x'>0,$ $t'>0,$ satisfying the system
\begin{gather}
\begin{gathered}
\left(\frac{\p}{\p x'}\frac{\p}{\p t'}+x't'\Delta_h\right )V+ 
x't' B\left(x',t',y, \frac{\p}{\p y}\right) V+
 C(x',t',y)t'\frac{\p V}{\p t'}+ 
 D(x',t',y) x'\frac{\p V}{\p x'}+\\  E(x',t',y)V=0 
\text{ in } \;\ x'>0, \;\ t'>0 \;\ \text{ and } \\
V(x',x',y)=f_1(x',y), \;\ \frac{\p V}{\p x'}(x',x',y)=f_2(x',y), \;\ x'>0,
\end{gathered}\label{eqmcq}
\end{gather}
where $B\left(x',t',y,\frac{\p}{\p y}\right)$ is an $N\times N$ matrix of
first order differential operators having derivatives in $y$ only,
$C,$ $D$ and $E$ are $N\times N$ matrices of functions. Moreover
$B,C,D$ and $E$ are smooth in $|x'|<T,$  $|t'|<T,$
$y \in \p X.$ If the data  $f_1$ and $f_2$ are such that
\begin{gather}
\begin{gathered}
\int_{0}^T\int_{\p X}\left(  x'| f_1|^2 +
  {x'}^3|\nabla_y f_1|^2 \right)\sqrt{h}({x'}^2,y) dx' dy 
<\infty, \;\
\int_{0}^T\int_{\p X}
 x'|f_2|^2 \sqrt{h}({x'}^2,y) dx' dy <\infty, 
\end{gathered}\label{condD}
\end{gather}
where $\left|\nabla_y f_1\right|^2=\sum_j |\nabla_y f_{1,j}|^2,$ 
then for $T$ small
\begin{gather}
\begin{gathered}
\int_{ \Omega}\left[ |V|^2+ x't'(x'+t')\left|\nabla_y V\right|^2+
x'\left|\frac{\p V}{\p x'}\right|^2 +t'\left|\frac{\p V}{\p t'}\right|^2 
\right]  \sqrt{h}(x't',y) dy dx' dt'
\leq \\
C(T)\int_{0}^{T}
 \int_{\p X}\left( x'\left|f_1\right|^2 + x'\left|f_2\right|^2 
+{x'}^3\left|\nabla_y f_1\right|^2\right)(x',y)\sqrt{h}({x'}^2,y) dy dx',
\end{gathered}\label{est3T}
\end{gather}
where $\left|\nabla_y V\right|^2=\sum_j |\nabla_y V_j|^2.$ 
\end{lemma}
\begin{proof} We begin by multiplying the system \eqref{eqmcq} by
 $x'\frac{\p V}{\p x'}-t'\frac{\p V}{\p t'}.$   Notice that the
operator $x'\frac{\p }{\p x'}-t'\frac{\p }{\p t'}$ is $\frac{\p}{\p t}$
written in these coordinates.  We obtain
\begin{gather}
\begin{gathered}
0=\frac{1}{2\sqrt{h}(x't',y)}\frac{\p}{\p t'}\left[ \left( 
x'\left|\frac{\p V}{\p x'}\right|^2+ {t'}^2x'|\nabla_y V|^2\right)
\sqrt{h}(x't',y)\right]
-\\  \frac{1}{2\sqrt{h}(x't',y)}\frac{\p}{\p x'}\left[\left( 
t'\left|\frac{\p V}{\p t'}\right|^2+ t'{x'}^2|\nabla_y V|^2\right)
\sqrt{h}(x't',y)\right] +
\\
+\sum_k x't' \delta_h\left(
\left(x'\frac{\p V_k}{\p x'}-t'\frac{\p V_k}{\p t'}\right)d_hV_k\right) 
+Q\left(V,x'\frac{\p V}{\p x'},t'\frac{\p V}{\p t'},
x't'\frac{\p V}{\p y_j}\right).
\end{gathered}\label{enest1}
\end{gather}
Here $\delta_h$ is the divergence operator in the $y$ variable, dual to $d_h$
with respect to the metric $h(x't',y),$
and $Q$ is a quadratic form with smooth coefficients.

Let
\begin{gather*}
\Omega_{\del}=[\del,T]\times [\del,T]\times \p X, \;\
\Omega_{\del}^+=\{(x',t',y)\in \Omega_\del; \;\ t'\geq x'\}, \;\ \text{ and }\\
\Omega_{a,b}=\{(x',t',y)\in \Omega, \;\ a\leq x'\leq t'\leq b\},
\end{gather*}
see figure \ref{enes1'}.
\begin{figure}[int1]
\epsfxsize= 5.0in
\centerline{\epsffile{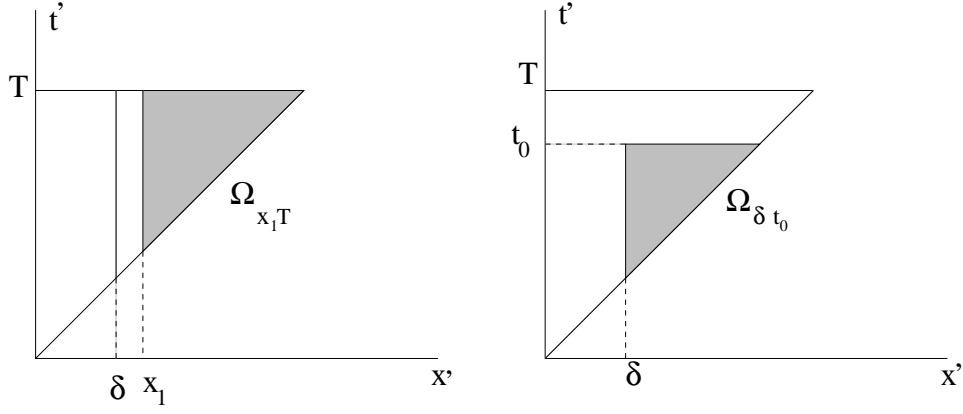}}
\caption{The regions $\Omega_{x_1,T}$ and $\Omega_{\del,t_0}.$}
\label{enes1'}
\end{figure}

Integrating \eqref{enest1} in 
$\Omega_{x_1,T},$ using the compactness of $\p X,$  the divergence 
theorem, and that the part of the first term in \eqref{enest1} which is inside 
the brackets is positive, we have
\begin{gather}
\begin{gathered}
\ha \int_{x_1}^{T} \int_{\p X}
\left( t'\left|\frac{\p V}{\p t'}\right|^2 +
t'{x'}^2|\nabla_y V|^2 \right)(x_1,t',y) 
\sqrt{h}(x_1t',y)dy dt'  
+\int_{\Omega_{x_1,T}} Q \;\ \sqrt{h}(x't',y) dy dx' dt' \\ \leq 
\ha \int_{x_1}^{T} \int_{\p X} \left( 
2x'\left|\frac{\p V}{\p x'}\right|^2
+ 2{x'}^3\left|\nabla_y V\right|^2\right)(x',x',y)\sqrt{h}({x'}^2,y) dy dx' .
\end{gathered}\label{ineq1}
\end{gather}

Proceeding similarly in the region $\Omega_{\delta, t_0}$ gives

\begin{gather}
\begin{gathered}
\ha \int_{\delta}^{t_0} \int_{\p X}
\left(x'\left|\frac{\p V}{\p x'}\right|^2 +
 x'{t'}^2|\nabla_y V|^2\right)(x',t_0,y) \sqrt{h}(x't_0,y) dy dt' + 
\int_{\Omega_{\delta,t_0}} Q \;\ \sqrt{h}(x't',y) dy dx' dt' \\ \leq 
\ha \int_{\del}^{t_0} \int_{\p X} \left( 
2x'\left|\frac{\p V}{\p x'}\right|^2
+ 2{x'}^3\left|\nabla_y V\right|^2\right)(x',x',y)\sqrt{h}({x'}^2,y) dy dx' .
\end{gathered}\label{ineq2}
\end{gather}

Next we integrate \eqref{ineq1} in $x_1\in [\delta,T]$ and
\eqref{ineq2} in $t_0\in [\delta,T]$ and add the results to get
\begin{gather}
\begin{gathered}
\ha \int_{\Omega_\del^+}
\left(  x't'(t'+x')|\nabla_y V|^2 +x'\left|\frac{\p V}{\p x'}\right|^2 
+t'\left|\frac{\p V}{\p t'}\right|^2 \right)
(x',t',y)\sqrt{h}(x't',y) dy dx' dt' + \\
\int_{\del}^{T} \int_{\Omega_{x_1,T}} Q \sqrt{h}(x't',y) dy dx' dt' dx_1+ 
\int_{\del}^{T} \int_{\Omega_{\delta,t_0}} Q
\sqrt{h}(x't',y) dy dx' dt' dt_0 \;\ \leq \\
 (T-\delta)\int_{\del}^{T} \int_{\p X} \left(
x'\left|\frac{\p V}{\p x'}\right|^2
+ {x'}^3\left|\nabla_y V\right|^2\right)(x',x',y)\sqrt{h}({x'}^2,y) dy dx' .
\end{gathered}\label{firsti} 
\end{gather}

With the exception of terms of $Q$ containing products involving $V_j,$
all the  others in the two middle terms of \eqref{firsti} are trivially
bounded by the terms in the first integral.  
To analyze the terms containing $V_j,$ we will bound the integral of $|V|^2$ in
$\Omega_\del^+.$
 So we write for $t'>x',$
\begin{gather*}
V_j(x',t',y)=V_j(t',t',y)-\int_{x'}^{t'} \frac{\p V_j}{\p s}(s,t',y) \; ds.
\end{gather*}
The  Cauchy-Schwartz inequality then gives
\begin{gather}
|V_j(x',t',y)|^2 \leq 2 |V_j(t',t',y)|^2 + 
2(t'-x') \int_{x'}^{t'} \left|\frac{\p V_j}{\p s}(s,t,y)\right|^2 \; ds. 
\label{term1}
\end{gather}
Hence, as $\sqrt{h}$ is bounded from above and below, there exists $C>0$ such 
that
\begin{gather*}
\int_{\Omega_\del^+}|V_j(x',t',y)|^2 \sqrt{h}(x't',y) \; dy dx' dt' =
\int_{\del}^{T}\int_{\del}^{t'} \int_{\p X} 
|V_j(x',t',y)|^2 \sqrt{h}(x't',y) \; dy dx' dt' \leq \\
C \left(\int_{\del}^{T} \int_{\p X} t' |V_j(t',t',y)|^2 \sqrt{h}({t'}^2,y) 
dy dt'+\int_{\del}^T \int_{\del}^{t'} \int_{\p X} t' 
\int_{x'}^{t'} \left|\frac{\p V_j}{\p s}(s,t',y)\right|^2 \sqrt{h}(st',y)
 ds \; dy dx' dt'\right).
\end{gather*}
Switching the order of integration of $y,$
 $x'$ and $s$ in the last integral, and
using that $t'<T,$  gives
\begin{gather}
\begin{gathered}
\int_{\Omega_\del^+} |V(x',t',y)|^2 \sqrt{h}(x't',y) \; dy dx' dt' \leq
C \int_{\del}^{T} \int_{\p X} t' |V_j(t',t',y)|^2 \sqrt{h}({t'}^2,y) \; dy dt'
+ \\ CT \int_{\Omega_\del^+}
s \left|\frac{\p V_j}{\p s}(s,t,y)\right|^2 \sqrt{h}(st',y) \; dy ds dt'.
\end{gathered}\label{estv2}
\end{gather}
Hence
\begin{gather*}
\int_{\del}^{T} \int_{\Omega_{x_1,T}} 
 \left|Q \sqrt{h}(x't',y)\right| dy dx' dt' dx_1+ 
\int_{\del}^{T} \int_{\Omega_{\delta,t_0}} 
 \left|Q \sqrt{h}(x't',y) \right|dy dx' dt' dt_0  \leq \\
C(T-\del) \int_{\Omega_\del^+} \left[ x'\left|\frac{\p V}{\p x'}\right|^2 
+t'\left|\frac{\p V}{\p t'}\right|^2+(x't')^2\left|\nabla_y V_j\right|^2\right]
\sqrt{h} (x't',y) dy dx' dt' + \\
C(T-\del)  \int_{\del}^{T} \int_{\p X} 
 {x'}\left| V\right|^2(x',x',y)\sqrt{h}({x'}^2,y) dy dx'.
\end{gather*}
After taking the limit as $\del \rightarrow 0,$
we find that there exists $K>0$ such that, for small $T$ 
\begin{gather}
\begin{gathered}
\int_{\Omega^+}
\left(  x't'(t'+x')|\nabla_y V|^2 +x'\left|\frac{\p V}{\p x'}\right|^2 
+t'\left|\frac{\p V}{\p t'}\right|^2 \right)
\sqrt{h}(x't',y) dy dx' dt' \leq \\
KT  \int_{0}^{T} \int_{\p X} \left( x'|V|^2 +
x'\left|\frac{\p V}{\p x'}\right|^2
+ {x'}^3\left|\nabla_y V\right|^2\right)(x',x',y)\sqrt{h}({x'}^2,y) dy dx'.
\end{gathered}\label{subsesv2} 
\end{gather} 
This does not quite give \eqref{est3T}, as the term in $|V|^2$ is not yet
included on the left hand side of the inequality.  However, taking
$\del\rightarrow 0$ in  \eqref{estv2},  substituting it in \eqref{subsesv2},
and choosing $T$ small enough,
 we obtain
\begin{gather}
\begin{gathered}
\int_{\Omega^+}
\left(|V|^2 +  x't'(t'+x')|\nabla_y V|^2 +x'\left|\frac{\p V}{\p x'}\right|^2 
+t'\left|\frac{\p V}{\p t'}\right|^2 \right)
\sqrt{h}(x't',y) dy dx' dt' \leq \\
K \int_{0}^{T} \int_{\p X} \left( x'|V|^2 +
x'\left|\frac{\p V}{\p x'}\right|^2
+ {x'}^3\left|\nabla_y V\right|^2\right)(x',x',y)\sqrt{h}({x'}^2,y) dy dx'.
\end{gathered} \label{est+}
\end{gather} 
As the operator in \eqref{eqmcq}
and the estimate \eqref{est+} remain of the same type after 
switching $x'$ and 
$t',$ this estimate also holds in the region below the diagonal. 
This ends the proof of the Lemma.
\end{proof}

Now we apply the Lemma \ref{enerest}
to prove Theorem \ref{rad}.   The goal is to
prove that for $f(x,y)$ smooth and compactly supported, the solution to
 \eqref{werad}
is smooth up to $\{x=0\}.$  We know from \eqref{supportv} 
that the change of variables \eqref{x't'} are
smooth on the support of $v,$ and we work in coordinates $(x',t').$  We will
show that $V(x',t',y)$ is smooth up to $\{x'=0\}\cup \{t'=0\}.$

  We first apply Lemma \ref{enerest} to the special case of
equation \eqref{eqq'}, noticing that in this case
$V$ is a single function instead of a 
vector. The data on the surface $\{x'=t'\}$ is given by \eqref{eqq'}, so
the integral on the right hand side of \eqref{est3T} is equal to
\begin{gather*}
\int_{0}^{T} 
\int_{\p X}  2x'|{x'}^{-n-1}f({x'}^2,y)|^2\sqrt{h}({x'}^2,y) dy dx'=
\int_{0}^{T^2} 
\int_{\p X} |f(x,y)|^2\sqrt{h}(x,y) \frac{dy dx'}{x^{n+1}}
\leq ||f||_{L^2(X)}^2.
\end{gather*}
Thus the right hand side of \eqref{est3T} is bounded
by the square of the norm of $f$ in $L^2(X,\dvol_g)$ and
\begin{gather*}
\int_{ \Omega_T}\left[ |V|^2+ x't'(x'+t')\left|\nabla_y V\right|^2+
x'\left|\frac{\p V}{\p x'}\right|^2 +t'\left|\frac{\p V}{\p t'}\right|^2 
\right]  \sqrt{h}(x't',y) dy dx' dt'
\leq  C ||f||_{L^2(X)}^2.
\end{gather*}

Now we want to obtain such energy estimates for the derivatives of $V$. We 
begin by analyzing the derivatives of $V$ in the $y$ variables and thus we
differentiate equation \eqref{eqq'}
with respect to $y.$
We  get a system of equations of the form
\begin{gather*}
\mcq \mcv= 0
\end{gather*}
where $\mcv$ is a $(n+1)\times 1$ vector whose transpose is
\begin{gather*}
\mcv^T=(V, \frac{\p V}{\p y_1}, ...., \frac{\p V}{\p y_n}),
\end{gather*}
and $\mcq$ is a matrix of second order operators with
principal part
\begin{gather*}
\mcq_2= \left(\frac{\p }{\p x'}\frac{\p }{\p t'} + x't'\Delta_h\right)
 \Id_{(n+1)\times (n+1)},
\end{gather*}
and lower order terms as in in \eqref{eqmcq}. So we conclude that

\begin{gather*}
\int_{ \Omega}\left[ |\mcv|^2+ x't'(x'+t')\left|\nabla_y \mcv\right|^2+
x'\left|\frac{\p \mcv}{\p x'}\right|^2 +t'\left|\frac{\p \mcv}{\p t'}\right|^2 
\right]  \sqrt{h}(x't',y) dy dx' dt'
\leq \\
C \left(||f||_{L^2(X)}^2+||\nabla_y f||_{L^2(X)}^2\right).
\end{gather*}

Using this argument repeatedly we conclude that

\begin{gather}
\begin{gathered}
\sum_{|\alpha|\leq k} \left[ \int_{ \Omega}\left[ 
\left|\left(\frac{\p}{\p y}\right)^{\alpha} V \right|^2+ 
x'\left|\frac{\p }{\p x'} \left(\frac{\p}{\p y}\right)^\alpha V \right|^2 
+ t'\left|\frac{\p}{\p t'}\left(\frac{\p}{\p y}\right)^\alpha V\right|^2 
\right]  \sqrt{h}(x't',y) dy dx' dt'\right]
\leq \\
C\sum_{|\alpha|\leq k+1} 
||\left(\frac{\p}{\p y}\right)^\alpha f||_{L^2(X)}^2, 
\;\ k \in \mn,
\end{gathered}\label{eqq'k}
\end{gather}

Next we use the equation to obtain information about 
the derivatives of $V$ with respect to $x'$ and $t'.$ It is convenient to
get rid of the first order terms in \eqref{eqq'}, and we do that by 
conjugating the operator with the factor $|h(x't',y)|^{-\oq}.$
Setting
\begin{gather}
 W(x',t',y)=|h|^{\oq}V(x',t',y), \label{defw}
\end{gather} then
$W$ satisfies
\begin{gather}
\left(\frac{\p^2}{\p x' \p t'}+ x't'\Delta_h +x't'B(x't',y,\frac{\p}{\p y})+
C(x't',y)\right)W=0,\label{eqq'W}
\end{gather}
where $B$ is a first order operator and $C$ is a smooth function. Since
$|h|$ is smooth and positive, $V$ and $W$ have the same regularity given by
\eqref{eqq'k}.  Moreover $W$ is also supported in $\{x'\geq \sqrt{x_0}\}\cup
\{t'\geq \sqrt{x_0}\}.$

  We have shown in
\eqref{eqq'k} that
 ${t'}^\ha \frac{\p }{\p t'}\frac{\p^\alpha }{\p y^\alpha}W \in L^2(\Omega)$
for any $\alpha.$ 
In particular
if $\Omega^+=\Omega \cap \{t'>x'\},$ 
${t'}^\ha \frac{\p }{\p t'}\frac{\p^\alpha }{\p y^\alpha}W \in L^2(\Omega^+).$ 
In view of the support of $W$ we have 
$\frac{\p }{\p t'}\frac{\p^\alpha }{\p y^\alpha}W \in L^2(\Omega^+).$ 

We deduce from \eqref{eqq'k} and
\eqref{eqq'W} that
\begin{gather}
 \left(\frac{\p}{\p y}\right)^\alpha W(x',t',y), \;\
\frac{\p}{\p t'} \left(\frac{\p}{\p y}\right)^\alpha W(x',t',y), \;\
\left(\frac{\p}{\p y}\right)^\alpha \frac{\p }{\p t'}\frac{\p }{\p x'}
  W(x',t',y)\in L^2(\Omega^+) \;\ \forall \;\ \alpha \in \mn^n
\label{bregw}
\end{gather}
Writing, for $t'>x',$
\begin{gather*}
\frac{\p}{\p x'}\left(\frac{\p}{\p y}\right)^\alpha W(x',t',y)=
 \frac{\p}{\p x'}\left(\frac{\p}{\p y}\right)^\alpha W(x',x',y) +
\int_{x'}^{t'}\frac{\p}{\p \mu} \frac{\p}{\p x'}\left(\frac{\p}{\p y}\right)^\alpha 
 W(x',\mu,y) \; d\mu,
\end{gather*}
we see that, for $t'>x',$
\begin{gather*}
\left|\frac{\p}{\p x'}\left(\frac{\p}{\p y}\right)^\alpha  W(x',t',y)\right|^2
\leq \\ 2 
\left|\frac{\p}{\p x'}\left(\frac{\p}{\p y}\right)^\alpha W(x',x',y)\right|^2 +
2(t'-x')\int_{x'}^{t'}\left|\frac{\p}{\p \mu}
 \frac{\p}{\p x'}\left(\frac{\p}{\p y}\right)^\alpha W(x',\mu,y)\right|^2 
 d\mu.
\end{gather*}
Integrating this in $\Omega^+,$  using that 
$\frac{\p}{\p x'}\frac{\p^\alpha}{\p y^\alpha}
W(x',x',y)$ is smooth with compact support, and \eqref{bregw}, we find that
\begin{gather}
 \frac{\p}{\p x'}\left(\frac{\p}{\p y}\right)^\alpha  W 
\in L^2(\Omega^+), \;\ \alpha \in \mn^n. \label{estax}
\end{gather}
So \eqref{bregw} and \eqref{estax} show that
\begin{gather*}
\left(\frac{\p}{\p y}\right)^\alpha  W , \;\
 \frac{\p}{\p x'}\left(\frac{\p}{\p y}\right)^\alpha  W \;\
\text{ and } \;\ \frac{\p }{\p t'}\left(\frac{\p}{\p y}\right)^\alpha  W 
\in L^2(\Omega^+), \;\ \alpha \in \mn^n.
\end{gather*}
Using the symmetry of $W$ with respect to the diagonal, we in fact have
\begin{gather*}
\left(\frac{\p}{\p y}\right)^\alpha  W , \;\
 \frac{\p}{\p x'}\left(\frac{\p}{\p y}\right)^\alpha  W \;\
\text{ and } \;\ \frac{\p }{\p t'}\left(\frac{\p}{\p y}\right)^\alpha  W 
\in L^2(\Omega), \;\ \alpha \in \mn^n.
\end{gather*}
Differentiating \eqref{eqq'W} first with respect to $y$ and then
with respect to $x'$ or $t'$ we find that 
\begin{gather*}
\frac{\p^2}{\p {x'}^2}\frac{\p}{\p {t'}}\left(\frac{\p}{\p y}\right)^\alpha  W
= \sum_{\beta}\sum_{ m=0,1} G_{\beta,m}(x',t',y)
\frac{\p^m}{\p {x'}^m}\left(\frac{\p}{\p y}\right)^\beta  W \in L^2(\Omega^+),
\\
\frac{\p^2}{\p {t'}^2}\frac{\p}{\p {x'}}\left(\frac{\p}{\p y}\right)^\alpha  W
= \sum_{\beta}\sum_{ m=0,1} F_{\beta,m}(x',t',y)
\frac{\p^m}{\p {t'}^m}\left(\frac{\p}{\p y}\right)^\beta  W \in L^2(\Omega^+),
\end{gather*}
with $F_{\beta,m}$ and $G_{\beta,m}$ smooth. Thus
\begin{gather*} 
\frac{\p^2}{\p {x'}^2}\frac{\p}{\p {t'}}\left(\frac{\p}{\p y}\right)^\alpha W,
\;\ 
\frac{\p^2}{\p {t'}^2}\frac{\p}{\p {x'}}\left(\frac{\p}{\p y}\right)^\alpha  W
\in L^2(\Omega^+).
\end{gather*}

Proceeding as above, we find that
\begin{gather*}
\frac{\p^2}{\p {x'}^2}\left(\frac{\p}{\p y}\right)^\alpha W,\;\
\frac{\p^2}{\p {t'}^2}\left(\frac{\p}{\p y}\right)^\alpha W, \;\
\text{ and } \;\ 
\frac{\p}{\p {x'}}\frac{\p}{\p {t'}}\left(\frac{\p}{\p y}\right)^\alpha W
\in  L^2(\Omega^+).
\end{gather*}

Using the symmetry of $W$ with respect to the diagonal we get that in fact
\begin{gather*}
\frac{\p^2}{\p {x'}^2}\left(\frac{\p}{\p y}\right)^\alpha W,\;\
\frac{\p^2}{\p {t'}^2}\left(\frac{\p}{\p y}\right)^\alpha W, \;\
\text{ and } \;\ 
\frac{\p}{\p {x'}}\frac{\p}{\p {t'}}\left(\frac{\p}{\p y}\right)^\alpha W
\in  L^2(\Omega).
\end{gather*} 
This argument can be repeated for all derivatives with respect to $x'$ and 
$t'$ and we conclude that
$W,$ and hence $V,$ is smooth in 
$[0,T]\times [0,T]\times\p X$ for $T$ small. Any
interval $[0,T']$ can be covered by small intervals in which the method
above can be applied.  So in fact this shows that the solution is
smooth in $[0,T]\times [0,T]\times\p X$ for any $T.$

Since $V$ is supported in $\{x'>\sqrt{x_0}\} \cup \{t'>\sqrt{x_0}\} $ 
and the change of coordinates \eqref{x't'} is smooth in this region, this
shows that $v$ 
has a smooth extension up to $\{x=0\}.$
This concludes the proof of Theorem~\ref{rad}.
\end{proof}
 
\section{The Radiation Fields And The Scattering Matrix}\label{Scat}

The spectrum of the Laplacian $\sigma(\Delta)$ was studied by Mazzeo and
Mazzeo and Melrose \cite{ma,ma1,mame}, see also section 3 of \cite{grazw} for 
a discussion. 
It consists of a finite pure
point spectrum $\sigma_{\pp}(\Delta),$ which is the set of $L^2(X)$ 
eigenvalues, and an absolutely continuous spectrum
$\sigma_{\ac}(\Delta)$ satisfying
\begin{gather}
\sigma_{\ac}(\Delta)= \left[ \frac{n^2}{4},\infty\right), \;\
\text{ and } \;\ \sigma_{\pp}(\Delta)\subset 
\left( 0, \frac{n^2}{4}\right). \label{spect}
\end{gather}

This gives a decomposition 
\begin{gather*}
L^2(X)= L^2_{\pp}(X) \oplus L^2_{\ac}(X),
\end{gather*}
where $L^2_{\pp}(X)$ is the finite dimensional 
space spanned by the eigenfunctions and
$L^2_{\ac}(X)$ is the space of absolute continuity, which is the orthogonal
complement of $L^2_{\pp}(X).$

 With  this choice of the spectral parameter, $\frac{n^2}{4}+\la^2,$
 we have that if
$\Im \la \not=0,$ then $\frac{n^2}{4}+\la^2 \not \in [\frac{n^2}{4}, \infty),$
while if $\Im \la <-\frac{n}{2}$ it follows that 
$\frac{n^2}{4}+\la^2 \not \in [0, \infty).$ 
The eigenvalues of $\Delta$ are  of finite multiplicity and
are described by points on
the line  $\Re \la=0$ and $-\frac{n}{2}<\Im \la <0.$
As the Laplacian is a non-negative operator, the spectral theorem
gives that the resolvent
\begin{gather}
R(\frac{n}{2}+i\la)=\left(\Delta-\frac{n^2}{4}-\la^2\right)^{-1}
: L^2(X)\longrightarrow L^2(X), \;\  \text{ provided } \;\
\Im \la <-\frac{n}{2}. \label{resb}
\end{gather}
 It was shown in
\cite{mame} that it can be meromorphically continued to 
$\mc\setminus\frac{i}{2} \mn$ as an operator acting on appropriate spaces.

Let
\begin{gather*}
H_E(X)=\{(f_1,f_2): \;\ f_1, \; f_2\in L^2(X), \;\ \text{ and } \;
 d f_1 \in L^2(X) \}.
\end{gather*}

For $w_0, \;\ w_1 \in C_0^\infty(\intx),$ we define
the energy of $w=(w_0,w_1)$ by
\begin{gather}
||w||_E^2=\ha \int_{X} \left( |d w_0|_g^2 -\frac{n^2}{4}|w_0|^2
+|w_1|^2 \right)\;\ d\vol_g, \label{nore}
\end{gather}
where $|d w_0|_g$ denotes the length of the co-vector with respect to the
metric induced by $g$ on $T^*X.$ 
But $||w||_E^2$ is only positive
when $w_0$ is in the space of absolute continuity of $\Delta,$ and only then 
it defines a norm. 
We denote
\begin{gather*}
\mcp_{\ac} : L^2(X) \longrightarrow L^2_{\ac}(X)
\end{gather*}
the corresponding projector. Let 
\begin{gather*}
E_{\ac}(X)= \mcp_{\ac}\left(H_E(X)\right)=\text{ Range of the projector } \;
\mcp_{\ac} \; \text{ acting on } \;\ H_E(X).
\end{gather*}
$E_{\ac}(X)$ is a Hilbert space equipped with the norm \eqref{nore}.

One can use integration by parts to prove that if $u(t,z)$ satisfies 
\eqref{we2}, then
$$||\left(u(t,\bullet),D_tu(t,\bullet)\right)||_E=
||\left(u(0,\bullet),D_tu(0,\bullet)\right)||_E.$$
Since $W(t)$ commutes with $\mcp_{\ac}$ this  gives, by for example slightly 
modifying the proof 
of Proposition 2.24 of \cite{fried1}, that the map  $W(t)$ defined by
\begin{gather}
\begin{gathered}
W(t):  C_0^\infty(\intx)\times C_0^\infty(\intx)
 \longrightarrow
C_0^\infty(\intx)\times C_0^\infty(\intx) \\
W(t)\left(f_1,f_2\right)=\left(u(t,z), D_t u(t,z)\right), \;\ 
 t \in \mr 
\end{gathered}\label{group}
\end{gather}
induces a strongly continuous group of unitary operators. 
$$W(t):E_{\ac}(X) \longrightarrow E_{\ac}(X), \;\  t\in \mr.$$

By changing $t\leftrightarrow t-\tau,$ one has that $\mcr_{\pm}$ 
satisfy
\begin{gather}
\mcr_{\pm} \circ \left(W(\tau)f\right)(s,y)=\mcr_{\pm} f(s+\tau,y), 
\;\ \tau \in \mr.
 \label{transl}
\end{gather}

So Theorem \ref{rad}
shows that $\mcr_{\pm}$ are ``twisted'' translation representations of
the group $W(t).$ That is, if one sets 
$\widetilde{\mcr_{\pm}}(f)(s,y)=\mcr_{\pm} f(- s,y),$ then
\begin{gather}
\widetilde{\mcr_{\pm}}( W(\tau))= T_\tau \widetilde{{\mcr_{\pm}}}, \label{transl1}
\end{gather}
where $T_\tau$ denotes right translation by $\tau$ in the $s$ variable. So
 $\widetilde{\mcr_{\pm}}$ are  translation representers in the sense of 
Lax and Phillips. 
Moreover  we will prove
\begin{thm}\label{inv} The maps $\mcr_{\pm}$ induce isometries
\begin{gather*}
\mcr_{\pm} : E_{\operatorname{ac}}(X)\longmapsto L^2(\mr \times \p X),
\end{gather*} 
where $L^2(\mr \times \p X)$ is defined with respect to $h_0$ fixed in 
\eqref{hmet1}.
\end{thm}

We deduce from Theorem \ref{inv} that the  {\it scattering operator}
\begin{gather}
\begin{gathered}
\mcs: L^2(\mr \times \p X) \longrightarrow L^2(\mr\times \p X) \\
\mcs=\mcr_+\circ \mcr_{-}^{-1}
\end{gathered}\label{sco}
\end{gather}
is unitary  in $L^2(\p X \times \mr),$ and in view of \eqref{transl1}, it
commutes with translations. 
  This implies that the Schwartz kernel 
$\mcs(s,y,s',y')$ of $\mcs$ satisfies
\begin{gather*}
\mcs(s,y,s',y')=\mcs\left(s-s',y,y'\right),
\end{gather*}
and thus is a convolution operator.
The {\it scattering matrix} is defined by conjugating $\mcs$ with the partial 
Fourier transform in the $s$ variable
\begin{gather*}
\mca=\mcf \mcs \mcf^{-1}.
\end{gather*}
$\mca$ is a unitary operator in $L^2(\p X \times \mr)$ and, since $\mcs$ acts 
as a 
convolution in the variable $s,$ $\mca$ is a multiplication in the 
variable $\la,$ i.e. it satisfies
\begin{gather}
\mca F(\la,y)= \int_{\p X} \mca(\la,y,y')  F(\la,y') \; d\vol_{h_0}(y').
\label{scatmat1}
\end{gather}

We will
prove  that the stationary and dynamical definitions of the scattering
matrix are equivalent:
\begin{thm}\label{equiv}  With $x$ given by \eqref{hmet1} and $\la\not=0,$
 the Schwartz kernel of the map $A(\la)$ defined by \eqref{scatmat},
 is equal to $\mca(\la,y,y'),$ defined in \eqref{scatmat1}.
\end{thm}

To prove Theorem \ref{inv} and  Theorem \ref{equiv}
we will use the connection between the wave equation, the resolvent,
and the Eisenstein function from \cite{jsb,mame}.

\section{The Radiation Fields and the Eisenstein Function}\label{mmc}

The following is an
 important fact concerning the behavior of the solution  to the wave
equation:
\begin{prop}\label{l2gro} Let $f=(f_1,f_2)\in H_E(X)$ and let
$u(t,z)$ be the solution to \eqref{we2} with initial data 
$f.$ There exists $C=C(f)>0$ such that 
\begin{gather}
\int_X |\frac{\p u(t,z)}{\p t}|^2 \;\ \dvol_g(z) \leq C e^{\frac{n}{2} t}, 
\;\ t >0. \label{l2gron}
\end{gather}
\end{prop}
\begin{proof} In view of  \eqref{spect} we can write

\begin{gather*}
f_j(z)= \sum c_{j,k} \phi_k(z) + g_j(z)
\end{gather*}
where $\phi_k$ is an eigenfunction associated with an eigenvalue
$\sigma_k=\frac{n^2}{4}+\la_k^2\in (0,\frac{n^2}{4})$ and $g_j$ is the 
projection of $f_j$ onto $L^2_{\ac}(X).$

We then have that 
\begin{gather*}
u(t,z)= u_{\ac}(t,z)+\sum_{k}
\left( c_{1,k} \cos(\la_k t)+ c_{2,k} \frac{1}{\la_k}
\sin(\la_k t) \right)\phi_k(z)
\end{gather*}
where $u_{ac}$ is the solution to \eqref{we2} with data $g=(g_1,g_2).$

As 
\begin{gather*}
u_{\ac}(t,z)= \cos t \sqrt{\Delta-\frac{n^2}{4}} \; g_1 +
\left(\Delta-\frac{n^2}{4}\right)^{-\ha}\sin t \sqrt{\Delta-\frac{n^2}{4}} \; 
g_2, \end{gather*}
it follows that
\begin{gather*}
\mcp_{\ac}u(t,z)= \cos t \sqrt{\Delta-\frac{n^2}{4}} \mcp_{\ac}f_1 +
\left(\Delta-\frac{n^2}{4}\right)^{-\ha}\sin t \sqrt{\Delta-\frac{n^2}{4}} 
\mcp_{\ac} f_2=u_{\ac}(t,z).
\end{gather*}
In this case the energy 
\begin{gather*}
\int_{X} \left(|\frac{\p u_{\ac}(t,z)}{\p t}|^2 +|\nabla_y u_{\ac}(t,z)|^2-
\frac{n^2}{4}|u_{\ac}(t,z)|^2\right) \dvol_g(z)=\\
\int_{X} \left(|g_2(z)|^2 +|\nabla_y g_1(z)|^2-
\frac{n^2}{4}|g_1(z)|^2\right) \dvol_g(z)=E_0>0,
\end{gather*}
and in particular this gives
\begin{gather*}
\int_X |\frac{\p u_{\ac}(t,z)}{\p t}|^2 \;\ \dvol_g(z) \leq C.
\end{gather*}
Since $0>\Im \la_k>-\frac{n}{2},$ 
the other part of the solution obviously satisfies \eqref{l2gron} and this
proves the proposition.
\end{proof}

Next we present an elementary and useful Lemma. The proof we  
gave of this result in the first version of this paper was incorrect. One of 
the referees 
pointed out the mistake and kindly provided us with  the proof we present.
\begin{lemma}\label{linalg} Let $\mch$ be an infinite dimensional
 Hilbert space, let $H$ be a subspace of $\mch$ of finite dimension and
let $H^\bot$ be the orthogonal to $H.$ If
 $E\subset \mch$ is a dense subspace of $\mch,$ then 
$E\cap H^\bot \not=\emptyset,$ and $E\cap H^\bot$ is dense in $H^\bot.$
\end{lemma}
\begin{proof}  
 We begin with the case $\dim H=1$ and we will prove that
$\left(E\cap H^\bot\right)^\bot=H.$   It is easy to see that 
$H\subset \left(E\cap H^\bot\right)^\bot.$  
To prove  that $\left(E \cap H^\bot\right)^\bot\subset  H$ we let
 $H=\pan\{\phi\},$ $\phi\not=0.$
As $E$ is dense, pick $f_1, f_2\in E$ such  that $\lan f_j, \phi \ran \not=0,$
$j=1,2.$ Then
$$ f=f_2-\frac{\lan f_2,\phi\ran}{\lan f_1,\phi\ran}  f_1 \in E\cap H^\bot.$$
Thus, if 
$u \in \left(E \cap H^\bot\right)^\bot,$ then 
$ \lan u, f \ran =0.$ But  this can  be  rewritten as
\begin{gather*}
\lan f_2, u-\frac{\lan f_1, u\ran}{\lan f_1, \phi\ran} \phi\ran =0, \;\
\text{for all }  \;\ f_2 \;\ \text{ with } \;\ \lan f_2,\phi\ran \not=0
\end{gather*}
Since $E$ is dense, this in  fact holds for every $f_2$ and
thus $u=\frac{\lan f_1, u\ran}{\lan f_1, \phi\ran} \phi.$ Therefore
$\left(E \cap H^\bot\right)^\bot\subset  H.$

The general case follows by induction. Suppose the result  is true for
$\dim H=N-1$  and let $\dim H=N.$ Then
\begin{gather*}
H= H_{N-1} \oplus D, \;\ H_{N-1}=\pan\{\phi_1,..., \phi_{N-1}\}
\;\ \text{ and } \;\ D=\pan\{\phi_N\}, \;\ \lan \phi_i,\phi_j\ran=0,
\;\ 1\leq i,j \leq N.
\end{gather*}

As proved above, $E_1=E\cap D^\bot$ is dense in $\mch_1=D^\bot.$
By assumption the result holds in  dimension $N-1,$
so $E_1\cap H_{N-1}^\bot$ is dense in the orthogonal to $H_{N-1}$ in
$D^\bot,$ which is the orthogonal to $H_N$ in $\mch.$ 

As $E_1\cap H_{N-1}^\bot=E\cap H_N^\bot,$ $E\cap H_N^\bot$ is dense
in $H_N^\bot$ in $\mch.$

\end{proof}

The most important consequence of this is
\begin{cor}\label{density} For $L^2_{\ac}(X)$ defined as above,
$C_0^\infty(\intx)\cap L^2_{\ac}(X)$ is dense in $L^2_{\ac}(X).$
\end{cor}
Now we are ready to prove the first mapping property of $\mcr_+.$

\begin{prop}\label{radfl2} Let $u$ be the solution to \eqref{we2} with
data $f=(f_1,f_2)\in C_0^\infty(\intx) \cap E_{\ac}(X).$ Then 
$\mcr_+ f(s,y)\in L^2(\mr\times \p X)$ and
\begin{gather*}
||\mcr_+ f||_{L^2(\mr\times \p X)}\leq 2 ||f||_E.
\end{gather*}
\end{prop}
\begin{proof} 
We modify the proof of Lemma 2.6 of \cite{fried1}. Let $u$ be the solution
of \eqref{we2} with initial data equal to $f.$ For $t, T>0$ fixed, 
consider the integral
\begin{gather*}
E_{T,t}=\int_{\{z=(x,y)\in X; \; t+\log x<T\}} 
\left|\frac{\p u}{\p t}(x,y,t)\right|^2 \;\ d\vol_g.
\end{gather*}

For $t+\log x <T,$ it follows that $x \rightarrow 0$ as 
$t\rightarrow \infty.$ So, for $t$ large,
we may work in local coordinates where 
\eqref{hmet1} holds.
First we prove that
\begin{gather}
\lim_{t \rightarrow \infty} E_{T,t}= 
\int_{-\infty}^{T} \int_{\p X} \left|\mcr_+ f\right|^2\;\ d\vol_{h_0}\; ds , \;\ 
\forall \; T\in \mr.
\label{lim1p}
\end{gather}
To see this we set
\begin{gather*}
u(t,x,y)=x^{\frac{n}{2}}v(t+\log x,x,y),
\end{gather*} 
and since $s=t+\log x,$
we have $\frac{\p u}{\p t} = x^{\frac{n}{2}} \frac{\p v}{\p s}.$
As $ds=\frac{dx}{x},$
it follows that
$d\vol_g=\sqrt{|h|}(x,y)\frac{dx dy}{x^{n+1}}=   \sqrt{|h|}(e^{(s-t)},y)
e^{n(t-s)}\; ds\; dy.$ Since the initial data is compactly supported,
there exists $s_0$ such that $v=0$ for $s<s_0.$ Thus
\begin{gather*}
E_{T,t}=\int_{s_0}^{T}\int_{\p X}
\left|\frac{\p v}{\p s}\right|^2(e^{s-t},y,s)
\sqrt{|h|(e^{(s-t)},y)} \; dy \;ds.
\end{gather*}
From Theorem \ref{rad} the convergence is uniform and thus we obtain
\begin{gather*}
\lim_{t \rightarrow \infty} E_{T,t}=
 \int_{-\infty}^{T}\int_{\p X}
\left|\frac{\p v}{\p s}(0,s,y)\right|^2 \sqrt{|h|(0,y)} \; ds dy= 
 \int_{-\infty}^{T}\int_{\p X}
\left|\mcr_+ f(s,y)\right|^2 \; d\vol_{h_0}\; ds <\infty \;\ \; \forall \; T\in \mr.
\end{gather*}

For any $(f_1,0)\in E_{\ac}(X),$
$\int_{X} \left(|d f_1|_{g}^2-
\frac{n^2}{4}| f_1|^2\right)\; d\vol_g >0.$ 
Thus the result follows by conservation of energy.
\end{proof}

In view of Corollary \ref{density} this can be restated as
\begin{cor}\label{l2ext} The maps $\mcr_{\pm}$ defined in \eqref{mcr}
and \eqref{mcr-} extend from 
$(C_0^\infty(\intx) \times C_0^\infty(\intx)) \cap E_{\ac}(X)$
by continuity to maps
\begin{gather*}
\mcr_{\pm}: E_{\ac}(X) \longrightarrow L^2(\mr \times \p X).
\end{gather*}
\end{cor}

The estimate \eqref{l2gron} shows that we can take partial
Fourier-Laplace 
transform with  respect to $t$ of the forward fundamental solution of
the wave equation  $U(t,z,z')$ and thus, 
for $\Im \la <-\frac{n}{2},$ we denote
\begin{gather*}
R\left(\frac{n}{2}+i\la\right)=\left(\Delta-\la^2-\frac{n^2}{4}\right)^{-1}=
\int_{\mr} U(t,z,z') e^{-i\la t} dt
\end{gather*}

It is easy to see that if $u(t,z)$ satisfies \eqref{we2} 
then $V(t,z)=H(t)u(t,z)$  satisfies
\begin{gather*}
\left(D_{t}^2-\Delta+\frac{n^2}{4}\right)V(t,z)=if_2(z)\delta(t)
+if_1(z)D_t\delta(t),  \text{ on } \mr\times \intx  \\
V(t,z)=0 \text{ for } t<0,
\end{gather*}
and  $\lim_{x\rightarrow 0} x^{-\frac{n}{2}} D_s V(s-\log  x, x,y)=
\lim_{x\rightarrow 0} x^{-\frac{n}{2}} D_s u(s-\log  x, x,y).$

So,
\begin{gather}
R\left(\frac{n}{2}+i\la\right)
(if_2+i\la f_1)
=\int_{\mr} V(t,z) e^{-i\la t} dt, \;\ \text{ with } \;
f=(f_1,f_2) \in C_0^\infty(X) \text{ and } \;\ \Im \la < -\frac{n}{2}.
 \label{ft1}
\end{gather}

In fact, if $f=(f_1,f_2) \in C_0^\infty(X)\cap E_{\ac}(X),$ conservation of 
(positive) energy gives that \eqref{ft1} holds
for $\Im\la \leq 0.$
Using local coordinates $z=(x,y)$ near $\p X,$
we are interested in 
$\lim_{x\downarrow 0} x^{-\frac{n}{2}}D_sV(s-\log x,x,y).$
From \eqref{ft1} we obtain, for $\Im \la < -\frac{n}{2},$ 
\begin{gather}
\begin{gathered}
\int_{\mr}x^{-\frac{n}{2}}D_sV(s-\log x,y) e^{-i\la s} ds=
\int_{\mr}x^{-\frac{n}{2}}D_tV(t,x,y) e^{-i\la t-i\la \log x} dt=\\
x^{-\frac{n}{2}-i\la}\la\; R\left(\frac{n}{2}+i\la\right)(if_2+i\la f_1)(x,y).
\end{gathered}\label{ft2}
\end{gather}

Next we study the behavior of the Schwartz kernel 
\begin{gather}
x^{-\frac{n}{2}-i\la} R\left(\frac{n}{2}+i\la,z,z'\right),\;\ z=(x,y), \;\ 
z'=(x',y'), \;\
\Im\la<-\frac{n}{2},
\text{ as } x\downarrow 0. \label{limres}
\end{gather}  
 Mazzeo and Melrose show in \cite{mame} 
that $R(\frac{n}{2}+i\la)$ has a meromorphic continuation from 
$\Im \la<-\frac{n}{2}$
 to $\mc\setminus \frac{i}{2}\mn.$ We briefly recall their construction and
use it to study \eqref{limres}.

Locally, in the interior of $X \times X,$ and for $\Im \la <<0,$ 
 $R(\frac{n}{2}+i\la)$ is pseudo-differential operator, so
its kernel is singular at the diagonal
$$D=\{(x,y,x',y') \in X \times X; x=x', y=y'\}.$$
To understand the behavior of the kernel of $R(\frac{n}{2}+i\la)$ up to
$$D_{\p X}=D \cap( \p X \times \p X),$$ and for other values of $\la,$
 Mazzeo and Melrose blow-up the intersection $D_{\p X} .$ 
This can be done in an invariant way, but in local coordinates
this can be  seen as introduction of polar coordinates around
$D_{\p X} .$  Taking coordinates $(x,y)$ and $(x',y')$ 
in a product decomposition of each copy of $X$ near $\p X,$   
the ``polar coordinates' are then given
 by  
\begin{gather}
R=[x^2+{x'}^2+ {|y-y'|}^2]^\ha, \;\ \rho=\frac{x}{R}, \;\
\rho'=\frac{x'}{R} \;\ \omega=\frac{y-y'}{R} \label{coord}
\end{gather}

A function is smooth in the
space $X \times_{0} X$ if it is smooth in 
 $(R,\rho, \rho',y, \omega)$  about
$D_{\p X} .$ 
As a set, $X \times_{0} X$ is $X \times X$ with 
$D_{ \p X}$ replaced by the interior pointing portion of its normal
bundle. Let
$$ \beta :X \times_{0} X \longrightarrow X \times X $$
denote the blow-down map.

\begin{figure}[int1]
\epsfxsize= 3.5in
\centerline{\epsffile{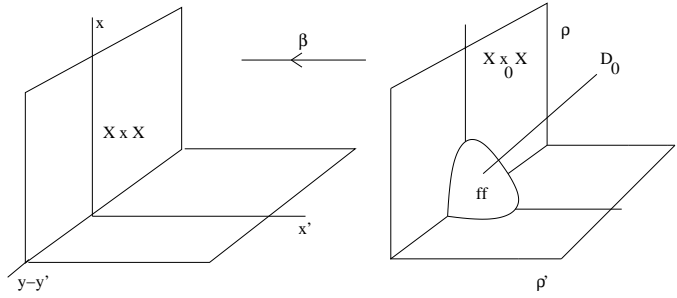}}
\caption{The blown up manifold $X\times_0 X.$}
\label{blow}
\end{figure}

The function $R$ is a defining function for a new face, which is called
the front face, $\ff.$ This is the lift of 
 $D_{\p X}=D \cap \left( \p X \times \p X\right).$ 
The functions $\rho$ and $\rho'$ are then
defining functions for the other two boundary faces which are called the
top face $\mathcal{T},$ and bottom face $\mathcal{B},$ respectively, i.e. 
\begin{gather*}
\ff=\{R=0\}, \;\ \mathcal{B}=\{\rho'=0\}, \;\
\mathcal{T}=\{\rho=0\}.
\end{gather*}
See Figure \ref{blow}, which is taken from section 3 of \cite{mame}. In 
$X \times_{0} X$ 
the lift of the diagonal of $X\times X$ only meets the 
boundary $\ff$ and is disjoint from the other two boundary faces.

The main result from \cite{mame} needed here is that the lift of the kernel of
the resolvent satisfies
\begin{gather}
\beta^*R(\frac{n}{2}+i\la)=R_1(\frac{n}{2}+i\la)+R_2(\frac{n}{2}+i\la), \label{res1}
\end{gather}
where $R_1$  is conormal of order $-2$  to the lifted diagonal,
 $D_0,$ and smooth up to the front face, and, most importantly for our 
purposes, vanishes to infinite order at the top and bottom faces.
The second part, $R_2,$ is of the form
\begin{gather}
R_2= \rho^{\frac{n}{2}+i\la}{\rho'}^{\frac{n}{2}+i\la}
F(\frac{n}{2}+i\la,\bullet),
\label{res2}
\end{gather}
where
$F(\frac{n}{2}+i\la,\bullet)
 \in C^{\infty}\left(X\times_0 X \right)$
is meromorphic in $\la\in \mc\setminus\frac{i}{2}\mn,$ and holomorphic
and in $\la \in \Im \la< -\frac{n}{2}.$

Now we use the construction above to verify that
${x'}^{-\frac{n}{2}-i\la}R(\frac{n}{2}+i\la,z,z')|_{\{x=0\}}$ is well defined.
This is carried out in Proposition 4.1 of \cite{jsb}, and we briefly 
describe it.

We first look at the lift of $x^{-\frac{n}{2}-i\la}R(\frac{n}{2}+i\la,z,z')$ under the 
blow-down map $\beta.$  Recalling that $x=R\rho,$ we deduce from \eqref{res1} 
and \eqref{res2}, and the
fact that $R_1(\frac{n}{2}+i\la)$ vanishes to infinite order at the top and 
bottom faces, that
\begin{gather}
\begin{gathered}
\beta^*\left(x^{-\frac{n}{2}-i\la}R(\frac{n}{2}+i\la,z,z')\right)|_{\rho=0}= 
(R\rho)^{-\frac{n}{2}-i\la}\left(R_1(\frac{n}{2}+i\la)+R_2(\frac{n}{2}+i\la)\right)|_{\rho=0}= \\
\left(R^{-\frac{n}{2}-i\la}{\rho'}^{\frac{n}{2}+i\la}
F_2(\frac{n}{2}+i\la,\rho,\rho',R,\omega)\right)|_{\rho=0}
\end{gathered}\label{rest}
\end{gather}
Since $F_2(\frac{n}{2}+i\la, \bullet)$
is smooth,   the
restriction of $\rho^{\frac{n}{2}+i\la}F_2$ to $\rho=0$ is 
well defined. 

Notice that the composition of the restriction to $\rho=0$ with the
blow-down map $\beta$ gives a map
\begin{gather*}
\beta_1: \p X \times_0 X \longrightarrow \p X \times X,
\end{gather*}
where $\p X \times_0 X$ denotes the manifold $\p X \times X$ with the 
submanifold $\{x'=0, y=y'=0\}$ blown-up. By an abuse of notation we continue 
to use $R$ and $\rho'$ to denote the restrictions of the corresponding 
functions to $\rho=0.$ In other words we consider
\begin{gather*}
R=[{x'}^2+ {|y-y'|}^2]^\ha, \;\ \rho'=\frac{x'}{R} \;\ \omega=\frac{y-y'}{R}
\end{gather*}

Therefore \eqref{rest} above gives a well defined distribution on the 
manifold $\p X \times_0 X,$ namely
\begin{gather*}
F=R^{-i\la}{\rho'}^{\frac{n}{2}+i\la}
F_2(\frac{n}{2}+i\la,\rho,\rho',R,\omega)|_{\rho=0}
\end{gather*}
By switching the variables $x$ and $x'$ in the proof of Proposition 4.1 
of \cite{jsb} it can be shown that the push forward of $F$ to
 $\p X \times X$ under the map $\beta_1$ is well defined.
In view of \eqref{res1} and
\eqref{res2} this
is essentially the same thing we want to do here.
Therefore
the restriction $x^{-\frac{n}{2}-i\la}R(\frac{n}{2}+i\la,z,z')|_{x=0}$ is well defined 
and we shall denote it by 
\begin{gather}
E(\frac{n}{2}+i\la,y,z')\defi x^{-\frac{n}{2}-i\la}R(\frac{n}{2}+i\la,z,z')|_{x=0}={\beta_1}_{*}F, \;\ \la \in \mc. \label{defe}
\end{gather}
It is clear from \eqref{res2} that $E(\frac{n}{2}+i\la,y,z')$ has a conormal 
singularity at $\{x'=0,\; y=y'\}.$ 
This is the transpose of the Eisenstein function, or Poisson operator,
see for example \cite{gui,gz0,jsb}. 

This defines a meromorphic family of operators
\begin{gather}
\begin{gathered}
E(\frac{n}{2}+i\la) : C_{0}^{\infty}\left(\intx\right) \longrightarrow 
C^{\infty}\left(\p X\right) \\
E(\frac{n}{2}+i\la)f(y)= \int_{X} E(\frac{n}{2}+i\la,y,z')f(z') \; 
d\vol_g(z'), \;\ \la \in \mc\setminus\frac{i}{2}\mn.
\end{gathered}\label{mereisn}
\end{gather}

So we conclude from \eqref{defe} and
\eqref{ft1} that 
\begin{gather}
\begin{gathered}
\int_{\mr} e^{-i\la s}\mcr_+(f_1,f_2)(y,s)ds=
\left.\int_{\mr} x^{-\frac{n}{2}}
(D_sV)(s-\log x,z) e^{-i\la s} ds\right|_{x=0}= \\ 
i \la \int_{X} E(\frac{n}{2}+i\la,y,z')(f_2(z')+\la f_1(z'))\; d\vol_g(z'), \\ 
\;\ \Im \la < -\frac{n}{2} \text{ if } \;\ f_1, f_2 \in C_0^\infty(X),\;\
\text{ and } \;\ 
 \Im \la \leq 0 \;\ \text{ if } \;\ f_1, f_2 \in C_0^\infty(X)\cap E_{\ac}(X). 
\end{gathered}\label{einrad+}
\end{gather}
In view of \eqref{mereisn}  the right hand side of this
equation has a meromorphic continuation to 
$\mc\setminus\frac{i}{2}\mn,$ so
the left hand side can be meromorphically continued to 
$\mc\setminus\frac{i}{2}\mn.$

Recall from \eqref{radfl2} that $\mcr_+(f_1,f_2)\in L^2(\mr\times \p X)$ if
$(f_1,f_2)\in E_{\ac}(X).$ Thus the left hand side of \eqref{einrad+} is
well defined when $\Im \la=0$ 
for such initial data. We want to understand the extension
of the right hand side for this type of data, and we proceed as in \cite{gui}.
An application of Green's identity,
see for example the proof of Proposition 2.1 of \cite{gui}, gives
\begin{gather}
R(\frac{n}{2}+i\la,z,z')-R(\frac{n}{2}-i\la,z,z')
= -2i\la \int_{\p X} 
E(\frac{n}{2}+i\la,y,z)E(\frac{n}{2}-i\la,y,z')\; d\vol_{h_0}(y),\label{resid}
\end{gather}
provided $\frac{n}{2}\pm i\la$ 
 are not poles of the resolvent $R(\bullet,z,z').$

On the other hand, by  using \eqref{resid} and Stone's formula, see
 the proof of Proposition 2.2 of \cite{gui},  we deduce that the map 

\begin{gather*}
E: C_0^\infty(\intx) \longrightarrow 
C^{\infty}(\p X \times \mr_+) \\
\;\ \;\ \;\ \phi \longmapsto 
\sqrt{\frac{2}{\pi}} \int_{X} E(\frac{n}{2}+ i\la,y,z')\phi(z')\; d\vol_g(z'), 
\;\ \la >0
\end{gather*} 
induces a surjective
 isometry of the space of absolute continuity of the Laplacian
\begin{gather}
E: L^2_{\operatorname{ac}}(X) \longrightarrow 
L^2\left(\mr^+; L^2\left(\p X\right),
 \la^2 d\la \right). \label{isometry}
\end{gather}
Moreover it is an spectral representation in the sense that
\begin{gather}
E\Delta= \left(\frac{n^2}{4}+\la^2\right) E. \label{spec}
\end{gather}

A similar analysis can be carried out for the backward radiation field.
We know that
\begin{gather*}
R\left(\frac{n}{2}-i\la\right)=\left(\Delta-\la^2-\frac{n^2}{4}\right)^{-1}=
\int_{\mr} U_{-}(t,z,z') e^{-i\la t} dt, \;\ \Im \la >\frac{n}{2},
\end{gather*}
and hence 
\begin{gather}
\begin{gathered}
\mcf\left( \mcr_{-}(f_1,f_2)\right)(y,\la)=
i\la E\left(\frac{n}{2}-i\la\right)\left(f_2+\la f_1\right), \\
\Im \la>\frac{n}{2}\;\ \text{ if } \;\ (f_0,f_1) \in C_0^\infty(X), \;\ 
\Im \la \geq 0 \;\ \text{ if } \;\ 
(f_0,f_1) \in C_0^\infty(X)\cap E_{\ac}(X).
\end{gathered}\label{einrad-}
\end{gather}

Now we are ready to prove Theorem \ref{inv}. 

\begin{proof}
Using \eqref{isometry} we observe that, for $\la \in (0,\infty)$
and $f, g \in L^2_{\operatorname{ac}}(X),$
\begin{gather}
\lan f,g\ran_{L^2(X)}=
\frac{2}{\pi} \lan \la E\left(\frac{n}{2}\pm i\la\right) f, 
\la E\left(\frac{n}{2}\pm i\la\right) g\ran_{L^2(\mr_+\times \p X)}.
\label{Eq1}
\end{gather}
On the other hand, as $\overline{E\left(\frac{n}{2}+i\la\right)}=
E\left(\frac{n}{2}-i\la\right)$ when $\la\in (0,\infty),$
\begin{gather*}
||E\left(\frac{n}{2}+i\la\right)(\la f_2+\la^2 f_1)||^2_{L^2(\mr\times \p X)}=
\\
\int_{-\infty}^{\infty} E\left(\frac{n}{2}+i\la\right)(\la f_2+\la^2f_1)
 E\left(\frac{n}{2}-i\la\right)(\la \overline{f}_2+\la^2\overline{f}_1) \;\ 
d\la d\vol_{h_0}= \\
\int_0^\infty E\left(\frac{n}{2}+i\la\right)(\la f_2+\la^2f_1)
E\left(\frac{n}{2}-i\la\right)(\la \overline{f}_2+\la^2\overline{f}_1)\;\ 
d\la d\vol_{h_0}+ \\
\int_{0}^\infty E\left(\frac{n}{2}-i\la\right)(-\la f_2+\la^2f_1)
E\left(\frac{n}{2}+i\la\right)(-\la \overline{f}_2+\la^2\overline{f}_1)\;\ 
d\la d\vol_{h_0}= \\
2\int_0^\infty \left|E\left(\frac{n}{2}+i\la\right)\la f_2\right|^2 \; 
d\la d\vol_{h_0}+
2\int_0^\infty \left|E\left(\frac{n}{2}+i\la\right)\la^2 f_1\right|^2\; 
d\la d\vol_{h_0}.
\end{gather*}
Equation \eqref{Eq1} shows that
\begin{gather*}
\int_0^\infty \left|E\left(\frac{n}{2}+i\la\right)\la f_2\right|^2 \; 
d\la d\vol_{h_0}=
\frac{\pi}{2} ||f_2||_{L^2(X)}^2.
\end{gather*}
But by \eqref{spec}
\begin{gather*}
\la^2 E\left(\frac{n}{2}+i\la\right) f_1=
E\left(\frac{n}{2}+i\la\right)\left(\Delta -\frac{n^2}{4}\right)f_1,
\end{gather*}
and thus
\begin{gather*}
\int_0^\infty \left|E\left(\frac{n}{2}+i\la\right)\la^2 f_1\right|^2\; 
d\la d\vol_{h_0}=
\lan \la E\left(\frac{n}{2}+i\la\right) f_1,
\la E\left(\frac{n}{2}+i\la\right)\left(\Delta-\frac{n^2}{4}\right)
 f_1\ran_{L^2(\mr \times \p X)}.
\end{gather*}
Again by \eqref{Eq1}
\begin{gather*}
\lan \la E\left(\frac{n}{2}+i\la\right) f_1, 
\la E\left(\frac{n}{2}+i\la\right)\left(\Delta-\frac{n^2}{4}\right)
 f_1\ran_{L^2(\mr \times \p X)}=
\frac{\pi}{2} \lan f_1,\left(\Delta-\frac{n^2}{4}\right) 
f_1\ran_{L^2(X)}= \\
\frac{\pi}{2} \int_{X} \left( |d f_1|^2 -\frac{n^2}{4}|f_1|^2\right)\; 
d\vol_{g}.
\end{gather*}
So we conclude that for $(f_1,f_2)\in E_{\ac}(X),$
\begin{gather*}
||E\left(\frac{n}{2}+i\la\right)(\la f_2-\la^2 f_1)||^2_{L^2(\mr \times \p X)}=
2\pi ||(f_1,f_2)||^2_{H_E(X)}
\end{gather*}
Plancherel's theorem and \eqref{einrad+} show that $\mcr_+$ is norm preserving,
 and
therefore its  range is closed. It is clearly dense, 
otherwise there would be $\phi\in L^2(\mr\times \p X)$ orthogonal to the 
range of $\mcr,$ but by \eqref{isometry}, $\phi=0.$
This concludes the proof of Theorem \ref{inv}.
\end{proof}

 In view of equations \eqref{einrad+} and \eqref{einrad-},
 Theorem \ref{equiv} is equivalent to
\begin{prop}\label{mapscatp}
For $\la \in \mr,$ $\la\not=0,$ the scattering matrix is the unitary 
operator $$A(\la): L^2( \p X)\longrightarrow L^2( \p X)$$ 
 that satisfies 
\begin{gather}
\begin{gathered}
E\left(\frac{n}{2}+i\la\right)= A(\la)E\left(\frac{n}{2}-i\la\right).
\end{gathered}\label{mapscat}
\end{gather}
\end{prop}
\begin{proof}
We need three things:  equations  \eqref{defe}, \eqref{resid}, and
the expansion for $E\left(\frac{n}{2}+i\la,y,x',y'\right)$ as 
$x'\rightarrow 0,$ which can be found in 
Propositions 4.1 and 4.2 of \cite{jsb}, i.e
\begin{gather}
E\left(\frac{n}{2}+i\la,y,x',y'\right)=
\frac{1}{2i\la} {x'}^{-i\la}\delta(y,y')+ 
\frac{1}{2i\la} {x'}^{i\la}A(\la)(y,y') + o(x), \;\  \text{ as } 
x'\downarrow 0. \label{expjsb}
\end{gather}
Now multiplying \eqref{resid} by ${x'}^{-i\la}$ and applying \eqref{defe}
and \eqref{expjsb}, the result follows.

This is a known argument and can be found for example
in the proofs of Proposition 2.5 of \cite{gz0}, or Corollary 2.6 of 
\cite{gui}, which deal with Riemann surfaces. Here we used
the results of \cite{jsb} where needed.
\end{proof}

\section{The support theorem}\label{supth}

The main goal of this section is to prove
\begin{thm}\label{L5} 
If $f\in L_{\ac}^2(X)$  and $\mcr_{+}(0,f)(s,y)$ is supported in $s\geq s_0,$
with $s_0<<0,$
then $f$ is  supported in $\{x\geq e^{s_0}\}.$
\end{thm}

Theorem \ref{L5}  is a ``support theorem'' in the terminology of 
Helgason \cite{helgason1,helgrt} and
is a generalization of Theorem 3.13 of \cite{lptr}, which is 
Theorem \ref{L5} for the hyperbolic space $\mh^3.$
Helgason \cite{helg5,helgrt} proved such a result for functions that are 
compactly supported, but for more general symmetric spaces.

It is important to observe, as emphasized by Lax in \cite{lax}, that this
theorem does not have an analogue in (asymptotically) Euclidean space. In that 
case the function
$f$ has to be rapidly decaying at infinity-- see Theorem 2.6 and Remark 2.9 of 
\cite{helgrt}-- which would correspond
to infinite order of vanishing of $f$ at $x=0.$ Here the only requirement is 
that
$f\in L^2_{\ac}(X).$  However, in coordinates \eqref{hmet1}
the distance along a geodesic
that approaches $\p X$ perpendicularly is $-\log x,$ so the requirement that 
$f\in L^2(X)$ already imposes an  exponential of decay of $f$ near the
boundary.

The first step in the proof of Theorem \ref{L5} is 

\begin{lemma}\label{LL5} 
If $f\in L_{\ac}^2(X)$  and $\mcr_{+}(0,f)(s,y)$ is supported in $s\geq s_0,$
then $f$ is compactly supported.
\end{lemma}
\begin{proof} 
Without loss of generality we may assume that $f\in C^\infty(\intx).$ Indeed,
we recall from equations \eqref{einrad+} and \eqref{spec} 
that the Fourier transform in $s$ of the forward radiation field
is a spectral representation of $\Delta-\frac{n^2}{4}.$  Then if
$\mcf$ denotes the Fourier transform in the $s$ variable then
\begin{gather}
q(\la^2)\mcf\mcr_+(0, f)(\la,y)= \mcf\mcr_+\left(0,
q\left(\Delta-\frac{n^2}{4}\right)f\right), \;\ f \in L^2_{\ac}(X), \;\
\;\ \forall \;\ q \in \mcs(\mr).
\label{Ela}
\end{gather}

If $\phi \in C_0^\infty(\mr),$ is an even function, which we choose to be
 supported in
$(-\eps,\eps),$  then
 its Fourier transform  $\widehat{\phi}$  is also an even function,
and there exists a smooth function  $\psi_1 \in \mcs(\mr)$ such that
\begin{gather}
\widehat{\phi}(\la)=\psi_1(\la^2). \label{phipsi}
\end{gather}
Therefore by \eqref{Ela}
\begin{gather}
{\phi}*\mcr_{+}(0,f)= 
\mcr_{+}\left( 0,\psi_1
\left(\Delta-\frac{n^2}{4}\right)f\right). \label{7.2new}
\end{gather}
 As $\frac{\p^{2k}}{\p s^{2k}} \phi*\mcr_{+}(0,f)\in L^2(\mr\times \p X),$
$k=1,2,...,$ then
\begin{gather*}
\left(\Delta-\frac{n^2}{4}\right)^k
\psi_1\left(\Delta-\frac{n^2}{4}\right)f\in L^2(X).
\end{gather*}
$\Delta-\frac{n^2}{4}$ is a standard elliptic operator in the interior of $X,$
so $\psi_1\left(\Delta-\frac{n^2}{4}\right)f\in C^\infty(\intx).$ Moreover,
by elliptic regularity  for totally characteristic operators,
 see Theorem 3.8 of
\cite{ma2},
\begin{gather}
\left(xD_x, xD_{y_1}, ..., xD_{y_n}\right)^\alpha 
\psi_1\left(\Delta-\frac{n^2}{4}\right)f\in L^2(X), \;\ 
\alpha\in \mn^{n+1}.\label{regf}
\end{gather}
If $\phi\in C_0^\infty(\mr)$ is even, and $\mcr_+(0,f)(s,y)=0$ for $s<s_0,$
$\phi*\mcr_+ (0,f)(s,y)=0$ for $s < s_0-\eps.$
So from now on we will assume that $f,$ instead of
$\psi_1\left(\Delta-\frac{n^2}{4}\right)f,$
 satisfies \eqref{regf}. Since the solution to the Cauchy problem is smooth
for all (finite) times,  the solution $V$ to \eqref{eqq'} is smooth in 
$(0,T)\times (0,T)\times\p X.$  We do not know {\it a priori} that
$V$ is smooth up to $\{x'=0\}\cup \{t'=0\}.$ We proved in Theorem \ref{rad}
that this is true if $f\in C_0^\infty(\intx),$ but here $f$  is not yet
known to be
compactly supported and may be singular at $\{x'=t'=0\}.$ 

To illustrate our method, let us assume for a moment that $V$ is smooth up
to $\{x'=0\}\cup \{t'=0\}.$  Using the equation, and that $V(x',x',y)=0,$
one can prove that
if $\mcr_+(0,f)(s,y)=0$ for $s<s_0,$  then all derivatives of $V$
vanish at $\{x'=0\}\cup \{t'=0\}.$ Therefore we
can extend $V$ smoothly across  $\{t'=0\}\cup \{x'=0\}$ as $V=0.$ 
We then want to use a uniqueness theorem to conclude that $V=0$
near $\{x'=t'=0\}.$ In particular this will imply that $f(x,y)=0$ in a 
neighborhood of $\{x=0\}$ and we will be done.

The operator $P'$ in \eqref{eqq'}
extends (although not uniquely) to a neighborhood of $\{x'=t'=0\},$ however 
notice that the coefficients of the terms  of \eqref{eqq'}
involving second order derivatives 
in $y$ vanish to second order at $x'=t'=0,$  so H\"ormander's uniqueness 
theorem cannot be used to guarantee that $V=0$ near $\{t'=x'=0\}.$
 The uniqueness theorem which deals with
the Cauchy problem for second order operators with this type of degeneracy
is due to Alinhac, Theorem 1.1.2 of  \cite{Ali}. Notice that, although
$P'$ is real, it is not
of principal type at $\{x'=t'=0\},$ so the result  of Lerner and Robbiano, see
\cite{lero} or Theorem 28.4.3 of \cite{Ho}, cannot be applied either.

The principal symbol of $P'$ is
\begin{gather}
p=\sigma_2(P')=-\xi\tau - x't'h(x't',y,\eta), \;\ h(x't',y,\eta)=
\sum_{ij} h^{ij}(x't',y)\eta_i\eta_j.\label{psp'}
\end{gather}
If $H_p$ denotes the Hamilton vector field of $p$ and $\phi=x'+t',$ then
\begin{gather}
H_p=-\tau\frac{\p}{\p x'}-\xi\frac{\p}{\p t'}-x't'H_h+
(t'h+x'{t'}^2h_{1})\frac{\p}{\p \xi}+(x'h+t'{x'}^2h_{1})\frac{\p}{\p \tau},
\label{hamp}
\end{gather}
where $H_h$ denotes the Hamilton vector field of $h$ in the variables
$y$ and $\eta$ only,  and $h_{1}(x't',y,\eta)$ denotes the derivative of 
$h$ in the first variable. Hence
\begin{gather*}
H_p^2\phi=-(x'+t')\left[h(x't',y,\eta)+ x't'h_{1}(x't',y,\eta)\right],
\end{gather*}

 Since $h(0,y,\eta)$ is non-degenerate, it follows that for
$x',t'$ small
\begin{gather*}
-\phi H_p^2\phi \geq \ha \phi^2 h.
\end{gather*}
 We can therefore apply Theorem 1.1.2 of \cite{Ali} to $P'$ and 
$\phi$ (with  $\Lambda_S=\emptyset,$ and $a=C=0$ as in the statement of
Theorem 1.1.2 of \cite{Ali})
to guarantee that $V=0$ in a neighborhood of $\{x'=t'=0\}.$ 
See figure \ref{enes1}. In particular this  shows
that $f(x,y)=0$ near $x=0,$ so  $f\in C_0^\infty(\intx),$ as we want to prove.
 
We emphasize that the key reason
this works is the fact that the extension of $V$ would be
 supported in the wedge $\{x'\geq 0, t'\geq0\},$ so the intersection of
the support of $V$ and $\{\phi=0\}$ is compact, as required by Theorem 1.1.2
of \cite{Ali}.

The main difficulty to apply this method to prove Theorem \ref{L5} is that
$V$ is not known to be smooth
up to $\{x'=0\}\cup \{t'=0\},$ and Theorem 1.1.2 of \cite{Ali}
is proved for smooth functions.  We have to show that our assumptions
guarantee enough regularity of $V$ to make Alinhac's argument work.
Fortunately the operator $P'$ in 
\eqref{eqq'}  is much simpler than
the general case of \cite{Ali}, and the necessary Carleman estimates  are
relatively simple to obtain.  

We will work with $W,$ which is defined
in \eqref{defw}, instead of $V.$ The advantage is that $W$ satisfies
\eqref{eqq'W} which has no first order derivatives in $x'$ or $t'.$

The proof is  divided in two steps: \newline
{\it Step 1:} We will use that $W(x',x',y)=0,$ $x'>0,$
and $\mcr_+(0,f)=0,$ to show that
$W$ can be extended as a (distribution) solution to \eqref{eqq'W} in a 
neighborhood of
$\{x'=t'=0\}$ vanishing when either $x'<0$ or $t'<0.$  \\
{\it Step 2:} We use the method of proof of Theorem 1.1.2 of \cite{Ali}
to show that $W=0$ in a neighborhood of $\{x'=t'=0\}.$  

\begin{figure}[int1]
\epsfxsize= 4.0in
\centerline{\epsffile{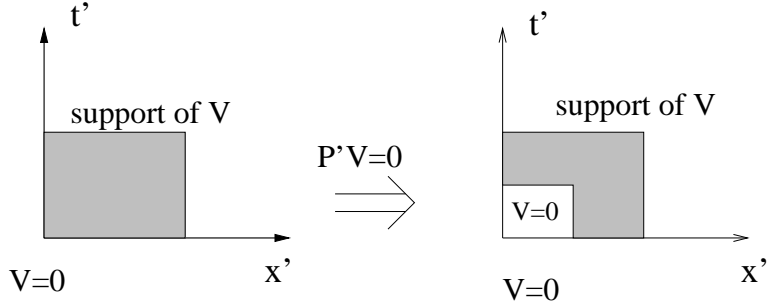}}
\caption{Unique continuation for $V.$ $P'$ is degenerate at 
$x'=t'=0,$ but $V$ is supported in the wedge $\{x'\geq 0, \; t'\geq 0\}$ }
\label{enes1}
\end{figure}

\noindent{ \it Proof of Step 1.} 
We will show that the conditions 
$W(x',x',y)=0,$ $x'>0,$ and $\mcr(0,f)=0$ for $s<s_0,$ imply that
 $W$ has an extension
$\tW$ satisfying
\begin{gather}
\tW \in H^{2k}
\left((-T,T)\times(-T,T); H^{-2k}(\p X)\right), \;\
k\in \mn \label{genreg}
\end{gather}
and
\begin{gather}
\begin{gathered}
\left(\frac{\p^2}{\p x' \p t'}+ x't'\Delta_h +x't'B(x't',y,\frac{\p}{\p y})+
C(x't',y)\right)\tW=0 \text{ in } \;\ \widetilde{\Omega} \\
\tW(x',t',y)=W(x',t',y) \;\ \text{ if } x'>0, \;\ t'>0 \\
\tW(x',t',y)=0 \;\ \text{ if } x'<0, \;\ t'<0. \\
\end{gathered}\label{eqext}
\end{gather}
Here we have extended the operators $\Delta_h,$ $B$ and the coefficient $C$ 
smoothly across $\{x'=0\}\cup \{t'=0\}.$

It was proved in Lemma \ref{enerest}
 that when $f\in L^2(X),$ and $T$ is small, the solution $W$ of
\eqref{eqq'W} satisfies 
\begin{gather}
\begin{gathered}
W,\; {x'}^\ha D_{x'} W ,\;  {t'}^\ha D_{t'} W \in 
L^2\left([0,T]\times [0,T] \times \p X)\right),\\
[x't'(x'+t')]^\ha D_{y_j}W \in L^2\left([0,T]\times [0,T] \times \p X\right), \;\
1\leq j \leq n.
\end{gathered}\label{regw}
\end{gather}
We emphasize that, due to the smoothness of $f,$ $W$ is smooth in 
the region $\{x'>0, t'>0\}.$

We will show \eqref{genreg} for $k=1.$
For that, let $\phi\in C^\infty(\p X)$  and
\begin{gather}
G(x',t')=\int_{\p X} W(x',t',y)\phi(y)\; dy. \label{defg}
\end{gather}

Using equation \eqref{eqq'W} and \eqref{regw} we deduce that
\begin{gather}
\frac{\p^2  W(x',t',y)}{\p x'\p t'} \in 
L^2\left([0,T]\times [0,T]; H^{-1}(\p X)\right).\label{dxdt}
\end{gather}
Differentiating \eqref{eqq'W} in $x'$ we have
\begin{gather}
\begin{gathered}
\frac{\p^2}{\p {x'}^2} \frac{\p}{\p t'} W = \\
-t'\Delta_h W -x't' Q W - t' \Delta_h x'\frac{\p}{\p x'}W
-t' B W - x't' B_1 W - t' B x'\frac{\p}{\p x'}W - C_1 W-
C \frac{\p}{\p x'}W,
\end{gathered}\label{difw'}
\end{gather}
where $Q$ and $B_1$ are second and first order operators respectively,
involving $y$ derivatives only. Using \eqref{regw}
we find that
\begin{gather}
{x'}^{\ha}\frac{\p^2}{\p {x'}^2}\frac{\p}{\p t'} W  \in
L^2\left([0,T]\times [0,T]; H^{-2}(\p X)\right). \label{dxdt'}
\end{gather}
Since $W(t',t',y)=0,$ $t'>0,$  equation \eqref{eqq'W} shows that 
$\frac{\p}{\p x'}\frac{\p}{\p t'} W(t',t',y)=0,$ $t'>0$ and thus
we can write
\begin{gather*}
\frac{\p^2}{\p x' \p t'} G(x',t')=
-\int_{x'}^{t'} \frac{\p^2}{\p s^2}\frac{\p}{\p t'} G(s,t')\;\ ds, \;\
\text{ if } t'>x'>0.
\end{gather*}
Therefore
\begin{gather*}
\left|\frac{\p}{\p x'}\frac{\p}{\p t'} G(x',t')\right|^2=
\left|\int_{x'}^{t'} 
\frac{\p^2}{\p s^2}\frac{\p}{\p t'} G(s,t')\;\ ds\right|^2 \leq
\log\left(\frac{t'}{x'}\right) 
\int_{x'}^{t'}s\left| 
\frac{\p^2}{\p s^2}\frac{\p}{\p t'} G(s,t')\right|^2\;\ ds.
\end{gather*}
So from \eqref{dxdt'} we obtain for $T$ small and $\del<\ha,$
\begin{gather*}
\int_0^{T} \int_0^{t'} {x'}^{-2\del}
\left|\frac{\p}{\p x'}\frac{\p}{\p t'} G(x',t')\right|^2 \;\ dx'\; dt' 
\leq \\ \int_{0}^{T}\int_{0}^{t'} \frac{s^{1-2\del}}{(1-2\del)}\left(
(2-2\del)|\log s|+ \frac{1}{1-2\del}\right)s\left| 
\frac{\p^2}{\p s^2}\frac{\p}{\p t'} G(s,t')\right|^2 \; ds \; dt'< \\
||{x'}^{\ha}\frac{\p^2}{\p {x'}^2}
\frac{\p}{\p t'} W||^2_{L^2\left([0,T]\times [0,T]; H^{-2}(\p X)\right)}<\infty.
\end{gather*}
So 
\begin{gather}
{x'}^{-\del}\frac{\p}{\p x'}\frac{\p}{\p t'} G
\in L^2\left(\{t'\geq x'> 0\}\right), \;\
\delta<\ha,\label{dxdt''}
\end{gather}
with uniform bound up to $x'=0.$ This together with \eqref{dxdt'}  imply that
\begin{gather}
\frac{\p}{\p x'}
\left({x'}^{\gamma}\frac{\p}{\p x'}\frac{\p}{\p t'} W\right)
\in L^2\left(\{t'\geq x'> 0\}; H^{-2}(\p X)\right), \;\ \gamma>\ha.
\label{regamma}
\end{gather}
Thus the restriction 
$\left({x'}^{\gamma} \frac{\p}{\p x'}\frac{\p}{\p t'} W\right)(0,t',y),$ 
$t'>0,$  is well defined for any $\gamma>\ha.$ In particular this shows
that
\begin{gather}
\left( x'\frac{\p}{\p x'} W\right)(0,t',y)=0, \;\ \text{ in } 
(0,T)\times \p X. \label{rf=01}
\end{gather}
In view of \eqref{dxdt}, the restriction
$\frac{\p}{\p t'} W (0,t',y)$  is well defined.

Recall that in these coordinates, and in terms of $W,$ the radiation field is 
given by
\begin{gather}
\mcr_+(0,f)= \ha |h|^\oq (0,y) \left.\left( t'\frac{\p W}{\p t'} - 
x'\frac{\p W}{\p x'}\right)\right|_{x'=0}=0.\label{rf=0}
\end{gather}

So \eqref{rf=0} and \eqref{rf=01} show that
\begin{gather}
\frac{\p}{\p t'} W (0,t',y)=0, \;\ \text{ in } \;\ 
(0,T)\times \p X. \label{rf=02}
\end{gather}
Again applying the regularity of $W$ given by \eqref{dxdt} and using
\eqref{rf=02}, we have for $t'>0,$

\begin{gather*}
\left|\frac{\p G}{\p t'}(t',t')\right|^2 \leq
\left|\int_{0}^{t'}\frac{\p^2 G}{\p s\p t'}(s,t') \; ds\right|^2\leq
\frac{{t'}^{\delta+1}}{\delta+1}
\int_0^{t'} s^{-\del} \left|\frac{\p^2 G}{\p s\p t'}(s,t')\right|^2 \; ds, 
\;\ \del<\ha,
\end{gather*} 
and we find that
\begin{gather}
\int_{0}^{T} {t'}^{-1-\del}\left|\frac{\p G}{\p t'}(t',t')\right|^2
<\infty, \;\ \del<\tha.
\label{t-1}
\end{gather}
In particular this shows that 
\begin{gather}
\frac{\p W}{\p x'}(x',x',y), \;\ \frac{\p W}{\p t'}(x',x',y) \in
 L^2( [0,T] ; H^{-2}(\p X)). \label{t-11}
\end{gather}
We can then write
\begin{gather*}
\frac{\p W}{\p t'}(x',t',y)=\frac{\p W}{\p t'}(x',x',y)+
\int_{x'}^{t'}\frac{\p^2 W}{\p s \p t'}(s,t',y) \; ds \text{ and } \\
\frac{\p W}{\p x'}(x',t',y)=\frac{\p W}{\p x'}(x',x',y)+
\int_{x'}^{t'}\frac{\p^2 W}{\p x' \p s}(x',s,y) \; ds
\end{gather*}
and  use \eqref{t-11} to show that
\begin{gather}
\frac{\p W}{\p t'}(x',t',y) \;\ \text{ and } \;\
\frac{\p W}{\p x'}(x',t',y)  \in
L^2\left( [0,T]\times [0,T]\; H^{-2}(\p X)\right). \label{pxpt-reg}
\end{gather}
Therefore $W(0,t',y),$ $t'\in [0,T],$
is well defined and,  by \eqref{rf=02}, $W(0,t',y)=F(y).$
Since $W(x',x',y)=0,$ and we
want to construct a smooth extension to $G,$ we must have
\begin{gather}
W(0,t',y)=0, \;\ t'>0.\label{vanofw}
\end{gather}

Substituting \eqref{pxpt-reg} back in \eqref{difw'} we deduce that
\begin{gather}
\begin{gathered}
\frac{\p}{\p x'} \frac{\p}{\p x'}\frac{\p}{\p t'} W\in 
L^2\left([0,T]\times [0,T]; H^{-2}(\p X)\right) 
\end{gathered}\label{dxdt'1}
\end{gather}
Therefore 
$\frac{\p}{\p x'}\frac{\p}{\p t'} W(0,t',y)$ is well defined
and from \eqref{eqq'W} and \eqref{vanofw}
\begin{gather}
 \frac{\p}{\p x'}\frac{\p}{\p t'} W (0,t',y)=0, \;\ t'>0.\label{vanofpdw}
\end{gather}
But we need  to prove that
$\left. \frac{\p W}{\p x'} \right|_{x'=0}=0$ on $[0,T].$
 To do this 
we take the derivative of \eqref{eqq'W} with respect to $t'.$ 
\begin{gather}
\begin{gathered}
 \frac{\p}{\p x'}\frac{\p^2}{\p {t'}^2} W = \\
-x'\Delta_h W -x't' Q_1 W - x' \Delta_h t'\frac{\p}{\p t'}W
-x' B W - x't' B_2 W - x' B t'\frac{\p}{\p t'}W - C_2 W-
C \frac{\p}{\p t'}W,
\end{gathered}\label{difwt'}
\end{gather}
where as in \eqref{difw'}, $Q_1$ and $B_2$ are second and first order 
operators respectively,
involving $y$ derivatives only. Using \eqref{regw} and \eqref{pxpt-reg}
we find that
\begin{gather}
\begin{gathered} 
 \frac{\p}{\p x'}\frac{\p^2}{\p {t'}^2} W \in 
L^2\left( [0,T]\times [0,T]\; H^{-2}(\p X)\right).
\end{gathered}\label{pxpt2-reg}
\end{gather}
Therefore \eqref{pxpt-reg}, \eqref{dxdt'1} and \eqref{pxpt2-reg}
imply that
\begin{gather}
W \;\ \text{ and } \;\ \frac{\p}{\p x'}\frac{\p}{\p t'} W  \in 
H^1\left([0,T]\times [0,T]; H^{-2}(\p X)\right).\label{h1reg}
\end{gather}
From \eqref{difwt'} and \eqref{pxpt2-reg} we deduce that
 $\frac{\p^2 W}{\p {t'}^2}(0,t')|_{x'=0}$ is well defined
and  by \eqref{rf=02}, $\frac{\p^2 W}{\p {t'}^2}(0,t')|_{x'=0}=0.$ By symmetry
 $\frac{\p^2 W}{\p {x'}^2}(x',0)=0.$ Thus we can write
\begin{gather*}
\frac{\p^2 W}{\p {t'}^2}(x',t',y)=
\int_0^{x'} \frac{\p}{\p s}\frac{\p^2}{\p {t'}^2} W \; ds \\
\frac{\p^2 W}{\p {x'}^2}(x',t',y)=
\int_0^{t'} \frac{\p}{\p s}\frac{\p^2}{\p {x'}^2} W \; ds. 
\end{gather*} 
and then conclude that 
\begin{gather*}
\frac{\p^2 W}{\p {t'}^2}(x',t',y), \;\ \frac{\p^2 W}{\p {x'}^2}(x',t',y)
 \in L^2\left([0,T]\times [0,T]; H^{-2}(\p X)\right).
\end{gather*}
From \eqref{regw} and \eqref{pxpt-reg} we have
\begin{gather}
W\in H^2\left([0,T]\times [0,T]; H^{-2}(\p X)\right). \label{2ndreg}
\end{gather}
Moreover $W=0$ on $\{x'=0\}\cup \{t'=0\}.$

Now we write
\begin{gather*}
\frac{\p^2 G}{\p {t'}^2}(t',t')=\int_0^{t'} 
\frac{\p}{\p s}\frac{\p^2 G}{\p {t'}^2}(s,t')\;\ ds
\end{gather*}
As in the proof of \eqref{t-1}  we use \eqref{pxpt2-reg} to show that
$G(t',t')\in C^1([0,T])$ and, by symmetry,
\begin{gather*}
\frac{\p G}{\p t'}(0,0)=\frac{\p G}{\p x'}(0,0)=0.
\end{gather*}
From \eqref{h1reg} $\frac{\p W}{\p t'}(0,t',y)\in H^1([0,\eps]; H^{-2}(\p X)),$
From this and \eqref{rf=02} we deduce that 
\begin{gather}
\frac{\p  W}{\p t'}(0,t',y)=0, \;\
\frac{\p  W}{\p x'}(x',0,y)=0, \;\ \text{ in } \;\ [0,T] \times \p X.
\label{vanpw}
\end{gather}
This and the regularity of $W$ given by \eqref{2ndreg} are
 enough to guarantee that if
 $\tW(x',t',y)$ is the extension of $W$ to 
$(-T,T)\times(-T,T)\times\p X,$ 
with $\tW=0$ in 
$\{x'<0\}\cup \{t'<0\},$  then  it satisfies \eqref{eqext} and 
\begin{gather}
\tW \in  H^2\left([-T,T]\times [-T,T]; H^{-2}(\p X)\right).
\label{k=2}
\end{gather}

To prove \eqref{genreg} we substitute the regularity \eqref{k=2}
 back into \eqref{eqq'W} and iterate following the argument above.
This also shows that all derivatives of $\tW$ vanish at 
$\{x'=0\}\cup \{t'=0\}.$ This ends the proof of step 1.

\noindent{\it Proof of step 2:}  We want to show that \eqref{genreg} 
and \eqref{regw} are enough to apply Alinhac's theorem.

First we make the change of variables
\begin{gather*}
t'=\mu+\nu, \;\ x'=\mu-\nu.
\end{gather*}
From \eqref{eqq'W} 
$\tU(u,v,y)=\tW(\mu+\nu,\mu-\nu,y)$ satisfies
\begin{gather*}
P\tU=0, \;\ \text{ where } \\
P=\frac{\p^2 }{\p \mu^2}-\frac{\p^2 }{\p \nu^2}+
4(\mu^2-\nu^2)\sum_{i,j=1}^n h^{ij}(\mu^2-\nu^2,y)\frac{\p }{\p y_i}\frac{\p }{\p y_j}
+(\mu^2-\nu^2)\sum_{j=1}^n b_j(\mu^2-\nu^2,y)\frac{\p }{\p y_j}+ \\
 C(\mu^2-\nu^2,y)
\end{gather*}
As in section 4.1 of \cite{Ali}, let 
\begin{gather*}
Q=-\frac{\p^2 }{\p \nu^2}+
4(\mu^2-\nu^2)\sum_{i,j=1}^n h^{ij}(\mu^2-\nu^2,y)\frac{\p }{\p y_i}\frac{\p }{\p y_j},
\;\  \\
L=\frac{\p^2 }{\p \mu^2}+Q \;\ \text{ and } \;\
A=L+\frac{\gamma(\gamma-1)}{\mu^2} \;\ \text{ with } \;\ \gamma\geq \frac{1}{10}.
\end{gather*}
We could have taken $\gamma>\gamma_0>0,$ but to avoid one more parameter
we set $\gamma_0=\frac{1}{10},$ which is good enough for our purposes.
Let $Y_0\subset \p X $ be a neighborhood $y_0\in \p X$
and let
$C_0^{*\infty}([0,T)_\mu\times (-T,T)_\nu\times Y_0)$ denote the
set of functions in 
$C_0^{\infty}([0,T)_\mu\times (-T,T)_\nu\times Y_0)$ 
which vanish to infinite order at $\{\mu=0\}.$
Integration by parts gives that for every
$v \in C_0^{*\infty}([0,T)_\mu\times (-T,T)_\nu\times Y_0),$ 
\begin{gather}
\begin{gathered}
||\mu^{-\gamma}L \mu^{\gamma} v||^2=-2\gamma\Re\left(\mu Q_\mu v, 
\frac{v}{\mu^2}\right)
+4(\gamma^2+\gamma)\left|\left|\mu^{-1}v_\mu\right|\right|^2+||Av||^2+
(4\gamma^3-4\gamma^2-6\gamma)\left|\left|\mu^{-2}v\right|\right|^2\\
-2\gamma\left(Av,\mu^{-2}v\right).
\end{gathered}\label{eqa}
\end{gather}

Here $Q_\mu$ is the operator obtained by differentiating the coefficients of
 $Q$ with respect to $\mu,$ and
$||,||$ denotes the norm in $L^2((-T,T)\times(-T,T)\times Y_0).$  This is 
equation 4.11 of \cite{Ali} with $a=0$ and $a'=1.$

Hardy's inequality, see for example Lemma 5.3.1 of \cite{davies}, gives
that  $||\mu^{-1}v_\mu||^2\geq \frac{9}{4}||\mu^{-2}v||^2.$ Using this and
applying Cauchy-Schwartz inequality to the last 
term of \eqref{eqa} gives
\begin{gather}
||\mu^{-\gamma}L \mu^{\gamma} v||^2\geq -2\gamma\Re\left(\mu Q_\mu v, \mu^{-2}v\right)
+2\gamma^2\left|\left|\mu^{-1}v_\mu\right|\right|^2+
(4\gamma^3+\ha \gamma^2+2\gamma)\left|\left|\mu^{-2}v\right|\right|^2. \label{ineql}
\end{gather}
The derivate of $Q$ with respect to $\mu$ is
\begin{gather*}
Q_\mu=8\mu\sum_{i,j=1}^n h^{ij}(\mu^2-\nu^2,y)\frac{\p }{\p y_i}\frac{\p }{\p y_j}
+8\mu(\mu^2-\nu^2)
\sum_{i,j=1}^n h_{1}^{ij}(\mu^2-\nu^2,y)\frac{\p }{\p y_i}\frac{\p }{\p y_j},
\end{gather*}
where $h_{1}$ denotes  the derivative of $h$ in the first variable.
Since $h(0,y)$ is non-degenerate, then for $T$ small enough, 
\begin{gather}
\Re\left(\mu Q_\mu v, \mu^{-2} v\right)\leq 
-C ||\nabla_y v||^2. \label{ineqqt}
\end{gather}
Letting $\mu^\gamma v=w,$ using \eqref{ineql} and \eqref{ineqqt}
we arrive at the following Carleman estimate for the operator $L.$
\begin{gather*}
||\mu^{-\gamma}L w||^2\geq C\gamma||\mu^{-\gamma}\nabla_y w||^2+
\gamma^2\left|\left|\mu^{-1}(\mu^{-\gamma} w)_\mu\right|\right|^2+
\gamma^3\left|\left|\mu^{-2-\gamma}w\right|\right|^2, \\
\text{ for all } \;\ 
w \in C_0^{*\infty}([0,T)_\mu\times (-T,T)_\nu\times Y_0).
\end{gather*}
By the definition of $P$ we have that, for $T$ small, there exists $K$ 
depending on $C$, $T,$  but not on $\gamma$ such that
\begin{gather}
\begin{gathered}
K||\mu^{-\gamma}P w||^2\geq \gamma||\mu^{-\gamma}\nabla_y w||^2+
\gamma^2\left|\left|\mu^{-1}(\mu^{-\gamma} w)_\mu\right|\right|^2+
\gamma^3\left|\left|\mu^{-2-\gamma}w\right|\right|^2, \\
\text{ for all }
\;\ w \in C_0^{*\infty}([0,T)_\mu\times (-T,T)_\nu\times Y_0).
\end{gathered}\label{carle1}
\end{gather}

Let 
$$\chi(y) \in C_0^\infty(\mrn), \;\ \chi(0)=1, \;\
\int_{\mrn} \chi(y)  \; dy=1
\text{ and } \chi_\delta(y)=\delta^{-n}\chi\left(\frac{y}{\delta}\right),
\;\ \delta>0,$$
 and for $y_0\in Z_0\subset\subset Y_0,$ let
 $$\phi \in C_0^\infty((-T,T)_\mu\times (-T,T)_\nu\times Y_0),
\text{ with } \phi(\mu,\nu,y)=1 \text{ in }
(-\frac{T}{2},\frac{T}{2})_\mu
\times (-\frac{T}{2},\frac{T}{2})_\nu\times Z_0.$$ 
Let
$v_\delta=\chi_\del*'\phi \tU,$ where as in \cite{Hoo}, $*'$ means that
the convolution is taken in the variable $y$ only. In view of
\eqref{genreg} and the fact that all derivatives of $\tW$ vanish at
$\{x'=0\}\cup \{t'=0\},$
$v_\del\in C_0^{*\infty}([0,T)_\mu\times (-T,T)_\nu\times {Y}_0),$ 
for  $\del$ small enough. Therefore we can apply
\eqref{carle1} to $v_\del.$ Thus 
\begin{gather}
K||\mu^{-\gamma}P v_\del||^2\geq 
\gamma||\mu^{-\gamma}\nabla_y v_\del||^2+
\gamma^2\left|\left|\mu^{-1}(\mu^{-\gamma} v_\del)_\mu\right|\right|^2+
\gamma^3\left|\left|\mu^{-2-\gamma}v_\del\right|\right|^2. \label{conv1}
\end{gather}
Now we want to take the limit of \eqref{conv1} as $\delta\rightarrow 0.$
First we observe that, since $\tU$ is supported in $\mu\geq |\nu|$ and  
satisfies 
$P\tU=0,$ and
$\phi=1$ near $\mu=\nu=0,$ it follows that $\mu>\ha T$ on the support of
$P\phi \tU.$ From \eqref{regw} we have

\begin{gather}
||\mu^{-\gamma}P \phi \tU||^2 
\leq (\ha T)^{-\gamma}||P \phi \tU||^2 <\infty, \;\ \forall \;\ 
\gamma>\frac{1}{10}.
\label{tgamma}
\end{gather} 

On the other hand we write
\begin{gather*}
P v_{\del}=\chi_\del*'(P\phi \tU) + [P,\chi_\del *']\phi\tU
\end{gather*}
and observe that
\begin{gather}
\begin{gathered}
\left[P,\chi_\del*'\right]=  \\
\left[\sum_{i,j}
\left((\mu^2-\nu^2)\left(h^{ij}(\mu^2-\nu^2,y)\frac{\p}{\p y_i}\frac{\p}{\p y_j}
+b_{j}(\mu^2-\nu^2,y)\frac{\p}{\p y_j}\right) +
(\mu^2-\nu^2)B_1(\mu^2-\nu^2,y)\right), \chi_\del*'\right]. 
\end{gathered}\label{conv}
\end{gather}
Equation 2.4.18 and Theorem 2.4.3 of of \cite{Hoo} show
 that, for $u\in H^{s-1}(\mrn)$
and $a\in C_0^\infty(\mrn),$
\begin{gather*}
||a(u*\chi_\del)-(au)*\chi_\del||_{H^s}\leq C||u||_{H^{s-1}}, \;\
\text{ and } \\
a(u*\chi_\del)-(au)*\chi_\del \rightarrow 0, \;\ \text{ in } \;\ H^s(\mrn)
\;\ \text{ as } \;\ \del\rightarrow 0.
\end{gather*}
Applying this to our case, it follows from \eqref{conv} that, if

\begin{gather}
\begin{gathered}
||\mu^{-\gamma}(\mu^2-\nu^2)\nabla_y\left(\phi\tU\right)|| <\infty, \;
\text{ and } \;
||\mu^{-\gamma}(\mu^2-\nu^2)\phi\tU||_{L^2([0,T] \times [0,T]; H^{-1}(Y_0))}<\infty,
\end{gathered}\label{condition}
\end{gather}
then
\begin{gather*}
\lim_{\delta \rightarrow 0} ||\mu^{-\gamma}[P,\chi_\del'](\phi\tU)|| =0.
\label{tgamma1}
\end{gather*}

Thus, if $\phi\tU$ satisfies \eqref{condition} for a certain 
$\gamma>\frac{1}{10},$ it follows from \eqref{tgamma} 
that we can take the limit as $\delta\rightarrow 0$ in
\eqref{conv1} and therefore
\begin{gather}
\begin{gathered}
\infty> K||\mu^{-\gamma}P\phi \tU||^2\geq 
\gamma||\mu^{-\gamma}\nabla_y \phi\tU||^2+
\gamma^2\left|\left|\mu^{-1}(\mu^{-\gamma} \phi\tU)_\mu\right|\right|^2+
\gamma^3\left|\left|\mu^{-2-\gamma}\phi \tU\right|\right|^2.
\end{gathered} \label{conv2}
\end{gather}
We know from \eqref{regw} that
$||[\mu(\mu^2-\nu^2)]^\ha\nabla_y \phi\tU||<\infty$ and that 
$||\phi\tU||<\infty.$
Since $\phi\tU$ is  supported in $\mu \geq  |\nu|,$  we have
$$||\mu^{-\ha}(\mu^2-\nu^2) \nabla_y \tU||<\infty \text{ and }
||\mu^{-\ha}(\mu^2-\nu^2)\phi\tU||_{L^2([0,T] \times [0,T]; 
H^{-1}(Y_0))}<\infty.$$
 So applying \eqref{conv2} with 
$\gamma=\ha$ we obtain
$$||\mu^{-\ha}\nabla_y \phi\tU||<\infty, \;\ 
||\mu^{-\frac{5}{2}} \phi\tU||<\infty.$$
Thus we can apply \eqref{conv2} with $\gamma=2.$
 Repeating this argument we obtain
\begin{gather*}
\mu^{-\gamma}\nabla_y\left(\phi\tU\right), \;\ \mu^{-1}(\mu^{-\gamma}\phi\tU)_t, \;\
\mu^{-\gamma-2}\phi\tU\in L^2((-T,T)\times (-T,T)\times \p M),\;\
\forall \;\ \gamma>\frac{1}{10}.
\end{gather*}

Equation \eqref{conv2} shows in particular  that
\begin{gather*}
\gamma^3||\mu^{-\gamma-2}\tU||_{L^2((-T/2,T/2)\times 
(T,T)\times Z_0)} =\gamma^3||\mu^{-\gamma-2}\phi\tU||_{L^2((-T/2,T/2)\times 
(T,T)\times Y_0)}  \leq \\
\leq \gamma^3||\mu^{-\gamma-2}\phi\tU||_{L^2((-T,T)\times 
(T,T)\times Y_0)} 
\leq 
K||\mu^{-\gamma}P\phi\tU ||_{L^2((-T,T)\times (\eps,\eps)\times Y_0)}
\end{gather*}
Using that $\mu>\ha T$ on the support of $P\phi\tU$ we arrive at
\begin{gather*}
\left(\frac{T}{2}\right)^{-\gamma-2}
 \gamma^3
||\tU||_{L^2((-T/2,T/2)\times (-T/2,T/2)\times Z_0)}
\leq K \left(\frac{T}{2}\right)^{-\gamma}
|| P\phi\tU ||_{L^2((-T,T)\times (-T,T)\times Y_0)}.
\end{gather*}
So in particular
\begin{gather}
 \gamma^3
||\tU||_{L^2((-T/2,T/2)\times (-T/2,T/2)\times Z_0)}
\leq K \left(\frac{T}{2}\right)^{2}
|| P\phi\tU ||_{L^2((-T,T)\times (-T,T)\times Y_0)}.
\label{fincar}
\end{gather}
By \eqref{regw} and the definition of $P,$ the right hand side of 
\eqref{fincar} is finite.
Then letting $\gamma\rightarrow \infty$ gives $\tU=0$ in 
$(-T/2,T/2)\times (-T/2,T/2)\times Z_0.$
Since $\p X$ is compact, this argument can be applied to a partition of unity
consisting of finitely many functions to deduce that there exists
$\eps>0$ such that $\tU=0$ in
$(-\eps,\eps)\times (-\eps,\eps)\times \p X.$ In particular $f(x,y)=0$ if
$x<\sqrt{\eps}.$  This ends the proof of the Lemma.
\end{proof}

To conclude the proof of Theorem \ref{L5} we need
\begin{lemma}\label{spcomp} If $f \in L_{\comp}^2(\intx)$ and
$\mcr_+(0,f)$ is supported in $s>s_0,$ with $s_0<<0,$ then
$f$ is supported in $x>e^{s_0}.$
\end{lemma}
\begin{proof}  
Let 
$\phi \in C_0^\infty(\mr)$ be even,  supported  in $(-1,1),$ and 
$\phi_{\epsilon}(s)=\eps^{-n}\phi(s/\eps).$ Then $\phi_\eps*\mcr_+(0,f)$
is supported in $s>s_0-\eps.$ Let $\psi_1\in \mcs(\mr)$ and
$\psi_{1,\eps}\in \mcs(\mr)$ be the functions obtained from $\phi$ and
$\phi_\eps$ by \eqref{phipsi}. It follows from Lemma \ref{LL5} that
$\psi_{1,\eps} (\Delta-\frac{n^2}{4})f$ is compactly supported.  If we prove
that $\psi_{1,\eps}(\Delta-\frac{n^2}{4})f$ is supported in $x\geq e^{s_0-\eps}$
then by letting $\eps\rightarrow 0$ we have proved Lemma \ref{spcomp}.
So we may assume that $f\in C_0^\infty(\intx),$ and let us
say that $f(x,y)=0$ for $x<x_0.$ Then  \eqref{werad1}, which
is due to the finite speed of propagation, gives that

\begin{gather}
V(x',t',y)=0 \;\ \text{ if }  \;\ t'<\sqrt{x_0}. \label{suppv'}
\end{gather}

Moreover we know from Theorem \ref{rad} that $V$ is smooth up to $x'=0.$
Since by assumption $V(0,t',y)=0$ for $t'<t_0=\exp(\frac{s_0}{2}),$ then,
as proved above, $V$ extends to a smooth
solution to  \eqref{eqq'} for $x'<0$
which vanishes in $\{x'<0, \;\ t'<t_0\}.$ 

The principal symbol of $P'$  is given by \eqref{psp'}, so  it is easily seen
that the level surfaces of $\{\phi=- t'-x'\}$ are not characteristic
for $P'.$ The Hamiltonian of $p',$ $H_{p'},$ was computed in \eqref{hamp}.
Then
\begin{gather*}
p'(0,t',y,\xi,\tau,\eta)=-\xi\tau, \;\
\left(H_{p'}\phi\right)(0,t',y,\xi,\tau,\eta)=\tau+\xi, \;\ 
\text{ and }\\
\left(H_{p'}^2\phi\right)(0,t',y,\xi,\tau,\eta)=t'h(0,y,\eta),
\end{gather*}
and therefore, for $t'>0,$
\begin{gather*}
\text{if} \;\  p'(0,t',y,\xi,\tau,\eta)=
\left(H_{p'}\phi\right)(0,t',y,\xi,\tau,\eta)=0, \;\
\text{ with }  (\xi,\tau,\eta)\not=0 \;\  \text{ then } \\
\left(H_{p'}^2\phi\right)(0,t',y,\xi,\tau,\eta)>0.
\end{gather*}
This, according to section 28.4 of \cite{Ho}, implies that the level surfaces
of $\phi$ are strongly pseudo-convex with respect
to $P',$ as long as $t'>0.$
  So applying H\"ormander's theorem, Theorem 28.3.4 of
\cite{Ho}, to \eqref{eqq'} and  the surface $\{\phi=-\sqrt{x_0}\}$
we find that $V(x',t',y)=0$ near
$(0, \sqrt{x_0},y).$  Now we can repeat this argument to show that
$V=0$  in a neighborhood of the segment $[0,t_0].$ And we can 
proceed in the same way to conclude that in fact $V=0$ in
$\{(x',t',y): 0<x'<\sqrt{x_0}, \; 0<t'<t_0\}.$ By the symmetry of $V,$ it 
follows that
$V=0$ in $\{(x',t',y): 0<t'<\sqrt{x_0}, \; 0<x'<t_0\},$
see figure \ref{fig2n}.  This is the point where we use that the initial data
is $(0,f).$

\begin{figure}[int1]
\epsfxsize= 4.0in
\centerline{\epsffile{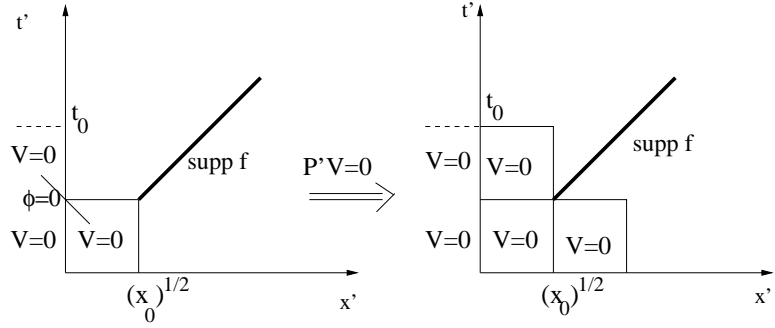}}
\caption{Unique continuation for $V.$ Away from $x'=t'=0,$
the level surfaces of $\phi$ are strongly pseudo-convex.}
\label{fig2n}
\end{figure}

Now going back to the variables $(x,t,y)$ which are given by
$x=x't',$ $t=\log t'-\log x',$ this implies that
the solution $u(x,t,y)$ of \eqref{we2} with initial data $(0,f)$ 
satisfies
\begin{gather}
u(x,t,y)=0 \;\ \text{ for }
\{(x,t,y):0<x<x_0, \;\ \log x_0-s_0<t< s_0-\log x_0 \}. \label{u=0}
\end{gather}

One particular case  of Tataru's theorem \cite{tat}, see also
\cite{horuni,ro,rozu,tat1}, states that if $u(t,z)$ is in
 $H^1_{\loc}$ and satisfies 
\begin{gather*}
\begin{gathered}
\left(D_t^2-\Delta-\frac{n^2}{4}\right)u(t,z)=0 \\
u(t,z)=0 \;\ -T<t<T, \text{ and } \;\ d(z,z_0)<\delta, \;\ \delta>0, \\
\end{gathered}
\end{gather*}
where $d(z,z_0)$ is the distance between $z$ and $z_0$ with respect to the 
metric $g,$ then 
\begin{gather}
u(t,z)=0  \text{ for } |t|+d(z,z_0)<T. \label{uniqt1}
\end{gather}

We then deduce from \eqref{u=0} and \eqref{uniqt1}
that $u(0,z)=D_tu(0,z)=0$ on the set of points $z=(x,y)$ at a distance
less than $s_0-\log x_0$ from the surface $\{x=x_0\}.$ In these coordinates
the distance is given by $\log x-\log x_0.$ So $f(x,y)=D_tu(0,x,y)=0$
for $\log x<s_0.$ This proves Lemma \ref{spcomp}.
\end{proof}

\section{The Inverse Problem}\label{prometdet}

We consider the inverse problem of determining the manifold and the 
metric from the scattering matrix $A(\la),$ {\it at all energies} 
$\la \in \mr\setminus 0.$ We remark that, as $A(\la)$ has a meromorphic
continuation to $\mc\setminus \frac{i}{2} \mn,$ see \cite{jsb}, 
this assumption implies that $A_1(\la)=A_2(\la)$
for $\la \in \mc\setminus\frac{i}{2}\mn.$

We are fixing a defining function $x$ of $\p X$ for which \eqref{hmet1} 
holds near $M$ and using it  in  \eqref{scatmat}
to define $A(\la).$  According to Lemma 2.1 of
\cite{graham},  $x$ is uniquely determined near $M$ by the choice of the
conformal representative $h_0=x^2 g|_{\p X}.$ So to define the scattering 
matrix we fix a  conformal representative 
$h_0.$ On the other hand, we recall from Theorem 1.1 and 
Corollary 1.1 of \cite{jsb} that $A(\la),$  for $\la$ fixed, determines
the conformal representative $h_0.$ So if we say that two asymptotically
hyperbolic metrics $g_1$ and $g_2$
on $X$ have the same scattering matrix, as defined by 
\eqref{scatmat}, it is implicit that we are fixing the same
 conformal representative for $x_j^2g_j|_{\p X},$ $j=1,2.$

 We will prove
\begin{thm}\label{metdet} Let $(X_1,g_1)$ and $(X_2,g_2)$
be asymptotically hyperbolic manifolds which
have the same boundary $\p X_1=\p X_2 =M.$ 
  Let $x_j\in C^\infty(X_j)$ be a defining function
of $M=\p X_j,$ $j=1,2,$ for which \eqref{hmet1} holds, and let 
$A_j(\la),$ $j=1,2,$ $\la \in \mr\setminus 0,$ be the corresponding
scattering matrices
defined in \eqref{scatmat} in terms of $x_j.$
Suppose that $A_1(\la)=A_2(\la)$ for every $\la \in \mr\setminus 0.$
Then there exists a diffeomorphism
$\Psi:X_1\longrightarrow X_2,$ smooth up to $M,$  such that
\begin{gather}
\Psi=\Id \text{ at }  M \;\ \text{ and } \;\ \Psi^*g_2= g_1.
\label{metdeteq}
\end{gather}
\end{thm}

As mentioned in the introduction, the proof is an application of the control 
method of Belishev \cite{be}, see also \cite{bk1,lassas,kaku,kaku1}.
We will also use a result, which is an application of this 
method, and is due to Katchalov and Kurylev \cite{kaku,kaku1}

First we construct
a diffeomorphism between neighborhoods of the boundary
that realizes \eqref{metdeteq} and later show that it can be extended to
a diffeomorphism between  the two manifolds.

We recall that Lemma 2.1 of \cite{graham} states that fixed a conformal
representative $h_0$ there exists a unique defining function $x_i,$ in a
neighborhood of $M$  for which  \eqref{hmet1} holds.
 However $x_i$ can be continued to $X_i,$ 
although not uniquely.
As explained in the paragraph before Theorem \ref{metdet}, it is implicit
that we are fixing a conformal representative 
\begin{gather}
h_0=x_1^2 g_1|_{M}=x_2^2 g_2|_{M}. \label{confrep}
\end{gather}

Note that, as observed in the  paragraph following the proof Lemma 2.1 
of \cite{graham}, for 
$\eps>0$ small, a defining function $x_i$ for which \eqref{hmet1} holds
near $M$  gives an identification of
$[0,\eps) \times \p M$ with a collar neighborhood 
$U_{i,\eps}\subset \overline{X_i}$ of $M$ by
\begin{gather}
\Psi_{i,\eps} : [0,\eps)\times M \longrightarrow U_{i,\eps}, \label{maps}
\end{gather}
where $\Psi_{i,\eps}(x,y)$ is the point obtained by flowing the integral 
curve of 
$\nabla_{x_i^2 g_i} x_i$ emanating from $y$ by $x$ units of time.
So  $x$ is the arc-length along the geodesics normal to the boundary $M$
with respect to the metrics $x_i^2 g_i,$ or the distance from the point
$(x,y)$ to $M$  with respect to this metric.  We can pick $\eps$ small such 
that these maps are diffeomorphisms. 
Then $x$  is a smooth defining function of $M,$
\begin{gather}
\begin{gathered}
x:  [0,\eps)\times M \longrightarrow \mr, \;\  \;\ 
x=x_i+o(x_i), \;\ \text{ and } \\
g_j=\frac{dx^2}{x^2}+\frac{h_j}{x^2}, \;\ j=1,2.
\end{gathered}\label{norm12}
\end{gather}

In fact for each $y\in M$ there exists $\eps_i(y)>0$ such that
for $x<\eps_i(y)$ the distance between $(x,y)$ and $M,$ with respect to
the metric $x_i^2g_i,$ is equal to $x.$ So 
 $\Psi_{i,\eps}$ extends to a diffeomorphism
\begin{gather}
\Psi_{i} : [0,\eps_i(y))\times M \longrightarrow X_i\setminus \Gamma_i,
\label{cut-loc}
\end{gather}
where $\Gamma_i$ is the cut-locus of $X_i$ with respect to $x_i^2g_i.$ It is
known that the set $\Gamma_i$ is a closed subset of measure zero
and $\Gamma_i\cap M=\emptyset.$ The
number $\eps$ in \eqref{norm12} is less than the smallest of the distances
 between $\Gamma_i$ and $M.$

 We will fix one such function throughout this section and we will prove

\begin{prop}\label{determ} Let $(X_1,g_1)$ and $(X_2,g_2)$ be asymptotically 
hyperbolic manifolds  satisfying the hypotheses of Theorem \ref{metdet}.
Then there exists $\eps>0$ such that in the product decomposition 
 $X\sim [0,\eps)\times \p X$ where \eqref{norm12} holds, $h_1=h_2.$
\end{prop}

The fact that
the metrics are equal in these coordinates imply that
$$\Psi_{1,\eps}^* \left(\left. g_1\right|_{U_{1,\eps}}\right) =
\Psi_{2,\eps}^* \left(\left. g_2\right|_{U_{2,\eps}}\right),$$
with $\Psi_{i,\eps}$ defined in \eqref{maps},
and therefore
$$\left[\Psi_{2,\eps}\circ \Psi_{1,\eps}^{-1}\right]^* 
\left(\left. g_2\right|_{U_{2,\eps}}\right)= 
\left. g_1\right|_{U_{1,\eps}}.$$

This  gives  a diffeomorphism  between neighborhoods of  the boundary
satisfying \eqref{metdeteq}.

\subsection{Preliminaries} Here we define some spaces that will
be useful in the proof of Proposition  \ref{determ}.

 Notice that
 if $u$ is a solution of \eqref{we2} with initial 
data $(0,f),$ then
$u(t,z)=-u(-t,z),$ $t\in \mr.$  Similarly, if $u$ is a solution of 
\eqref{we2} with initial data $(f,0),$ then
$u(t,z)=u(-t,z),$ $t\in \mr.$  This implies that 
\begin{gather}
\begin{gathered}
\mcr_+(0,f)(-s,y)=\mcr_-(0,f)(s,y), \;\
\text{ and } \;\ \mcr_+(f,0)(-s,y)=-\mcr_-(f,0)(s,y).
\end{gathered}\label{mc1}
\end{gather}

\begin{prop}\label{defm} Let $(X,g)$ be an asymptotically hyperbolic manifold
and let $x$ be a defining function  of $\p X$ for which \eqref{hmet1} holds.
For $F\in L^2(\mr\times \p X),$ let $F^*(s,y)=F(-s,y).$
Let $\mcr_{\pm}$ denote  the radiation fields defined with respect  to
$x,$ and let $\mcs=\mcr_+\mcr_{-}^{-1}$  be the scattering operator.
Let
\begin{gather}
 \mcm^f=\{ F\in L^2(\mr\times \p X): \;\ F=\mcs F^*\} \;\ \text{ and } \;\
\mcm^b=\{ F\in L^2(\mr\times \p X): \;\ F^*=\mcs F\}, \label{mcmdef0}
\end{gather}
where $f$  and $b$ stand for forward and backward. Then
\begin{gather}
\begin{gathered}
 \mcm^f=\{\mcr_{+}(0,f): \;\ f\in L_{\ac}^2(X)\} \;\ \text{ and } \;\
\mcm^b=\{\mcr_{-}(0,f): \;\ f\in L_{\ac}^2(X)\}
\end{gathered}\label{mcmdef}
\end{gather}
\end{prop}
\begin{proof} We will prove the first  equality. The proof of the second one 
is identical.  If
$F=\mcr_+(0,f),$ then according to \eqref{mc1}, $F^*=\mcr_-(0,f),$ and so
\begin{gather*}
\mcs F^*=\mcr_+\mcr_-^{-1}F^*=\mcr_+(0,f)=F.
\end{gather*}

Notice also that if
$F=\mcr_+(g,0),$ then according to \eqref{mc1}, $F^*=-\mcr_-(g,0),$ and so
\begin{gather*}
\mcs F^*=\mcr_+\mcr_-^{-1}F^*=-\mcr_+(g,0)=-F.
\end{gather*}
Conversely, let $F\in L^2(\mr \times \p X).$
We know from Theorem \ref{inv} that there exists $(f_1,f_2)\in E_{\ac}$
such that $F=\mcr_+(f_1,f_2).$  Let $F_1=\mcr_+(f_1,0)$ and 
$F_2=\mcr_+(0,f_2).$ If $F=\mcs F^*,$ then
\begin{gather*}
F=F_1+F_2=\mcs F_1^*+ \mcs F_2^*.
\end{gather*}
In view of the  discussion above, $\mcs F_2^*=F_2$ and
$\mcs F_1^*=-F_1.$ 
So $F_1+F_2=-F_1+F_2,$ and thus $F_1=0.$ By uniqueness, $f_1=0.$
\end{proof}
As ranges of  bounded operators, $\mcm^f$ and $\mcm^b$ are closed.  By 
Proposition \ref{defm} and Theorem \ref{equiv}
they are also defined in terms of the scattering matrix and the space
$L^2(\mr \times \p X),$ which depends on the volume element of the metric 
$h(0,y,dy).$  However, as discussed above,
the scattering matrix determines the
tensor $h(0,x,dy),$ so $\mcm^f$  and $\mcm^b$ are
 determined by $A(\la),$ $\la \in \mr\setminus 0.$
Thus we have

\begin{cor}\label{mcmd} Let $(X,g)$ be an asymptotically hyperbolic 
manifold.
Let $x$ be a defining function  of $\p X$ for which \eqref{hmet1} holds
and let $\mcr_{\pm}$ denote  the radiation fields defined with respect  to
$x.$ Then the spaces $\mcm^f$  and $\mcm^b$ defined in \eqref{mcmdef0}
are closed subspaces of
$L^2(\mr\times \p X),$ and so are Hilbert spaces with the inherited norm.
Moreover, $\mcm^f$ and $\mcm^b$ are
 determined by the scattering matrix $A(\la),$ for all
$\la \in \mr\setminus 0.$
\end{cor}

We will need
\begin{lemma}\label{charact} Let $(X,g)$ be an asymptotically
hyperbolic manifold. Let $\mcm^f$ and $\mcm^b$ be defined in \eqref{mcmdef0}.
 For $x_1\in (0,\eps),$ we have
\begin{gather}
\begin{gathered}
 \mcm^f(x_1)\defi
\{ F\in \mcm^f, \;\ F(s,y)=0 \;\ \text{ if } \;\ s <\log x_1\}=
\{\mcr_+(0,f): \;\ f \in L^2_{\ac}(X), \;\ f=0 \;\
\text{ if } \;\ x<x_1\}, \\
 \mcm^b(x_1)\defi\{ F\in \mcm^b, \;\
F(s,y)=0 \;\ \text{ if } \;\ s>-\log x_1\}=
\{\mcr_-(0,f): \;\ f \in L^2_{\ac}(X), \;\ f=0 \;\
\text{ if } \;\ x<x_1\}.
\end{gathered}\label{mcmdefx1}
\end{gather}

\end{lemma}
\begin{proof}
Finite speed of propagation guarantees that if 
$f \in L^2_{\ac}(X)$  and $f=0$
 in $\{ x<x_1\},$ then $F=\mcr_+(0,f)\in\mcm^f(x_1)$ and
$F=\mcr_-(0,f)\in \mcm^b(x_1).$  On the other hand,
if $F\in \mcm^f(x_1),$ in particular
$F\in\mcm^f,$ and thus by Proposition \ref{defm}, $F=\mcr_+(0,f)$
with $f\in L^2_{\ac}(X).$
Since $F(s,y)=0$ if $s<\log x_1,$ it follows from Theorem \ref{L5}
that
$f=0$ in $x<x_1.$

If $F\in \mcm^b(x_1),$ then $F=\mcr_-(0,g).$ So $F^*=\mcr_+(0,g)$ and
$F^*=0$ for $s<\log x_1.$ Therefore Theorem \ref{L5} gives that
$g=0$ in $x<x_1.$ 
\end{proof}

We emphasize that the proof of Lemma \ref{charact}
required the full power Theorem \ref{L5} and this is where the support theorem
enters in the study of the inverse problem.

In what follows, for $F\in \mcm^f,$ or $F\in \mcm^b,$
 we will use, by an abuse of notation,
$\mcr_\pm^{-1} F= f,$ and $\mcr_\pm f =F$  instead
of $\mcr_\pm^{-1}F=(0,f)$ and $\mcr_\pm(0,f)=F,$ respectively.

\subsection{The case of no eigenvalues}
To better explain our methods, we will first consider  the case where
the manifolds have no  eigenvalues.  
In this  particular case,
 Proposition \ref{determ} is an easy consequence of
\begin{prop}\label{Schk}
 Let $g_1$ and $g_2$ be asymptotically hyperbolic metrics
satisfying the hypotheses of Theorem \ref{metdet}. Moreover assume that
$\Delta_{g_i},$ $i=1,2,$ have no eigenvalues.
 Let $\mcr_{j,\pm},$ $j=1,2,$ denote
the corresponding forward or backward radiation fields defined in 
coordinates in which
\eqref{norm12} holds.  Then there  exists $\eps>0$ such that 
\begin{gather}
\begin{gathered}
|h_1|^{\oq}(x,y)
\mcr_{1,-}^{-1}F(x,y)=|h_2|^{\oq}(x,y)\mcr_{2,-}^{-1}F(x,y),
\;\ (x,y)\in (0,\eps)\times M, \;\ \forall \;\ F\in \mcm^b, \\
|h_1|^{\oq}(x,y)\mcr_{1,+}^{-1}F(x,y)=|h_2|^{\oq}(x,y)\mcr_{2,+}^{-1}F(x,y),
\;\ (x,y)\in (0,\eps)\times M, \;\ \forall \;\ F\in \mcm^f. \\.
\end{gathered}\label{equality}
\end{gather}
\end{prop}
Indeed, suppose Proposition \ref{Schk} has been proved. We will
prove Proposition \ref{determ}.
\begin{proof} Let $x$ be such that \eqref{norm12} holds.
We know that for any  $F\in \mcm^b,$ 
\begin{gather}
 \mcr_{j,-}^{-1}\left(\frac{\p^2}{\p s^2}F\right)=
\left(\Delta_{g_j}-\frac{n^2}{4}\right)\mcr_{j,-}^{-1} F. \label{f-1}
\end{gather}
So the first equation in \eqref{equality} implies that for
$(x,y)\in  (0,\eps)\times M,$ and any $F\in \mcm^b,$
\begin{gather}
\begin{gathered}
|h_1|^{\oq}(x,y)\mcr_{1,-}^{-1}F(x,y)=
|h_2|^{\oq}(x,y)\mcr_{2,-}^{-1}F(x,y), 
\text{ and } \\
|h_1|^{\oq}(x,y)
\left(\Delta_{g_1}-\frac{n^2}{4}\right)\mcr_{1,-}^{-1} F(x,y)=
|h_2|^{\oq}(x,y)\left(\Delta_{g_2}-\frac{n^2}{4}\right)\mcr_{2,-}^{-1} F(x,y).
\end{gathered}\label{key02}
\end{gather}
Set $\mcr_{1,-}^{-1}F=f.$ Since $F$ is arbitrary and 
the metrics have no eigenvalues, equations \eqref{key02} give, 
in particular,
\begin{gather}
\begin{gathered}
|h_1|^{\oq}(x,y)
\left(\Delta_{g_1}-\frac{n^2}{4}\right)f(x,y)=
|h_2|^{\oq}(x,y)
\left(\Delta_{g_2}-\frac{n^2}{4}\right)
\frac{|h_1|^{\oq}(x,y)}{|h_2|^{\oq}(x,y)}f(x,y), \\ 
\forall \; f \in C_0^\infty\left((0,\eps)\times M \right).
\end{gathered} \label{key2}
\end{gather} 
Therefore the operators on both sides of \eqref{key2} are equal. In particular
the coefficients of the principal parts of
$\Delta_{g_1}$ are equal to those
of $\Delta_{g_2},$ and hence the tensors $h_1$ and $h_2$ from \eqref{norm12}
 are equal. This proves Proposition \ref{determ}.
\end{proof}

We begin the proof of Proposition \ref{Schk} with

\begin{lemma}\label{proj} Let $(X,g)$ be an asymptotically
hyperbolic manifold such that $\Delta_g$  has no eigenvalues.
Let $x$ be such that \eqref{hmet1} holds in
$(0,\eps)\times M.$ For $x_1 \in (0,\eps),$ let 
$\mcp_{x_1}^b$ denote the orthogonal projector
\begin{gather*}
\mcp_{x_1}^b : \mcm^b \longrightarrow \mcm^{b}( x_1),
\end{gather*}
and let $\chi_{x_1}$ be the characteristic function of
 the set $\{x\geq x_1\}.$ 
Then  for every $f\in L_{\ac}^2(X)=L^2(X),$
\begin{gather*}
\mcp_{x_1}^b\mcr_{-}(0,f)=\mcr_{-}(0,\chi_{x_1}f)\in \mcm^{b}( x_1).
\end{gather*}
\end{lemma}
\begin{proof} 
Since $\mcp_{x_1}^b$ is a projector, then for all $G\in \mcm^b( x_1)$
\begin{gather*}
\lan \mcp_{x_1}^b \mcr_-(0,f), G\ran_{L^2(\mr\times \p X)}=
\lan \mcr_-(0,f), G\ran_{L^2(\mr\times \p X)}.
 \end{gather*}
In particular, since $L^2_{\ac}(X)=L^2(X),$ then 
for all $g\in C_0^\infty(\intx)$ supported in $\{x>x_1\},$
\begin{gather*}
\lan \mcp_{x_1}^b \mcr_-(0,f), \mcr_-(0,g)\ran_{L^2(\mr\times \p X)}=
\lan f, g \ran_{L^2(X)}.
 \end{gather*}

On the other hand, by Lemma  \ref{charact} there exists $f_{x_1}\in L^2(X)$ 
supported in $\{x\geq x_1\}$
 such that $\mcp_{x_1}^b\mcr_-(0,f)=\mcr_-(0,f_{x_1}).$ Therefore, for all
$g\in C_0^\infty(\intx)$ supported in $\{x>x_1\},$
\begin{gather*}
\lan \mcp_{x_1}^b \mcr_-(0,f), \mcr_-(0,g)\ran_{L^2(\mr\times \p X)}=
\lan \mcr_-(0,f_{x_1}), \mcr_-(0,g)\ran_{L^2(\mr\times \p X)}=
\lan f_{x_1},g\ran_{L^2(X)}.
\end{gather*}
Thus
\begin{gather*}
\lan f_{x_1},g\ran_{L^2(X)}= \lan f,g\ran_{L^2(X)}, \;\
\forall \;\ g\in L^2(X) \;\ \text{ supported in } \;\ \{x\geq x_1\},
\end{gather*}
and so $f_{x_1}=\chi_{x_1}f.$

This ends the proof of the Proposition.
\end{proof}

\begin{remark}\label{indep} Since the spaces $\mcm^f,$ $\mcm^b,$
 $\mcm^{f}(x_1)$ and 
$\mcm^{b}(x_1)$ are determined by the scattering matrix, so are the projectors
$\mcp_{x_1}^{f}$ and $\mcp_{x_1}^{b}.$
\end{remark}

We use this and properties of
 propagation  of  singularities for solutions  to the wave equation
to prove
\begin{lemma}\label{jump} Let $(X,g)$ be an asymptotically hyperbolic
manifold
which has no  eigenvalues. Let $x$ be  a boundary defining function
such that \eqref{hmet1} holds.  Let $\mcr_\pm$ be the radiation fields
with respect to $x.$
Then, there exists $\eps>0,$ such that for any $x_1\in(0,\eps/4),$ 
any  $F\in \mcm^b,$  any 
$\phi \in C_0^\infty(\mr)$ even, and  the corresponding
$\psi_1$ be defined by \eqref{phipsi}, 
\begin{gather}
\begin{gathered}
\mcr_{+}\mcr_{-}^{-1}\left(\mcp_{x_1}^b G\right)(s,y)
= \ha x_1^{-\frac{n}{2}}w(x_1,y)\frac{|h|^{\oq}(x_1,y)}{|h|^{\oq}(0,y)}
(s-\log x_1)_{+}^0 +\text{ smoother terms},\\
\text{ if  } \; s< \log  x_1+ \log 4, \;\ \text{ where } \;\ 
G=\phi*F, \;\ \text{ and } \;\ w=\mcr_-^{-1}G=
\psi_1\left(\Delta_{g}-\frac{n^2}{4}\right)\mcr_-^{-1}F.
\end{gathered} \label{key}
\end{gather}
Notice that by \eqref{7.2new} $G\in \mcm^b,$ and that, according to Remark 
\ref{indep},
the left hand side of \eqref{key} is determined  by the scattering
matrix of $g.$
\end{lemma}
\begin{proof}
We choose $\eps>0$  so that the maps $\Psi_{j,\eps}$ defined
by \eqref{maps} are diffeomorphisms. That is,  the distance from $\Gamma,$
the cut-locus, to
$M$ is greater than $\eps,$ and take $x_1<\frac{\eps}{4}.$
From Lemma \ref{proj} we know that
\begin{gather*}
(0,\chi_{x_1}w)=\mcr_{-}^{-1}\left(\mcp_{x_1}^{b}G\right),
\end{gather*}
with $\chi_{x_1}$ being  the characteristic function of  the set
$\{x\geq x_1\},$  which is the set of points whose distance  to
$M$ is greater than or equal to $x_1,$
and hence
\begin{gather*}
\mcr_{+}(0,\chi_{x_1}w)=
\mcr_{+}\mcr_{-}^{-1}\left(\mcp_{x_1}^{b}G\right).
\end{gather*}

So we want to analyze $\mcr_{+}(0,\chi_{x_1}w).$ 
Recall from the definition of
the forward radiation field that this amounts to finding the solution
$u$ to \eqref{we2}  with initial data $(0,\chi_{x_1}w),$ then taking
$v(s,x,y)=x^{-\frac{n}{2}} u(s-\log x, x,  y)$
and restricting $\frac{\p v}{\p s}$ to $x=0.$  We remark that, although this 
is the definition
of the forward radiation field for $C_0^\infty(\intx)$ data, it follows from
the discussion in section
\ref{supth}  that this also holds for initial data in $L^2_{\ac}(X)=L^2(X),$ 
 see equations \eqref{dxdt'} and \eqref{rf=01}. 
We are  concerned with the restriction of $\frac{\p v}{\p s}$ to 
$x=0,$ so we will only consider the behavior of $v$ for $s>\log x.$ 

The initial data $\chi_{x_1}w$ has a conormal singularity at
$S=\{s=\log x, \; x=x_1\},$ therefore the wavefront set of $v$ will be
 contained
in the flow-out of  $N^*S\cap \Sigma$, where $\Sigma$
is the characteristic variety  of $P,$ which is defined in  \eqref{neqq}.
The principal symbol of $P$ is
$p=-2\xi\sigma-x\xi^2-xh(x,y,\eta),$  $\sigma$ is the dual to $s.$ So the 
null bicharacteristics satisfy
\begin{gather*}
\dot{x}=-2\sigma-2x\xi, \;\ \dot{s}=-2\xi, \;\
 \dot{y}=-x\frac{\p  h}{\p \eta}, \;\
\dot{\xi}=\xi^2+h+x\frac{\p  h}{\p x}, \;\  
\dot{\eta}=x\frac{\p  h}{\p y}, \;\ 
\dot{\sigma}=0, \\
 x(0)=x_1, \;\ s(0)=\log x_1,  \;\
 y(0)=y_0, \;\ \xi(0)=\xi_0,  \;\  \eta_0=\eta_0,\;\
\sigma(0)=\sigma_0, \;\   2\sigma\xi+x\xi^2+xh(x,y,\eta)=0.
\end{gather*}
Note that since $p=0$ is satisfied for the  solutions  to this system, 
and $\sigma=\sigma_0\not=0$, we must have
\begin{gather*}
\xi=\frac{1}{x}\left(-\sigma_0 \pm (\sigma_0^2-x^2h)^\ha\right).
\end{gather*}
Since we  are concerned with the forward singularities, we must have
 $\sigma_0>0.$ Then it follows that $\xi\leq 0,$ and  thus $s$ is 
non-decreasing.

We  analyze the singularities that start over
$N^*S=\{x=x_1, \; s=\log x_1,  \; \eta_0=0\}.$
Using $x$ as a parameter, we find  two families of curves
\begin{gather*}
s=\log x_1, \; y=y_0, \;\  \xi=0, \;\ \eta=\eta_0, \; \sigma=\sigma_0 \;\ \text{ if }
 \; \xi_0=0, \;\ \text{ and  } \\
s=2\log x-\log  x_1, \; y=y_0, \;\ \eta=\eta_0, \;\ \sigma=\sigma_0, \;\
\xi=\frac{\xi_0 x_1}{x}, \;   \text{ if } \; \xi_0<0, \; \text{ and } \;
2\sigma_0+x_1\xi_0=0.
\end{gather*}
These curves make up the two characteristic surfaces of $P$
emanating from
$\{s=\log x,\; x=x_1\},$ 
which are $\Sigma^+=\{s=\log x_1\}$ and
$\Sigma^-=\{s=2\log x-\log x_1\}.$ So, in $s>\log x,$ $v$ is singular along
$\Sigma^+$ and $\Sigma^-.$  See figure \ref{propsing}.

\begin{figure}[int1]
\epsfxsize= 2.5in
\centerline{\epsffile{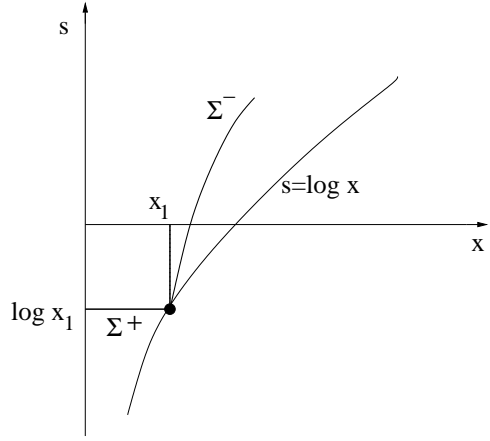}}
\caption{ The singularities of $v$ in $s>\log x.$}
\label{propsing}
\end{figure}

A bicharacteristic in $\Sigma^+,$ let us say, starting over $(x_1,y_0)$
will hit $\mr \times M$ at $(\log x_1, y_0),$ see figure \ref{propsing}.
 On the other  hand, since
$s>\log x,$ note that $\Sigma^-$ consists 
of null bicharacteristics
which go into the interior of $X,$ so they might also intersect the 
boundary. They could, however, become trapped in the interior and not reach
$M.$ Take one such curve, $\gamma,$ starting over $(x_1,y),$ $y\in M,$ which
intersects the boundary at another point $y'\in  M,$ $y'\not=y.$
The projection of $\gamma$ onto $X$ is a geodesic with 
respect to the metric $x^2g,$ which connects the points $(x_1,y)$ and 
$y'\in M$ and which is orthogonal to $\{x=x_1\}$ at $(x_1,y).$
To reach the point $y'$  this  geodesic has to pass through the cut-locus
$\Gamma$ and again reach the surface $\{x=x_1\}.$  Since $\eps$ is less than 
the distance from $\Gamma$ to $M$ and $x_1<\eps/4,$  and along $\gamma,$
for $x<\eps,$
$t=s-\log x=\log x-\log x_1,$ 
the geodesic reaches $\Gamma$ for  $t>\log \eps-\log (\eps/4)=\log 4.$
Thus it reaches
$\{x=x_1\}$ for $t=T>\log 4.$   Now we analyze the flow of $H_p$ starting
at the point $q$ where $\gamma$ intersects $\{x=x_1\}.$ 
 The surface $t=T$ becomes
$s=T+\log x$ and we  think of this singularity as starting
at $\{s=T+\log x, \;\ x=x_1\}.$
But, as observed above, $s$ is increasing along the part of $\gamma$
connecting $q$ to a point over $M.$ Then this singularity
will hit the boundary at $y'$ for $s>T+\log x_1>\log 4+\log x_1.$ 
Hence the singularity of $\mcr_+(0,\chi_{x_1} w)$ at
$s=\log x_1$ comes only from $\Sigma^+,$ and these can be computed explicitly.
Moreover $\mcr_+(0,\chi_{x_1} w)$ is supported in $s\geq \log x_1$ and
$$\mcr_+(0,\chi_{x_1} w) \in 
C^\infty\left((\log x_1, \log x_1+\log 4)\times M\right).$$

As commonly done in this type of problem,
we will find a conormal expansion for $v$ along $\{s=\log x_1\}$ and
$\{s=2\log x-\log x_1\}.$  We will construct
\begin{gather*}
V^+(s,x,y)\sim \sum_{j=1}^{\infty} v_j^+(x,y)\left(e^s- x_1\right)_{+}^j \;\
\text{ and } \;\
V^-(s,x,y)\sim 
\sum_{j=1}^\infty v_j^{-}(x,y) \left( x^2 - x_1e^s\right)_{+}^{j}, 
\end{gather*}
that is $V^+$ and $V^-$ are  the asymptotic sums of those series, 
such that, $V^\pm$  is supported in a neighborhood of $\Sigma^\pm$
(this can be arranged, as in the proof of Borel's lemma),
and in the set 
$$D_\eps=\{ \log x \leq s < \log x_1  + \log 4, \; 0 \leq x \leq \eps/2\},$$
\begin{gather*}
PV^{\pm} \in C^\infty(D_\eps), \;\
v(s,x,y) - V^+(s,x,y) - V^-(s,x,y) \in C^\infty(D_\eps),
\end{gather*}
where $P$ is defined in \eqref{neqq}. Moreover  they satisfy
\begin{gather}
\begin{gathered}
V^+(\log x,x,y)+ V^-(\log x, x,y)\in C_0^\infty, \\
\frac{\p V^+}{\p s}(\log x,x,y)+\frac{\p V^-}{\p s}(\log x,x,y)-
 x^{-\frac{n}{2}} w(x,y)(x-x_1)_+^0 \in C_0^\infty.
\end{gathered}\label{v+v-}
\end{gather}
Once this is accomplished, we then  have
\begin{gather*}
P(v-V^+-V^-) \in C^\infty(D_\eps), \\
(v-V^+-V^-)\restr_{s=\log  x} \in C_0^\infty(\intx), \;\ \;\
\frac{\p (v-V^+-V^-)}{\p  s}\restr_{s=\log  x} \in C_0^\infty(\intx).
\end{gather*}
Then the energy  estimates from section \ref{prad}, adapted to the case where
the right hand side is not equal to zero, but is smooth up to $M,$ gives  that 
$\left(v-V^+ -V^-\right)\restr_{x=0}$ is $C^\infty$ in 
$(-\infty,\log x_1 +\log 4) \times M.$

We  have
\begin{gather}
\begin{gathered}
V^+(\log x,x,y)+ V^-(\log x, x,y)
 \sim \sum_{j=1}^{\infty} \left( v_j^+(x,y)+x^jv_j^-(x,y)\right)
(x-x_1)_+^j\sim \sum_{j=1}^{\infty} Z_j(x_1,y)(x-x_1)_+^j, \\ \text{ with } \;\
Z_j(x_1,y)= \sum_{k+m=j, \; k\geq 0, m\geq 1} \frac{1}{k!} 
 \frac{\p^k}{\p x^k}\left(v_m^+ + x^mv_m^- \right)(x_1,y)
\end{gathered}\label{inc0}
\end{gather}
and
\begin{gather*}
\frac{\p V^+}{\p s}(s,x,y)+\frac{\p V^-}{\p s}(s,x,y)\sim \sum_{j=1}^\infty j e^s v_j^+(x,y)(e^s-x_1)_+^{j-1}
-\sum_{j=1}^\infty j x_1e^s v_j^-(x,y)(x^2-x_1e^s)_+^{j-1}.
\end{gather*}
So, in particular, when $s=\log x,$
\begin{gather}
\begin{gathered}
\frac{\p V^+}{\p s}(\log x,x,y)+\frac{\p V^-}{\p s}(\log x,x,y) \sim 
\sum_{j=1}^\infty \left( j x v_j^+(x,y)- j x_1 x^j v_j^-(x,y)\right)
(x-x_1)_+^{j-1} \sim \\
\sum_{j=1}^{\infty} M_j(x_1,y)(x-x_1)_+^{j-1}, \;\ \text{ where } \\
M_j(x_1,y)=\sum_{m+k=j, \; m\geq 1, k\geq 0} \frac{m}{k!} 
\frac{\p^k}{\p x^k}\left( xv_m^+-x_1x^m v_m^-\right)(x_1,y).
\end{gathered} \label{dvds1}
\end{gather}
Since
\begin{gather}
\begin{gathered}
 x^{-\frac{n}{2}} w(x,y)(x-x_1)_+^0 \sim
\sum_{j=1}^{\infty} \frac{1}{(j-1)!} 
\frac{\p^{j-1}}{\p x^{j-1}}
\left( x^{-\frac{n}{2}} w\right)(x_1,y)(x-x_1)_+^{j-1},
\end{gathered}\label{dvds}
\end{gather}
 we conclude  from  \eqref{inc0}, \eqref{dvds1}  and \eqref{dvds}
that conditions \eqref{v+v-} translate into
\begin{gather}
\begin{gathered}
Z_j(x_1,y)= \sum_{k+m=j, \; k\geq 0, m\geq 1} \frac{1}{k!} 
 \frac{\p^k}{\p x^k}\left(v_m^+ + x^mv_m^- \right)(x_1,y)=0, \;\ j\in \mn, \\
M_j(x_1,y)=\sum_{m+k=j, \; m\geq 1, k\geq 0} \frac{m}{k!} 
\frac{\p^k}{\p x^k}\left( xv_m^+-x_1x^m v_m^-\right)(x_1,y)=
\frac{1}{(j-1)!} \frac{\p^{j-1}}{\p x^{j-1}}\left(x^{-\frac{n}{2}} w\right)(x_1,y), \;\ j \in \mn.
\end{gathered}\label{dvds2}
\end{gather}
In particular,  the terms with $j=1$ satisfy
\begin{gather*}
v_1^+(x_1,y)+ x_1v_1^-(x_1,y)=0, \;\ \text{ and } \;\
x_1v_1^+(x_1,y)- x_1^2v_1^-(x_1,y)=x_1^{-\frac{n}{2}} w(x_1,y).
\end{gather*}
Hence
\begin{gather}
v_1^+(x_1,y)= -x_1 v_1^-(x_1,y)=
\ha x_1^{-\frac{n}{2}-1} w(x_1,y). \label{inicond}
\end{gather}

 Let
$$Q=P-2\frac{\p}{\p x}\frac{\p}{\p s}-A\frac{\p}{\p s}.$$
That is, $Q$ is the part of $P$ that does not have derivatives in $s.$
As discussed above, the singularity of $v$ in $\Sigma^-$  will
not hit the boundary for $s\in (\log x_1, \log x_1+\log 4).$ 
 Since we are interested on the singularity of 
$\mcr_{+}(0,\chi_{x_1}w)$ at 
$s=\log x_1,$ we will restrict our computations to $V^+,$ but keeping
in mind that, at least for $x<\eps$  and $s<\log  x_1+\log  4,$
similar computations also hold for $V^-.$  

Since $A$ and the coefficients of $Q$ do not depend on $s,$
\begin{gather}
\begin{gathered}
PV^+(s,x,y)\sim  x_1(2\frac{\p}{\p x}v_1^++Av_1^+)(e^s-x_1)_+^0 + \\
 \sum_{j=1}^{\infty}\left( (j+1)x_1\left(2\frac{\p}{\p x}+A\right)v_{j+1}^+
+\left(2j\frac{\p}{\p x}+j A+ Q\right)v_j^+\right)(e^s-x_1)_+^{j}\sim 0 \;\ \text{ in } \;\
x<x_1.
\end{gathered}\label{exppv}
\end{gather}
Since we want $PV^+(s,x,y) \in C^\infty$ in $\{ x< x_1\}$, all the coefficients
in \eqref{exppv} must be equal to zero,  and so we get  the transport equations
 for $v_j^+,$ $j=1,2,...,$ which
are
\begin{gather}
\begin{gathered}
x_1\left(2\frac{\p}{\p x}+A\right)v_1^+(x,y)=0, \;\ x< x_1,\\
\text{ and } \\
(j+1)x_1\left(2\frac{\p}{\p x}+A\right)v_{j+1}^+
+\left(2j\frac{\p}{\p x} +jA+Q\right)v_j^+=0, \;\ x< x_1,
\end{gathered} \label{transp}
\end{gather}
with the initial conditions $v_j^+(x_1,y),$ $j=1,2,...,$  given by 
\eqref{dvds2}.

Since $A=\ha\frac{1}{|h|}\frac{\p |h|}{\p x},$ we deduce from \eqref{transp} 
and \eqref{inicond} that
\begin{gather}
v_1^+(x,y)=\ha \frac{|h|^\frac{1}{4}(x_1,y)}{|h|^{\oq}(x,y)}\; 
 x_{1}^{-\frac{n}{2}-1}w(x_1,y), \;\ x\leq x_1. \label{solv1}
\end{gather}
But
\begin{gather}
\begin{gathered}
\frac{\p}{\p s} V^+(s,x,y)\sim \sum_{j=1}^\infty j e^s v_j^+(x,y)(e^s-x_1)_+^{j-1}
=\sum_{j=1}^\infty j (e^s-x_1+x_1) v_j^+(x,y)(e^s-x_1)_+^{j-1}
\sim  \\ x_1v^+_1(x,y)(e^s-x_1)_{+}^0
+\sum_{j=1}^{\infty} \left((j+1)x_1 v_{j+1}^+(x,y)+j v_{j}^+(x,y)\right)
(e^s-x_1)_{+}^j.
\end{gathered}\label{high}
\end{gather}
The highest  singularity of 
$\mcr_+(0, \chi_{x_1} w)=\left.\frac{\p v}{\p s}\right|_{\p X}$ at 
$s=\log x_1$ is 
the highest singularity of $\frac{\p}{\p s} V^+(s,0,y),$  which is
\begin{gather}
x_1v^+_1(0,y)(e^s-x_1)_{+}^0=\ha 
\frac{|h|^\frac{1}{4}(x_1,y)}{|h|^{\oq}(0,y)}\;  x_{1}^{-\frac{n}{2}}
w(x_1,y)(e^s-x_1)_+^0. \label{singp}
\end{gather}
This gives \eqref{key} and the Lemma is proved. 
\end{proof}

Now we can finish the proof of Proposition \ref{Schk}.
\begin{proof} We know from Lemma \ref{jump} that there exists $\eps>0$
such that for  any $\phi \in C_0^\infty(\mr)$
even, and $F\in \mcm^b,$  then
\begin{gather*}
\frac{|h_1|^{\oq}(x,y)}{|h_1|^{\oq}(0,y)}
\psi_1\left(\Delta_{g_1}-\frac{n^2}{4}\right)\mcr_{1,-}^{-1}F(x,y)=
\frac{|h_2|^{\oq}(x,y)}{|h_2|^{\oq}(0,y)}\psi_1
\left(\Delta_{g_2}-\frac{n^2}{4}\right)\mcr_{2,-}^{-1}F(x,y),
\;\ (x,y) \;\ (0,\eps)\times \p X,
\end{gather*}
where $\psi_1$ is determined by \eqref{phipsi}.  Now take a sequence
$\phi_m$ such that $\phi_m*F\rightarrow F$ in $L^2(\mr\times \p X)$
as $m\rightarrow \infty.$ Then, by the continuity of $\mcr_{j,-},\; $ $j=1,2,$
\begin{gather*}
\psi_{1,m}\left(\Delta_{g_j}-\frac{n^2}{4}\right)\mcr_{j,-}^{-1}F \rightarrow 
\mcr_{j,-}^{-1}F \text{ in } L^2(X_j) \text{ with respect to } g_j.
\end{gather*}
 Since $|h_1|(0,y)=|h_2|(0,y),$ \eqref{equality} holds
for the backward radiation field.

Let $F=\mcs^{-1}G=\mcr_{1,-}\mcr_{1,+}^{-1}G=\mcr_{2,-}\mcr_{2,+}^{-1}G$,
then
$$|h_1|^{\oq}(x,y)
\mcr_{1,+}^{-1}G=|h_2|^{\oq}(x,y)\mcr_{2,+}^{-1}G, \;\ \text{ in } \;\
 (0,\eps)\times \p X.$$ 
Since  $\mcs$ is unitary and $G$ is arbitrary, 
therefore \eqref{equality} holds for the forward radiation field.
This ends the proof of Proposition \ref{Schk}. 
\end{proof}

 Since $|h_1|=|h_2|$ the following  is an immediate  consequence  of
 Proposition \ref{Schk}.

\begin{prop}\label{detsk} Let $(X_1,g_1)$ and $(X_2,g_2)$ satisfy
the hypotheses of Theorem \ref{metdet} be such that $\Delta_{g_j},$
$j=1,2,$ have no eigenvalues. Then there exists $\eps>0$ such that
in  coordinates for which \eqref{norm12} holds, the 
Schwartz kernels of the radiation fields
$\mcr_{j,\pm}(s,y,z),$ for $(s,y)\in \mr\times M,$ 
$z\in (0,\eps) \times M,$  satisfy
\begin{gather}
\begin{gathered}
\mcr_{1,+}(s,y,z)=\mcr_{2,+}(s,y,z) \; \text{ and } \;
\mcr_{1,-}(s,y,z)=\mcr_{2,-}(s,y,z).
\end{gathered}\label{detesk1}
\end{gather}
\end{prop}

Next we want to show that $\Psi$ extends to a global diffeomorphism 
from $X_1$ to $X_2,$ as claimed in Theorem \ref{metdet}. We could follow the 
method of \cite{be,bk1,kaku}, see also \cite{kaku1,lassas}, and 
construct the diffeomorphism.  However it is easier to 
show that one can apply their result to construct $\Psi.$

\begin{prop}\label{bebk1}(\cite{be,bk1,kaku})
 Let $(Z,g)$ be a smooth compact Riemannian manifold with boundary
with  boundary $\p Z.$ Let
$\nu_{j},$ $\nu_j \leq \nu_{j+1},$
$j\in \mn,$ denote the Neumann eigenvalues of the operator 
$\Delta_{g},$ on $Z$ and let $\gamma_{j}$ 
denote the corresponding eigenfunctions.  
Then $(Z,g)$ is uniquely determined, modulo a diffeomorphism that is the 
identity at $\p Z,$  by the collection
\begin{gather}
\nu_{j} \;\ \text{and } \;\  \gamma_{j}|_{\p Z}, \;\ j\in \mn\setminus J,
\text{ where } \; J \text{ is a finite subset}. \label{nspd}
\end{gather}
\end{prop}
This result is proved in \cite{be,bk1} for $J=\emptyset$ and in 
\cite{kaku,kaku1} for $J\not=\emptyset.$

\begin{figure}[int1]
\epsfxsize= 4.0in
\centerline{\epsffile{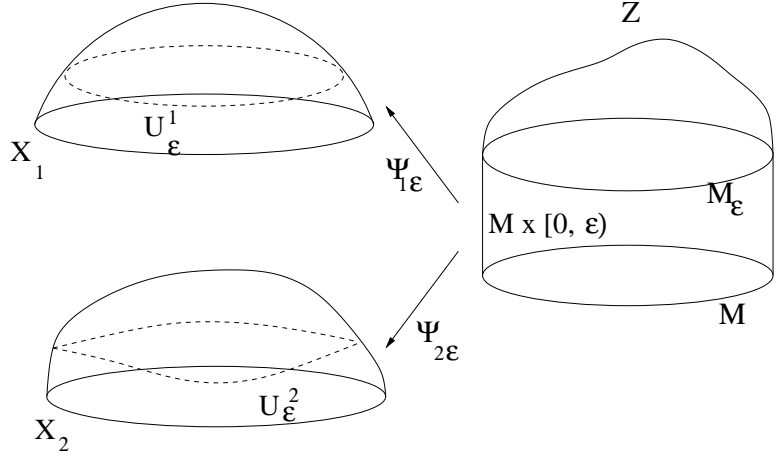}}
\caption{ The maps $\Psi_{j,\eps},$ $j=1,2$ and their extension to $Z.$}
\label{fig3n}
\end{figure}

So far we have shown that
there exists a diffeomorphism $\Psi$ of a neighborhood of $\p X$ such that
\begin{gather}
\Psi^*(g_1|_{U_{1,\eps}})= g_2|_{U_{2,\eps}} \;\ \text{ near }  \;\ M. 
\label{locdif}
\end{gather}
Observe that $X_{j,\eps}=X_j\setminus U_{j,\eps}$ are smooth compact manifolds
with boundary and their boundaries $\p X_{1,\eps},$ $\p X_{2,\eps}$
can be identified by the diffeomorphisms with  $M\times\{\eps\}=M_\eps,$ see
figure \ref{fig3n}. 

We think of $M_\eps$ as the boundary of some smooth compact manifold
$Z.$  We will prove that the fact that the
 two metrics have the same scattering matrix imply that
$\Delta_{g_i}$ in $X_{i,\eps},$ $i=1,2,$
have the same  Neumann spectral  data in $M_\eps,$ that is,
they have the same eigenvalues  and the
same traces of the eigenfunctions.
Then Proposition \ref{bebk1}, with $J=\emptyset,$
put together with Proposition \ref{determ} proves Theorem \ref{metdet}.
We  should also remark that method of proof of Proposition \ref{bebk1}
guarantees  that the resulting map  is $C^\infty.$

As in \cite{mafr}, we recall that
the graph of the Calder\'on projector of $\Delta_{g_j}-\la^2-\frac{n^2}{4}$ 
\; in $X_{j,\eps},$ $j=1,2,$ denoted by $C_{j,\la},$ is
the closed subspace of  $L^2(M_\eps)\times H^1( M_\eps)$ 
consisting of $(f,g) \in  L^2( M_\eps) \times H^{1}(M_\eps)$
such that there exists $u$   satisfying  
\begin{gather*}
\left(\Delta_{g_j}-\la^2-\frac{n^2}{4}\right)u=0 \;\ \text{ in } \;\ 
X_{j,\eps},  \; j=1,2,  \\
u|_{M_\eps}=f, \;\ \p_\nu u|_{M_\eps}=g.
\end{gather*}
Here  $\p_\nu  u$ denotes the normal derivative of $u$  at  $M_\eps.$
We will show that if $A_1(\la)=A_2(\la),$ 
$\la \in \mr\setminus{0},$ then $C_{1,\la}=C_{2,\la},$ 
$\la \in \mr\setminus{0}.$  But the Calder\'on projector depends continuously
on $\la$ and so does its graph. Thus $C_{1,\la}=C_{2,\la},$ $\la \in \mr.$

Since  $\Delta_{g_j},$ $j=1,2,$ has no point spectrum in $L^2(X_j),$
 the operators
$\Delta_{g_j}-\frac{n^2}{4},$ $j=1,2,$ are positive in $X_j.$  In particular
their  restriction to $X_{j,\eps}$ are also positive. 
Therefore
$\la^2+\frac{n^2}{4}$ is in the Neumann spectrum of 
$\Delta_{g_j}$ in $X_{j,\eps},$ if and only if
$C_{j,\la}$ contains a subspace of pairs of the form $(g,0),$ $g\not=0.$ 
Therefore, once we prove that
$C_{1,\la}=C_{2,\la},$ $\la \in \mr,$  then the eigenvalues of the
Neumann problem for
$\Delta_{g_1}$ in $X_{1,\eps},$ and the traces of the 
corresponding eigenfunctions
on $M_\eps,$ are equal to the corresponding ones
of $\Delta_{g_2}$ in $X_{2,\eps}.$  Then
Proposition \ref{bebk1}, with  $J=\emptyset,$ can be used to prove 
Theorem \ref{metdet}.  

To prove that
 $C_{1,\la}=C_{2,\la},$ 
$\la \in \mr\setminus{0},$
we apply the same argument used in the proof of Lemma 3.2, chapter 3.8 of
\cite{megs}, see also the proof of Lemma 2.1 of \cite{uh},
 to show that for any $\la\not=0,$ the set of functions given by
\begin{gather}
v_j(z,\la)=\int_{M} E_j^*\left(\frac{n}{2}+i\la\right)(z,y)\phi(y) , \;\ j=1,2,
\;\ \phi \in C^\infty(M), \label{einsf}
\end{gather}
where $E_j^*$ is the Eisenstein
Function, or Poisson
operator, which is the adjoint of the operator $E_j$ defined in
\eqref{defe}, is dense in the set of solutions of

\begin{gather}
\left(\Delta -\la^2-\frac{n^2}{4}\right)u=0 \;\ \text{ in } \;\ X_{j,\eps},\;\
j=1,2, \label{eins1}
\end{gather}
in the Sobolev space $H^k(X_{j,\eps}),$ for any $k\geq 2.$

The dual to $H^k(X_{j,\eps})$ can be identified with the space of
$f\in H^{-k}(X_j)$ supported in $X_{j,\eps}.$ If such  $f$ satisfies
$\lan f,v\ran=0$ for all $v$ given by \eqref{einsf}, then 
$E_j\left(\frac{n}{2}+i\la\right)f=0.$ This implies that
the function
$h=R\left(\frac{n}{2}+i\la\right)f$ vanishes to infinite order
at $\p X.$ Notice that $f$ is supported
in $X_{j,\eps},$ so $h$ is smooth near $\p X.$ Then unique 
continuation, see Theorem 14 of \cite{ma1}, implies
 that $h$ is supported in $X_{j,\eps}$ and
therefore $f=\left(\Delta_{g_j}-\la^2-\frac{n^2}{4}\right)h,$ with $h$ 
supported in $X_{j,\eps}.$ Then 
$\lan f,u\ran=0$ for any $u$ satisfying \eqref{eins1}.   

Since $E_j^*\left(\frac{n}{2}+i\la\right)(z,y)$ is the partial Fourier
transform in $s$ of $\mcr_{j,-}^{-1}(s,y,z),$ $j=1,2,$
Proposition \ref{detsk} implies that for all $v_j,$ $j=1,2,$ given by
\eqref{einsf}, which we know  is smooth in $X_j,$ satisfy
$v_1(z,\la)=v_2(z,\la),$ $z=(x,y)\in (0,\eps)\times M,$ 
$\la \in  \mr\setminus 0.$ Therefore
 their traces and normal derivatives at $M_\eps$
are equal, and the density  of this set implies  that
the same is true for solutions
of \eqref{eins1}. Thus $C_{1,\la}=C_{2,\la},$ 
$\la \in \mr\setminus{0}$ and this proves our claim.

\subsection{The general case} Now we remove the extra assumption on the
non-existence of eigenvalues.  Let $(X,g)$ be  an asymptotically
hyperbolic manifold.
 The only poles of the resolvent 
$R(\frac{n}{2}+i\la)=(\Delta_g-\frac{n^2}{4}-\la^2)^{-1}$
in $\{\Im \la<0\}$ correspond to the finitely many eigenvalues of
$\Delta_g.$
Proposition 3.6 of \cite{grazw} states
that  if $\la_0>0$ is such that $\frac{n^2}{4}-\la_0^2$ 
is an eigenvalue of
$\Delta_g$ then the scattering matrix has a pole at $-i\la_0$ and its
residue is  given  by
\begin{gather}
\begin{gathered}
\Res_{-i\la_0} A(\la)= \left\{\begin{array}{c} \Pi_{\la_0}, \;\ \text{ if } \;\
\la_0 \not\in \mn/2, \\
\Pi_{\la_0}-p_l, \;\ \text{ if } \;\ \la_0=\frac{l}{2}, \; l\in \mn, \end{array}\right.
\end{gathered}\label{grazw}
\end{gather}
where $p_l$ is a differential operator whose coefficients depend on 
the tensor  $h,$ defined in \eqref{norm12}, and its derivatives  at $\p X,$  
and 
the Schwartz kernel of
$\Pi_{\la_0}$ is
\begin{gather}
\begin{gathered}
K(\Pi_{\la_0})(y,y')=2\la_0\sum_{j=1}^N \phi_j^0\otimes\phi_j^0(y,y'), \;\
\phi_j^0(y)=x^{-\frac{n}{2}-\la_0}\phi_j(x,y)|_{x=0}.
\end{gathered} \label{grazw1}
\end{gather}
Here  $N$ is the multiplicity
of the eigenvalue $\frac{n^2}{4}-\la_0^2$ and $\phi_j,$ $1\leq j \leq N,$ are 
the corresponding orthonormalized eigenfunctions.

If two asymptotically hyperbolic manifolds $(X_1,g_1)$ and $(X_2,g_2)$
have the same scattering matrix $A_1(\la)=A_2(\la)$ for all 
$\la\in \mr\setminus 0,$ we know 
from \cite{jsb} that in coordinates
where \eqref{norm12} is satisfied, all derivatives of $h_1$ and $h_2$
agree at $x=0.$ 
Therefore  the operators $p_l$ in \eqref{grazw} are the same.
Thus, \eqref{grazw} and \eqref{grazw1},
and the meromorphic continuation of the scattering matrix, show that
$\Delta_{g_1}$ and $\Delta_{g_2}$ have the same eigenvalues, with
the same multiplicity. Moreover, \eqref{grazw1} implies that
if $\phi_j,$ and  $\psi_j,$ $ 1 \leq j  \leq N,$ are
orthonormal sets of eigenfunctions of $\Delta_{g_1}$ and $\Delta_{g_2},$ 
respectively,
corresponding to the eigenvalue $\frac{n^2}{4}-\la_0^2,$ 
then there exists  a constant orthogonal $N\times N$ matrix
$A$ such that $\Phi^0=A\Psi^0,$ where
 $(\Phi^0)^T=(\phi_1^0,\phi_2^0,...,\phi_N^0),$ and 
$(\Psi^0)^T=(\psi_1^0,\psi_2^0,...,\psi_N^0).$ 
So by redefining one set of eigenfunctions  from let  us say,
$\Psi$ to $A\Psi,$ where $\Psi^T=(\psi_1,\psi_2,...,\psi_N),$
 we may assume  that 
\begin{gather}
\phi_j^0(y)=\psi_j^0(y),  \;\  j=1,2,...,N. \label{eqeig}
\end{gather}
We remark that this does  not
change the orthonormality of  the eigenfunctions in $X_2$
 because $A$  is orthogonal.
Let  us denote
\begin{gather}
\begin{gathered}
\mu_j=\frac{n^2}{4}-\la_j^2, \;\ \phi_j, \;\ 1\leq j \leq L, \;\
\text{ all the eigenvalues and eigenfunctions of } \; \Delta_{g_1}, \\
\mu_j=\frac{n^2}{4}-\la_j^2, \;\ \psi_j, \;\ 1\leq j \leq L, \;\
\text{ all the eigenvalues and eigenfunctions of } \; \Delta_{g_2},
\end{gathered}\label{eigenf}
\end{gather}
with $\mu_1\leq \mu_2 \leq... \leq  \mu_L,$ and the eigenfunctions
$\phi_j$ and $\psi_j$ chosen to satisfy \eqref{eqeig}.

When there are eigenfunctions, 
Lemma \ref{proj} is no longer valid.  To present the correct statement we
begin with the following lemma, which was suggested and proved by
one of the referees,
\begin{lemma}\label{cjx0} Let $V$  be  a finite  dimensional  subspace  of
a Hilbert space $H$ and let $Q$ be  the  orthogonal projector  onto $V.$ 
 Let $P_t,$ $t\in  [0,a],$ be a strongly continuous
family of projections with $P_0=\Id$ and $P_t P_{t'}=P_t$  if $t\leq t'.$
Then there  exists $\eps  \in (0,a),$ and a unique continuous
family of bounded operators $T(t): H \longrightarrow V,$ $t\in [0,\eps],$
 such  that  $QP_t(\Id-T(t))=0.$ Moreover,
if $f_m\in H$ is a bounded sequence,  then there exists a subsequence
$f_{m'}$ such that $T(t)f_{m'}$ converges uniformly in $[0,\eps].$
\end{lemma}
\begin{proof}  As $Q(P_t)|_{V}$
is continuous and $Q(P_0)|_V=\Id|_V,$ there exists $\eps\in (0,a)$ such that
$QP_t: V \longrightarrow V$ is invertible on $V$ for $t\in  [0,\eps].$ If
$R(t)$ is its inverse, then $R(t)$ is continuous for $t\in [0,\eps],$ and
\begin{gather*}
QP_t(\Id-R(t)QP_t)=QP_t-(QP_t)R(t)QP_t=0.
\end{gather*}
So we take $T(t)=R(t)QP_t.$ Then $T(t)$ is  bounded and continuous in $t.$
If $\wtt(t)$ is another solution,
$QP_t(T(t)-\wtt(t))=0,$ so  multiplying by $R(t)$ on  the left,  we get
$T(t)-\wtt(t)=0.$

If $f_m$ is a bounded sequence of functions on $H,$ then $P_t f_m$
is bounded. From the strong continuity of $P_t,$ $P_t f_m$
is equicontinuous  and hence it has a convergent subsequence
$P_tf_{m'}.$ Thus $T(t)f_{m'}=R(t)QP_t f_{m'}$ converges. This ends the
proof of the lemma. 
\end{proof} 

The application of this lemma that we have in mind is:
\begin{cor}\label{cjx} Let $(X,g)$ be an asymptotically hyperbolic manifold
 and let $x$ be a defining function of $\p X$ for which \eqref{hmet1} holds.
Let $V=L^2_{\pp}(X),$  be the space
spanned by the set of  orthonormal 
eigenfunctions  of $\Delta_g,$ $\zeta_j \in L^2(X),$  $1 \leq j \leq N.$
Let $\chi_t$ denote  the characteristic function of the set
$\{x \geq t\}$ and let $P_t f=\chi_t f.$ 
Then there exists $\eps>0,$ and a unique continuous family of bounded 
operators $T(t):L^2(X) \longrightarrow  V,$ $t\in [0,\eps],$ such that
\begin{gather}
\lan \chi_{t}( f-T(t)f), \zeta_k\ran = 0, \;\ k=1,2,...,N. \label{cj}
\end{gather}
Moreover, any bounded sequence $f_m \in L^2(X),$ $m=1,2,...,$
has a subsequence $f_{m'}$ such that 
$T(t)(f_{m'})$ converges uniformly in $[0,\eps],$ $j=1,2,...,N.$
In this case the operators $R(t)$ and $QP_t$ are given by
\begin{gather}
R(t)=(R(t)_{ij}), \;\ R(t)_{ij}
=\lan \chi_t \zeta_i, \zeta_j \ran, \;\
(QP_tf)^T =\left(\lan \chi_t f, \zeta_1\ran,\lan \chi_t f, \zeta_2\ran,...,
\lan \chi_t f, \zeta_N\ran\right).\label{cjx1}
\end{gather}
\end{cor}

When there exist eigenfunctions, Lemma \ref{proj} has to be replaced by

\begin{lemma}\label{proj1} Let $(X,g)$ be an asymptotically hyperbolic
manifold and let $\{\zeta_j: 1\leq j \leq N\}$ denote the $L^2(X)$ 
eigenfunctions  of $\Delta_g.$
 Let  $\mcp_{x_1}^b$ denote the orthogonal projector
\begin{gather*}
\mcp_{x_1}^b : \mcm^b \longrightarrow \mcm^{b}( x_1).
\end{gather*}
Let $\chi_{t}$ be the characteristic function of
 the set $\{x\geq t\}.$
Then there exists $\eps>0$ such that for every $x_1\in (0,\eps)$ and
every $f\in L_{\ac}^2(X)$ 
\begin{gather}
\begin{gathered}
\mcp_{x_1}^b\mcr_{-}(0,f)=\mcr_{-}(0,f_{x_1}), \;\ \text{ where } \\
f_{x_1}(x,y)=\chi_{x_1}(x,y)\left( \Id - T(x_1)\right)f,
\end{gathered}\label{f1f2}
\end{gather}
and $T(x_1)$ is the family of operators given by Corollary \ref{cjx}.
\end{lemma}
\begin{proof}  We follow the proof of Lemma \ref{proj}.
We know from Lemma \ref{charact} that there exists $f_{x_1}\in L^2_{\ac}(X)$
supported in $\{x \geq x_1\}$ such that
$\mcp_{x_1}^b\mcr_{-}(0,f)=\mcr_{-}(0,f_{x_1})$ and that
\begin{gather*}
\lan f_{x_1}-f, g \ran =0 \;\ \text{ for all }
g \in C_0^\infty(\intx)\cap L^2_{\ac}(X) \;\ \text{ supported in}
\{x_1 \leq x \}.
\end{gather*}
Therefore there exist constants $\alpha_j(x_1),$ $j=1,2,...,L,$ such that
$$ f_{x_1}= \chi_{x_1}\left( f + \sum_{j=1}^N \alpha_j(x_1) \zeta_j\right).$$
Since $f_{x_1}\in L^2_{\ac}(X_1),$ it follows that
 $\lan f_{x_1}, \zeta_k\ran=0,$ $1 \leq k \leq N.$  Then, by uniqueness,
$\sum_{j=1}^N\alpha_{j}(x_1)\zeta_j=-T(x_1)(f),$ is the family
of operators given by Corollary \ref{cjx}.
\end{proof}

If there are eigenvalues, then Lemma \ref{jump} is replaced by

\begin{lemma}\label{jump1}   Let $(X,g)$ be an
asymptotically hyperbolic manifold and let $x$ be a defining function of
$\p X$ such that \eqref{hmet1} holds. Let $\mcr_{\pm}$ be the  radiation fields
with respect to $x.$ Then
there exists $\eps>0$ such that for all
 $F\in \mcm^b,$  $\phi \in C_0^\infty(\mr)$ even, and 
$\psi_1$  defined by \eqref{phipsi}, and any $x_1\in(0,\eps/4),$ 
\begin{gather}
\begin{gathered}
\mcr_{+}\mcr_{-}^{-1}\left(\mcp_{x_1}^b G\right)(s,y)
= \\ \ha x_1^{-\frac{n}{2}}\left(\Id-T(x_1)\right)w(x_1,y)
\frac{|h_1|^{\oq}(x_1,y)}{|h_1|^{\oq}(0,y)}
(s-\log x_1)_{+}^0 + \text{ smoother terms},\\
\text{ if } \; s< \log x_1+\log 4, \; \text{ where } \;\ 
G=\phi*F, \;\ \text{ and } w_1=\mcr_{-}^{-1}G=
\psi_1\left(\Delta_{g}-\frac{n^2}{4}\right)\mcr_{-}^{-1}F.
\end{gathered} \label{key1}
\end{gather}
Here $T(x_1)$ is the operator given by Corollary \ref{cjx}.
Notice that by \eqref{7.2new}, $G=\phi*F\in \mcm^b,$ and that
 the left hand side of  \eqref{key1} is  determined  by the scattering
matrix $A(\la).$
\end{lemma}
\begin{proof} We follow step by step
the proof of Lemma \ref{jump}. We know that
the eigenfunctions are smooth in the interior of the manifold. So the only
thing that changes is the singularity of the initial data in
equation \eqref{dvds2}. The initial data  of $v_1^+,$ which was given by
\eqref{inicond}, must, according to \eqref{f1f2},
 be replaced by
$$v^+(x_1,y)= \ha x_1^{-\frac{n}{2}-1}\left(\Id - T(x_1)\right)w_1(x_1,y).$$
 \end{proof}

As in the proof of Proposition \ref{Schk}, we take a sequence $\phi_m$ such 
that 
$\phi_m*F\rightarrow F$ in  $L^2(\mr\times M),$ in Lemma \ref{jump1}, and
we obtain
\begin{prop}\label{Schk1}
 Let $(X_j,g_j),$ $j=1,2,$ be asymptotically hyperbolic manifolds
satisfying the hypotheses of Theorem \ref{metdet}.
 Let $\mcr_{j,-},$ $j=1,2,$ denote
the corresponding backward radiation fields defined in coordinates in which
\eqref{norm12} holds.  Let $T_j(t),$ $j=1,2,$ be the operators
given by Corollary \ref{cjx}. Then  there exists $\eps>0$ such that
for any $F\in \mcm^b,$
\begin{gather}
\begin{gathered}
|h_1|^{\oq}(x,y)
\left( \Id- T_1(x)\right)\mcr_{1,-}^{-1}F(x,y)
=|h_2|^{\oq}(x,y)\left( \Id- T_2(x)\right) \mcr_{2,-}^{-1}F(x,y) \;\
\text{ in } (0,\eps)\times M.
\end{gathered}\label{equality1}
\end{gather}
\end{prop}

We  deduce from \eqref{equality1} and \eqref{f-1} 
that for every  $F\in \mcm^b,$
\begin{gather}
\begin{gathered}
|h_1|^\oq (\Id -T_1)\mcr_{1,-}^{-1} F(x,y)=
|h_2|^\oq (\Id -T_2)\mcr_{2,-}^{-1} F(x,y), 
\;\ (x,y) \in (0,\eps) \times M,
\;\ \text{ and } \\
|h_1|^\oq (\Id -T_1)\left(\Delta_{g_1}-\frac{n^2}{4}\right)\mcr_{1,-}^{-1}
F(x,y)=
|h_2|^\oq (\Id -T_2)\left(\Delta_{g_2}-\frac{n^2}{4}\right)\mcr_{2,-}^{-1}
F(x,y), \;\ (x,y) \in (0,\eps) \times M.
\end{gathered}\label{equality2}
\end{gather}

We will use this to  conclude that the 
tensors $h_1$ and $h_2$ are equal, and  thus prove Proposition \ref{determ}.
 This is based on the following lemma, which is a direct consequence of
Corollary \ref{cjx}.

\begin{lemma}\label{tco}  For $\del>0,$ let $K_\del$
be the characteristic function of $\{x\leq \del\}.$ Then there exists $\eps>0$
such that for every $\del\in  (0,\eps),$
the operators $T_j,$ $j=1,2,$ defined by Corollary \ref{cjx} satisfy
\begin{gather*}
K_\del T_j : L^2(X_j) \longrightarrow
L^2([0,\del] \times M), \;\ j=1,2, 
\end{gather*}
where $L^2$ is defined with
respect to the metric $g_j,$ and moreover are compact.
\end{lemma}

We  can then prove Proposition \ref{determ}.  
\begin{proof}We just need to observe  that
it follows from \eqref{equality2} and Lemma \ref{tco} that for $\eps>0$ small,
\begin{gather}
\begin{gathered}
K_\eps |h_1|^\oq\mcr_{1,-}^{-1} F(x,y)= 
K_\eps|h_2|^\oq \mcr_{2,-}^{-1} F(x,y)
+K_\eps R F(x,y), \;\ R \; \text{ compact,  and } \\
K_\eps |h_1|^\oq\left(\Delta_{g_1}-\frac{n^2}{4}\right)\mcr_{1,-}^{-1}
F(x,y)=
K_\eps |h_2|^\oq \left(\Delta_{g_2}-\frac{n^2}{4}\right)\mcr_{2,-}^{-1}
F(x,y)+K_\eps SF(x,y), \;\ S  \text{ compact. }
\end{gathered}\label{equality3}
\end{gather}
Let $\delta \in (0,\eps)$ and set $\mcr_{1,-}^{-1}F=f.$ We find that for all  
$f \in L^2_{\ac}(X_1),$   with respect to $g_1,$ 
\begin{gather}
\begin{gathered}
K_\del |h_2|^{\oq} \left(\Delta_{g_2}-\frac{n^2}{4}\right) 
\frac{|h_1|^{\oq}}{|h_2|^{\oq}} f=
K_\del 
|h_1|^{\oq} \left(\Delta_{g_1}-\frac{n^2}{4}\right)f + K_\del \widetilde{S}f,
\;\  \widetilde{S} \;\ \text{ compact.}
\end{gathered}\label{equality4}
\end{gather} 
Let $Q$ be the differential operator
\begin{gather*}
Q f=\left(|h_2|^{\oq} \left(\Delta_{g_2}-\frac{n^2}{4}\right) 
\frac{|h_1|^{\oq}}{|h_2|^{\oq}}-
|h_1|^{\oq} \left(\Delta_{g_1}-\frac{n^2}{4}\right)\right)f.
\end{gather*}
Since any 
$f \in L^2(X_1),$ with respect to $g_1,$ can be written as
$f= \mcp_{\ac} f + \sum_{j=1}^N \lan f, \phi_j\ran \phi_j,$
it follows from \eqref{equality4} that
\begin{gather*}
K_\del Qf = K_\del Q\mcp_{\ac} f + 
\sum_{j=1}^N  \lan f, \phi_j\ran K_\del Q \phi_j=
K_\del \widetilde{S} \mcp_{\ac} f +\sum_{j=1}^N  \lan f, \phi_j\ran 
K_\del  Q \phi_j, \;\ f \in L^2(X_1).
\end{gather*}
The  first  term  of the sum is compact because $\widetilde{S}$ is, and
the second
term is of finite rank. Therefore $K_\del Q: L^2(X_1) \longrightarrow
H^{-2}([0,\del] \times M)$ is compact. But $Q$ is also a second order 
differential 
operator. This implies that $K_\del Q=0.$ See  for example exercise 6.2  
on page 52 of \cite{taylor}.   Therefore  the tensors $h_1$ and $h_2$
are equal in  $[0,\del]  \times M,$ and this proves Proposition \ref{determ}.
\end{proof}

Now we need to show that the diffeomorphism  can be extended to  the whole
manifold. We will use the same method as in the case of no eigenvalues.
 We  have shown  that  $\Delta_{g_1}=\Delta_{g_2}$ in 
coordinates \eqref{norm12} in $(0,\eps)\times M.$ Using
 \eqref{eqeig} and the equation satisfied by the eigenfunctions, it 
can be shown that \eqref{eqeig} is satisfied to  infinite order, that is
all derivatives of the ``rescaled eigenfunctions''  agree at $\{x=0\}.$
Therefore unique continuation  for this type of operators, see
 Theorem 14 of \cite{ma1}, shows that there exists $\eps>0,$ such that the 
eigenfunctions of $\Delta_{g_1}$ and those of  $\Delta_{g_2}$ are equal
in $[0,\eps)\times M.$ That is
\begin{gather}
\phi_j(x,y)=\psi_j(x,y), \;\ 1 \leq  j \leq L, \;\ (x,y)\in [0,\eps)\times M.
 \label{eigfn}
\end{gather}

We deduce from the first equation in \eqref{equality2} that the Schwartz 
kernels  $\mcr_{j,-}^{-1}(s,y,z),$ $(s,y) \in \mr\times M$ and $z=(x,y')
\in (0,\eps) \times M$ satisfy
\begin{gather}
(I-T_1)\mcr_{1,-}^{-1}(s,y,z) = (I-T_2)\mcr_{2,-}^{-1}(s,y,z). \label{alm0}
\end{gather} 

Since $T_j,$ $j=1,2,$ is a linear operator, we  may
take Fourier transform  in $s$ of \eqref{alm0}
and deduce that the Schwartz kernels  $E_j^*(\frac{n}{2}+i\la)(z,y),$
$j=1,2,$ $(\la,y)\in \mr\setminus 0 \times M,$  and 
$z=(x,y') \in (0,\eps) \times M$ satisfy 
\begin{gather}
(I-T_1)E_1^*(\frac{n}{2}+i\la)(z,y)=(I-T_2)E_2^*(\frac{n}{2}+i\la)(z,y).
 \label{alm1}
\end{gather}

We will use this to prove

\begin{prop}\label{difexten2} There exists $\eps>0$ such that
for every 
$\la \in \mr\setminus 0,$ $y\in M$ and
$z=(x,y')\in (0,\eps)\times M,$
\begin{gather}
E_1^*(\frac{n}{2}+i\la)(z,y)=E_2^*(\frac{n}{2}+i\la)(z,y).\label{cauchy}
\end{gather}
\end{prop}
\begin{proof}
Let $\phi \in C^\infty(M)$ and let $v_j(z,\la),$ $j=1,2,$ be the functions 
given by \eqref{einsf}.
We will show that there exists $\eps>0$ such that for
 $\delta \in (0,\eps),$ there exists
 $\La=\La(\del)>0$ such that
\begin{gather}
v_1(\del,y',\la)=v_2(\del,y',\la), \;\ \forall \;\ y'\in M, \;\
\text{ and } |\la|>\La.
\label{cauchyd}
\end{gather}
Since by the analytic continuation of 
$E_j(\frac{n}{2}+i\la),$ $\la\in \mr\setminus 0,$ 
$v_j(z,\la),$ $j=1,2,$ is real analytic in $\la\in \mr\setminus 0,$
for each $z,$ 
it follows  that \eqref{cauchyd} holds for every $\la\in \mr\setminus 0.$
Since $\del$ is arbitrary, \eqref{cauchy} follows.

Equation \eqref{alm1} implies that
\begin{gather}
\begin{gathered}
v_1(z,\la)-v_2(z,\la)=T_1v_1(z,\la)-T_2v_2(z,\la), \;\ 
z=(x,y') \in (0,\eps)\times M,
\;\ \la \in \mr\setminus 0.
\end{gathered} \label{alm1v1v2}
\end{gather}

Let us denote
\begin{gather*}
\Phi(x,y)=(\phi_1(x,y),...,\phi_L(x,y)), \;\
\Psi(x,y)=(\psi_1(x,y),...,\psi_L(x,y)), \\
T_1v_1(z,\la)=\sum_{j=1}^L C^1_j(x,\la) \phi_j(x,y), \;\
T_2v_2(z,\la)=\sum_{j=1}^L C^2_j(x,\la) \psi_j(x,y), \;\
\end{gather*}
where, by  equation \ref{cjx1},
\begin{gather*}
C^{j}(x,\la)\defi C^{j}(x,v_j(z,\la))=
\left[R^j(x)\right]^{-1} F^{j}(x,\la),\;\ j=1,2, \text{ where }\\
F^1(x,\la)^T\defi F^1(x,v_1(z,\la))^T=
(\lan\chi_x  v_1(z,\la), \phi_1\ran, ..., \lan \chi_x v_1(z,\la), \phi_L\ran),
 \\
 F^2(x,\la)^T\defi F^2(x,v_2(z,\la))^T=
(\lan\chi_x  v_2(z,\la), \psi_1\ran, ..., \lan \chi_x v_2(z,\la), \psi_L\ran) 
\\
R^1(x)=(R_{ij}^1(x)), \;\ R_{ij}^1(x)=\lan \chi_x \phi_j,\phi_j \ran, \;\
R^2(x)=(R_{ij}^2(x)), \;\ R_{ij}^2(x)=\lan \chi_x \psi_j,\psi_j \ran, \;\
\end{gather*}

We know that $\lan \phi_i, \phi_j\ran=\lan \psi_i, \psi_j\ran =\del_{ij},$
we also  know  that if $x<\eps,$ then  $(1-\chi_x)\phi_i=(1-\chi_x)\psi_i.$
But,
\begin{gather*}
\lan \chi_x \phi_i,\phi_j\ran=
\lan \phi_i,\phi_j\ran-\lan (1-\chi_x) \phi_i,\phi_j\ran=
\lan \phi_i,\phi_j\ran-\lan (1-\chi_x) \phi_i,(1-\chi_x)\phi_j\ran=\\
\lan \psi_i,\psi_j\ran-\lan (1-\chi_x) \psi_i,(1-\chi_x)\psi_j\ran=\lan \psi_i,\psi_j\ran-\lan (1-\chi_x) \psi_i,\psi_j\ran=
\lan \chi_x \psi_i,\psi_j\ran.
\end{gather*}
Hence
\begin{gather}
R_{ij}^1(x)=\lan \chi_x \phi_i,\phi_j\ran=R_{ij}^2(x)=
\lan \chi_x \psi_i,\psi_j\ran.\label{a1=a2}
\end{gather}
  We denote 
$R(x)=R^1(x)=R^2(x).$ We can pick $\eps$ small so that $R^{-1}(x)$ is 
uniformly bounded for $x\in [0,\eps].$ 

Recall that $\mu_j=\frac{n^2}{4}-\la_j^2$ is the eigenvalue corresponding to $\phi_j,$ as
 defined in \eqref{eigenf}.  For $\del<\eps,$ the divergence theorem
gives
\begin{gather*}
\lan \chi_\del v_1(z,\la), \phi_j(z)\ran = 
\frac{1}{\mu_j} \lan \chi_\del v_1(z,\la), \Delta_{g_1}\phi_j(z)\ran= \\
\frac{1}{\mu_j}(\frac{n^2}{4}+\la^2)\lan \chi_\del v_1(z,\la), \phi_j(z)\ran
+ \frac{1}{\mu_j}
\int_{M} \left[v_1(\del,y,\la) \frac{\p \phi_j}{\p x}(\del,y)-
\frac{\p v_1}{\p x}(\del,y,\la) \phi_j(\del,y)\right] 
\frac{\sqrt{h}(\del,y)}{\del^{n-1}}\; dy.
\end{gather*} 
Doing the  same computation  for $\lan \chi_\del v_2, \psi_j\ran,$ 
using that $\mu_j=\frac{n^2}{4}-\la_j^2$ and \eqref{eigfn}, we obtain
\begin{gather}
\begin{gathered}
\lan \chi_\del v_1(z,\la), \phi_j(z)\ran = 
-\frac{1}{\la_j^2+\la^2}
\int_{M} \left[v_1(\del,y,\la) \frac{\p \phi_j}{\p x}(\del,y)-
\frac{\p v_1}{\p x}(\del,y,\la) \phi_j(\del,y)\right]\frac{\sqrt{h}(\del,y)}{\del^{n-1}}\; dy,\\
\lan \chi_\del v_2(z,\la), \psi_j(z)\ran = 
-\frac{1}{\la_j^2+\la^2}
\int_{M} \left[v_2(\del,y,\la) \frac{\p \phi_j}{\p x}(\del,y)-
\frac{\p v_2}{\p x}(\del,y,\la) \phi_j(\del,y)\right] 
\frac{\sqrt{h}(\del,y)}{\del^{n-1}}\; dy.
\end{gathered}\label{inner}
\end{gather}

We conclude
from \eqref{inner} that
there exists $K_0(\del)>0,$ such that
\begin{gather}
|F^1(\del,\la)-F^2(\del,\la)|<\frac{K_0(\del)}{\la^2} 
\left[ \sup_{M} |v_1(\del,y,\la)-v_2(\del,y,\la)| +
\sup_{M} \left| \frac{\p  v_1}{\p x}(\del,y,\la)-
\frac{\p  v_2}{\p x}(\del,y,\la)\right|\right]. \label{8.50new}
\end{gather}
This and \eqref{a1=a2} imply that there exists $K_1(\del)>0,$ such that 
\begin{gather}
\begin{gathered}
\sup_{M}\left|T_1v_1(\del,y,\la)-T_2 v_2(\del,y,\la)\right|\leq \\
\frac{K_1(\del)}{\la^2}
\left[ \sup_{M} |v_1(\del,y,\la)-v_2(\del,y,\la)| +
\sup_{M} \left| \frac{\p  v_1}{\p x}(\del,y,\la)-
\frac{\p  v_2}{\p x}(\del,y,\la)\right|\right].
\end{gathered} \label{almo}
\end{gather}

Setting $x=\delta$ in  \eqref{alm1v1v2} and using
\eqref{almo} we get
\begin{gather}
\sup_{M}|v_1(\del,y,\la)-v_2(\del,y,\la)| \leq 
\frac{K_1(\del)}{\la^2}
\left[ \sup_{M}|v_1(\del,y,\la)-v_2(\del,y,\la)|+
\sup_{M}|\frac{\p v_1}{\p x}(\del,y,\la)-
\frac{\p v_2}{\p x}(\del,y,\la)|\right].\label{almo2}
\end{gather}

To  estimate the terms involving $\frac{\p v_j}{\p x}(\del,y,\la),$
we differentiate \eqref{alm1v1v2} in $x.$  We get that
\begin{gather}
\frac{\p v_1}{\p x}(x,y,\la)-\frac{\p v_2}{\p x}(x,y,\la)=
\frac{\p T_1 v_1}{\p x}(x,y,\la)-
\frac{\p T_2 v_2}{\p x}(x,y,\la). \label{differ}
\end{gather}
Since $\Phi=\Psi$  for $x<\eps,$ we have
$T_j v_j(x,y,\la)=R^{-1}(x) F^j(x,\la)\cdot \Phi(x,y),$ and thus
\begin{gather}
\begin{gathered}
\frac{\p }{\p x}T_j v_j(x,y,\la)= \\
\left[ \frac{d}{d x} R^{-1}(x)\right] F^j(x,\la)\cdot \Phi(x,y)
+ R^{-1}(x) \frac{d }{d x}F^j(x,\la)\cdot \Phi(x,y)+
R^{-1}(x)F^j(x,\la) \cdot \frac{\p \Phi}{\p x}(x,y).
\end{gathered}\label{differ1}
\end{gather}
The first and the third terms are much like $T_jv_j(z,\la),$ and because of
\eqref{8.50new}, satisfy
an estimate like \eqref{almo}. The second term has to be considered
separately.  Notice that
\begin{gather*}
\left(\frac{d}{d x} F^j(x,\la)\right)^T=
\left(\frac{d }{d x}\lan\chi_x v_j,\phi_1\ran,
...\frac{d }{d x}\lan\chi_x v_j,\phi_L\ran\right),
\end{gather*}
and
\begin{gather*}
\frac{d}{d x} \lan\chi_x v_j,\phi_k\ran\restr_{x=\del}=-
\frac{1}{\del^{n+1}}
\int_{M} v_j(\del,y,\la) \phi_k(\del,y,\la) \sqrt{h}(\del,y,\la) dy.
\end{gather*}

Then
\begin{gather}
\frac{d}{d x} \left(\lan\chi_x v_1,\phi_k\ran
- \lan\chi_x v_2,\phi_k\ran\right)\restr_{x=\del}=-
\frac{1}{\del^{n+1}}
\int_{M}\left( v_1(\del,y,\la)-v_2(\del,y,\la)\right)
 \phi_k(\del,y,\la) \sqrt{h}(\del,y,\la) dy. \label{middle}
\end{gather}

Then there exists $K_2(\del)$  such that
\begin{gather}
\begin{gathered}
\left|\frac{d }{d x}F^1(\del,\la)-\frac{d }{d x}F^2(\del,\la)\right|
\leq  \frac{K_2(\del)}{\la^2} \sup_{M} |v_1(\del,y,\la)-v_2(\del,y,\la)|.
\end{gathered}\label{differ2}
\end{gather}
So from \eqref{8.50new},
 \eqref{differ}, \eqref{differ1} and \eqref{differ2} we find that 
there exists $K_3(\del)$ such that
\begin{gather}
\sup_{M}|\frac{\p v_1}{\p x}(\del,y,\la)-
\frac{\p v_2}{\p x}(\del,y,\la)| \leq 
\frac{K_3(\del)}{\la^2}
\left[ \sup_{M}|v_1(\del,y,\la)-v_2(\del,y,\la)|+
\sup_{M}|\frac{\p v_1}{\p x}(\del,y,\la)-
\frac{\p v_2}{\p x}(\del,y,\la)|\right].\label{almo3}
\end{gather}

So we conclude from \eqref{almo2} and \eqref{almo3} that
\begin{gather*}
\sup_{M}|v_1(\del,y,\la)- v_2(\del,y,\la)| +
\sup_{M}|\frac{\p v_1}{\p x}(\del,y,\la)-
\frac{\p v_2}{\p x}(\del,y,\la)| 
\leq  \\
\frac{K_1(\del)+K_3(\del)}{\la^2}
\left[ \sup_{M}|v_1(\del,y,\la)-v_2(\del,y,\la)|+
\sup_{M}|\frac{\p v_1}{\p x}(\del,y,\la)-
\frac{\p v_2}{\p x}(\del,y,\la)|\right]
\end{gather*}

Taking $\La^2>K_1(\del)+K_3(\del),$ equation \eqref{cauchyd} follows.
This ends the proof of the Proposition.
\end{proof}

As in the case of no eigenvalues, this can be used to show that the
graphs of the Calder\'on projectors satisfy $C_{1,\la}=C_{2,\la},$ 
$\la \in \mr\setminus 0$ and by continuity of $C_{j,\la},$ if
$\la \in \mr.$ 
However, in this case, the operators
$\Delta_{g_j}-\frac{n^2}{4},$  $j=1,2,$ are not positive, and
this only determines the
Neumann eigenvalues, and the traces of  the corresponding eigenfunctions,
of $\Delta_{g_j}$
for eigenvalues $\nu_j\geq \frac{n^2}{4}.$  But since the set
$J=\{j: \nu_j\leq \frac{n^2}{4}\}$ is finite, we may again apply
Proposition \ref{bebk1}, this time with $J\not=\emptyset.$

Observe that for $\del<\eps$   fixed,
\eqref{cauchyd}, and  the density of $v_j(\la,z),$ $j=1,2,$
imply  that $C_{1,\la}=C_{2,\la},$ $\la^2>\Lambda^2$ in $X_{j,\del}.$
This gives the
spectral  data with eigenvalues $\nu_j>\Lambda$  and we could apply
Proposition \ref{bebk1} already at this  stage. 

Then Proposition \ref{bebk1} combined 
Proposition \ref{determ}  proves Theorem \ref{metdet}. Again, we remark
that  the method of proof of Proposition  \ref{bebk1} guarantees  that
the map  is $C^\infty.$

\section{Acknowledgments}

This paper is dedicated to the memory of F.G. Friedlander,
who started the study of radiation fields and whose 
ideas were often used here. I also had the pleasure of having several 
discussions about the subject with him.

I thank R. Melrose, who was a student of
Friedlander, for bringing the paper \cite{fried2} to my attention and for
several instructive conversations.
I also thank C. Guillarmou, Y. Kurylev, 
R. Mazzeo , P. Stefanov, S. Tang and J. Wunsch for useful conversations and 
for comments on the paper.

I am very grateful to the referees for carefully reading the manuscript.
They found several mistakes in previous versions, and  made numerous 
suggestions  to improve the exposition. 

Part of this work was carried out while I was visiting the
Universidade Federal de Pernambuco, Brazil. I thank them, and especially
Prof. Paulo Santiago, for their hospitality.

\end{document}